# CONTINUUM MANY DIFFERENT THINGS: LOCALISATION, ANTI-LOCALISATION AND YORIOKA IDEALS

MIGUEL A. CARDONA, LUKAS DANIEL KLAUSNER, AND DIEGO A. MEJÍA

Abstract. Combining creature forcing approaches from [KS12] and [FGKS17], we show that, under CH, there is a proper $\omega^\omega$-bounding poset with $\aleph_2$-cc that forces continuum many pairwise different cardinal characteristics, parametrised by reals, for each one of the following six types: uniformity and covering numbers of Yorioka ideals as well as both kinds of localisation and anti-localisation cardinals, respectively. This answers several open questions from [KM22].

## 1. Introduction

The second and third authors [KM22] used a CS (countable support) product of lim sup creature-like forcings to prove that, consistently, uncountably many uniformity numbers of Yorioka ideals and uncountably many anti-localisation numbers $\mathfrak{v}^{\exists}_{c,h}$ are pairwise different. This paper is a continuation of [KM22] and solves several open questions.

For an increasing $f\colon \omega \to \omega$, the *Yorioka ideal parametrised by $f$* is denoted by $\mathcal{I}_f$ (see details in Definition 2.2). This is a $\sigma$-ideal on the Cantor space $2^\omega$ with a Borel base, which represents a fragment of the $\sigma$-ideal $\mathcal{SN}$ of strong measure zero subsets of $2^\omega$. In fact, $\mathcal{I}_f \subseteq \mathcal{N}$ (where $\mathcal{N}$ denotes the $\sigma$-ideal of measure zero subsets) and $\mathcal{SN} = \bigcap\{\mathcal{I}_f \mid f\colon \omega \to \omega \text{ increasing}\}$. These were introduced by Yorioka [Yor02] to characterise $\operatorname{cof}(\mathcal{SN})$ and to show that no inequality between $\operatorname{cof}(\mathcal{SN})$ and $\mathfrak{c} := 2^{\aleph_0}$ can be proved in ZFC.

Figure 1 illustrates the known relations between the cardinal characteristics associated with the Yorioka ideals and those in Cichoń's diagram. See details in [Osu06, KO08, Osu08, CM19].

The localisation and anti-localisation cardinals are defined as follows.

**Notation 1.1.** We fix the following terminology.

(1) Given a formula $\phi$ and a set $D$ (usually countably infinite), "$\forall^\infty \ell \in D\colon \phi$" means that all but finitely many members of $D$ satisfy $\phi$; "$\exists^\infty \ell \in D\colon \phi$" means that infinitely many members of $D$ satisfy $\phi$. We write "$\forall^\infty \ell\colon \phi$" and "$\exists^\infty \ell\colon \phi$" when the domain of discourse $D$ is clear.

2010 *Mathematics Subject Classification.* Primary 03E17; Secondary 03E05, 03E35, 03E40.
The first author was partially supported by the Austrian Science Fund (FWF) P30666 and is the recipient of a DOC fellowship of the Austrian Academy of Sciences at the Institute of Discrete Mathematics and Geometry, TU Wien, and the third author was supported by the Grant-in-Aid for Early Career Scientists 18K13448, Japan Society for the Promotion of Science.

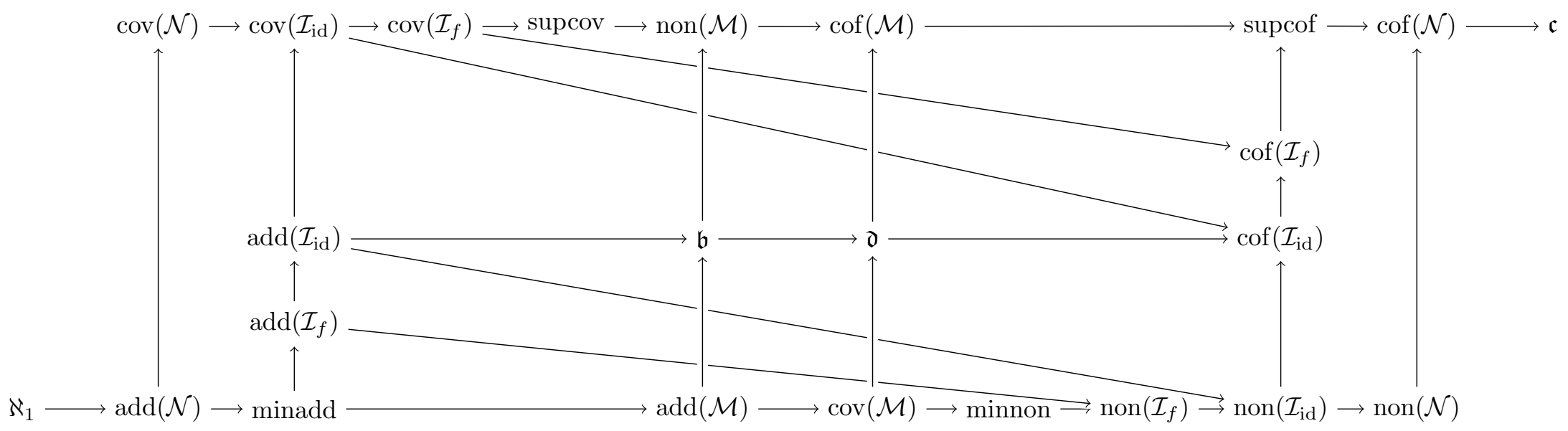

FIGURE 1. Cichoń's diagram with the cardinal characteristics associated with the Yorioka ideals. Here, minadd is the minimum over $\mathrm{add}(\mathcal{I}_f)$ for increasing $f\colon \omega \to \omega$; the cardinals supcov, minnon and supcof are defined likewise.

(2) When $b = \langle b(\ell) \mid \ell \in D\rangle$ is a sequence of non-empty sets and $h\colon D \to \omega$, let

$$\prod b := \prod_{\ell \in D} b(\ell), \qquad \mathcal{S}(b,h) := \prod_{\ell\in D} [b(\ell)]^{\leq h(\ell)}.$$

(3) A *slalom* is a function $\varphi\colon D \to [\omega]^{<\omega}$ for some countably infinite $D$.

(4) When $x$ and $\varphi$ are functions with domain $D$, we write
   (i) $x \in^* \varphi$ iff $\forall^\infty \ell \in D\colon x(\ell) \in \varphi(\ell)$, which is read as "$\varphi$ *localises* $x$";
   (ii) $x \in^\infty \varphi$ iff $\exists^\infty \ell \in D\colon x(\ell) \in \varphi(\ell)$. The expression $x \notin^\infty \varphi$ is read as "$\varphi$ *anti-localises* $x$".

(5) When $x, y \in \omega^D$, $x \leq^* y$ means $\forall^\infty \ell\colon x(\ell) \leq y(\ell)$.

**Definition 1.2.** Let $b = \langle b(\ell) \mid \ell \in D\rangle$ be a sequence of non-empty sets and let $h\colon D \to \omega$. We define the following cardinal characteristics:

$$\mathfrak{c}^\forall_{b,h} := \min\left\{|C| \;\middle|\; C \subseteq \mathcal{S}(b,h),\ \forall x \in \prod b\ \exists \varphi \in C\colon x \in^* \varphi\right\},$$

$$\mathfrak{v}^\forall_{b,h} := \min\left\{|F| \;\middle|\; F \subseteq \prod b,\ \neg\exists \varphi \in \mathcal{S}(b,h)\ \forall x \in F\colon x \in^* \varphi\right\},$$

$$\mathfrak{c}^\exists_{b,h} := \min\left\{|F| \;\middle|\; F \subseteq \mathcal{S}(b,h),\ \forall x \in \prod b\ \exists \varphi \in F\colon x \in^\infty \varphi\right\},$$

$$\mathfrak{v}^\exists_{b,h} := \min\left\{|C| \;\middle|\; C \subseteq \prod b,\ \forall \varphi \in \mathcal{S}(b,h)\ \exists x \in C\colon x \notin^\infty \varphi\right\}.$$

The first two cardinals are called *localisation cardinals*, while the latter two are called *anti-localisation cardinals.*

Notice that $\mathfrak{c}^\forall_{b,h}$ and $\mathfrak{c}^\exists_{b,h}$ are increasing on $b$ and decreasing on $h$ in the sense that if $\forall^\infty \ell\colon |b(\ell)| \leq |b'(\ell)|$ and $h' \leq^* h$, then $\mathfrak{c}^\forall_{b,h} \leq \mathfrak{c}^\forall_{b',h'}$ and $\mathfrak{c}^\exists_{b,h} \leq \mathfrak{c}^\exists_{b',h'}$. On the other hand, $\mathfrak{v}^\forall$ and $\mathfrak{v}^\exists$ are decreasing on $b$ and increasing on $h$. Also note that $\mathfrak{c}^\exists_{b,h} \leq \mathfrak{c}^\forall_{b,h}$ and $\mathfrak{v}^\forall_{b,h} \leq \mathfrak{v}^\exists_{b,h}$. See more in Lemma 2.10 and Figure 2.

We are interested in these cardinals when $D = \omega$. However, for technical reasons, in this work we would like to allow $D$ to be any countably infinite set (and we note that this does not have any relevant effects, anyway); but for this introductory section, let $D = \omega$.

The localisation and anti-localisation cardinals have appeared in many contexts. The following are well-known characterisations (see e. g. [CM23, Sections 4 and 5]):

- (Bartoszyński [Bar84]) $\operatorname{add}(\mathcal{N}) = \mathfrak{v}^{\forall}_{\omega,h}$ and $\operatorname{cof}(\mathcal{N}) = \mathfrak{c}^{\forall}_{\omega,h}$ when $h$ diverges to infinity (here $\omega$ is interpreted as the constant sequence $\omega$), and
- (Bartoszyński [Bar87], Miller [Mil82]) $\operatorname{non}(\mathcal{M}) = \mathfrak{c}^{\exists}_{\omega,h}$ and $\operatorname{cov}(\mathcal{M}) = \mathfrak{v}^{\exists}_{\omega,h}$ when $h \geq^* 1$, where $\mathcal{M}$ denotes the $\sigma$-ideal of meagre subsets of $2^\omega$.

In general, when $b(n)$ is infinite for infinitely many $n$, the localisation and anti-localisation cardinals coincide with other well-known cardinal characteristics (for an overview, see [CM19, Section 3] and [CM23]), so the interesting case is when $0 < h(n) < |b(n)| < \omega$ for all (but finitely many) $n$.[1] Miller [Mil81] proved $\operatorname{non}(\mathcal{SN}) = \min\{\mathfrak{v}^{\exists}_{b,h} \mid b \in \omega^\omega\}$ when $h \geq^* 1$, and $\operatorname{add}(\mathcal{M}) = \min\{\mathfrak{b}, \operatorname{non}(\mathcal{SN})\}$, which can be dualised to $\operatorname{cof}(\mathcal{M}) = \sup(\{\mathfrak{d}\} \cup \{\mathfrak{c}^{\exists}_{b,h} \mid b \in \omega^\omega\})$ (see [CM19, Theorem 3.23]).[2] On the other hand, due to Bartoszyński and Shelah [BS92] (see also [CM23, Section 7]), $\mathfrak{v}^{\forall}_{b,h} \leq \operatorname{non}(\mathcal{E})$ and $\operatorname{cov}(\mathcal{E}) \leq \mathfrak{c}^{\forall}_{b,h}$ whenever $\prod_{n<\omega} \frac{h(n)}{b(n)} = 0$, where $\mathcal{E}$ is the $\sigma$-ideal generated by the measure zero closed subsets of $2^\omega$; hence the localisation cardinals (and some variations thereof) play an important role for characterising $\operatorname{non}(\mathcal{E})$ and $\operatorname{cov}(\mathcal{E})$.

Figure 2 illustrates the provable inequalities among localisation cardinals (when $h$ diverges to infinity), anti-localisation cardinals (when $h(n)/b(n)$ converges to 0), and the cardinals in Cichoń's diagram (see [CM19] for a summary). In addition, as hinted in [GS93] and proved in [KM22, Lemma 2.3] and [CM23, Lemma 6.2]: If $\sum_{n<\omega} \frac{h(n)}{b(n)}$ converges, then $\operatorname{cov}(\mathcal{N}) \leq \mathfrak{c}^{\exists}_{b,h}$ and $\mathfrak{v}^{\exists}_{b,h} \leq \operatorname{non}(\mathcal{N})$; and if $\sum_{n<\omega} \frac{h(n)}{b(n)}$ diverges, then $\operatorname{cov}(\mathcal{E}) \leq \mathfrak{v}^{\exists}_{b,h}$ and $\mathfrak{c}^{\exists}_{b,h} \leq \operatorname{non}(\mathcal{E})$.[3]

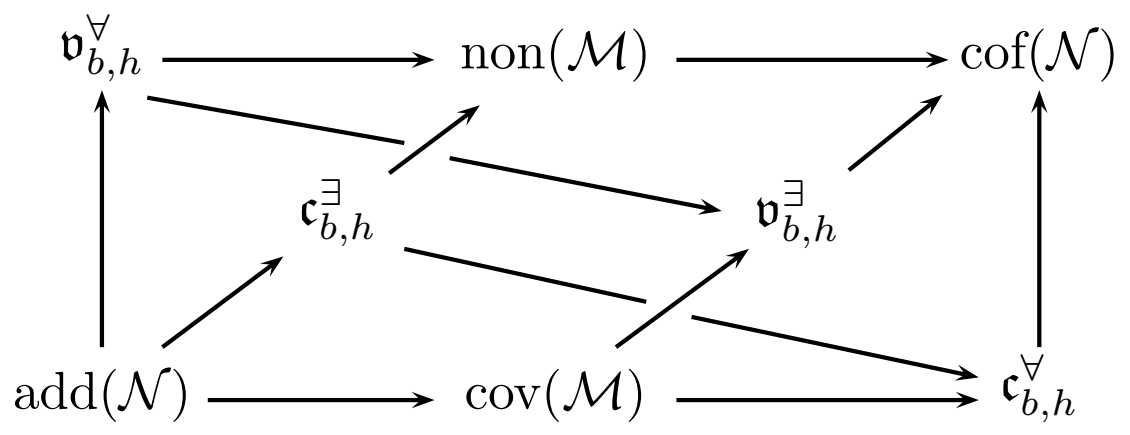

FIGURE 2. The ZFC-provable inequalities among (non-trivial) localisation cardinals, anti-localisation cardinals and the cardinals in Cichoń's diagram.

The purpose of this work is to prove the following result.

**Main Theorem.** *Assume* CH. *Then there is an $\aleph_2$-cc $\omega^\omega$-bounding proper poset forcing that there are continuum many pairwise different cardinal characteristics of each one of the six types $\mathfrak{c}^{\forall}_{b,h}$, $\mathfrak{c}^{\exists}_{b,h}$, $\mathfrak{v}^{\forall}_{b,h}$, $\mathfrak{v}^{\exists}_{b,h}$, $\operatorname{non}(\mathcal{I}_f)$ and $\operatorname{cov}(\mathcal{I}_f)$.*

[1] There are also other trivial cases: If $h$ does not diverge to infinity, then $\mathfrak{c}^{\forall}_{b,h} = \mathfrak{c}$ and $\mathfrak{v}^{\forall}_{b,h}$ is finite (Goldstern and Shelah [GS93]); and if the quotient $h(n)/b(n)$ does not converge to 0, then $\mathfrak{c}^{\exists}_{b,h}$ is finite and $\mathfrak{v}^{\exists}_{b,h} = \mathfrak{c}$ (see [CM19, Section 3]). More details can be found in [CM23, Section 3].

[2] However, it is consistent that $\operatorname{cof}(\mathcal{M}) < \operatorname{add}(\mathcal{SN})$, see [GJS93], so we cannot prove $\operatorname{cof}(\mathcal{M}) = \max\{\mathfrak{d}, \operatorname{cov}(\mathcal{SN})\}$.

[3] Since $\mathcal{E} \subseteq \mathcal{M} \cap \mathcal{N}$, it is clear that $\operatorname{cov}(\mathcal{E})$ is above the covering numbers of $\mathcal{N}$ and $\mathcal{M}$, and that $\operatorname{non}(\mathcal{E})$ is below their uniformity numbers.

**Observation 1.3** (about the notation)**.** The cardinal $\mathfrak{c}^{\forall}_{b,h}$ appears first (explicitly) in [GS93], denoted by $\mathbf{c}(b,h)$, where the $\mathbf{c}$ seem to stand for "cover" (according to [GS93, Definition 0.3]). Years later, Kellner starts using the notation $c^{\forall}(b,h)$ for "$\forall$-cover" [Kel08, Definition 2.1], and also mentions $c^{\exists}(b,h)$ (without definition), which was later introduced in [KS09] while changing the notation to $c^{\forall}_{b,h}$ and $c^{\exists}_{b,h}$. This notation is kept in [KS12], and later in [KO14] the fraktur font is used for the first time. Close to the first notation, $\mathfrak{c}_{b,h}$ is used for $\mathfrak{c}^{\forall}_{b,h}$ in [GK21].

The first and third authors studied these cardinals under the relational systems *localisation* and *anti-localisation* (see Section 2), but with different notation: The bounding and dominating numbers for localisation are denoted by $\mathfrak{b}^{\mathrm{Lc}}_{b,h}$ and $\mathfrak{d}^{\mathrm{Lc}}_{b,h}$, which correspond to $\mathfrak{v}^{\forall}_{b,h}$ and $\mathfrak{c}^{\forall}_{b,h}$, respectively; and the bounding and dominating numbers for anti-localisation are denoted by $\mathfrak{b}^{\mathrm{aLc}}_{b,h}$ and $\mathfrak{d}^{\mathrm{aLc}}_{b,h}$, which correspond to $\mathfrak{c}^{\exists}_{b,h}$ and $\mathfrak{v}^{\exists}_{b,h}$, respectively. Most of the work involving these authors uses this notation.

The notation $\mathfrak{v}^{\forall}_{b,h}$ and $\mathfrak{v}^{\exists}_{b,h}$ is used for the first time in [KM22], where the fraktur "v" is taken from the second letter of the word "evasion" (since the anti-localisation relation is a form of evasion). We decided to keep the notation of [KM22] for this paper.

In joint work in progress with Slovak researchers, we use the notation

$$\mathfrak{sl}_{\mathrm{t}}(b,h,\mathrm{Fin}) = \mathfrak{c}^{\forall}_{b,h}, \qquad \mathfrak{sl}^{\perp}_{\mathrm{t}}(b,h,\mathrm{Fin}) = \mathfrak{v}^{\forall}_{b,h},$$
$$\mathfrak{sl}_{\mathrm{e}}(b,h,\mathrm{Fin}) = \mathfrak{c}^{\exists}_{b,h}, \text{ and } \qquad \mathfrak{sl}^{\perp}_{\mathrm{e}}(b,h,\mathrm{Fin}) = \mathfrak{v}^{\exists}_{b,h}.$$

This is influenced by [Š23] towards the study of a more general version of localisation and anti-localisation cardinals, where $\mathrm{Fin} := [\omega]^{<\aleph_0}$ can be replaced by any ideal on $\omega$.

**Motivation.** In the direction of this paper's topic, Goldstern and Shelah [GS93] used a CS product of lim sup tree posets to force that $\aleph_1$ many cardinals of the type $\mathfrak{c}^{\forall}_{b,h}$ are pairwise different, answering a question of Blass [Bla93]. Here, "lim sup tree" refers to trees like in Sacks' and Miller's perfect trees posets, where splitting can be delayed from any node. This result was improved by Kellner [Kel08], who showed that, consistently, continuum many cardinals of the type $\mathfrak{c}^{\forall}_{b,h}$ are pairwise different. Later, Kellner and Shelah [KS09] included $\mathfrak{c}^{\exists}_{b,h}$ in this line of research: They introduced a CS product-like *liminf* forcing construction with *decisive creatures with halving* to force that there are $\aleph_1$ many pairwise different cardinals of each of the types $\mathfrak{c}^{\forall}_{b,h}$ and $\mathfrak{c}^{\exists}_{b,h}$, which was improved in [KS12] to continuum many using the same forcing technique, but with creatures that determine the parameters $b$ and $h$ generically. Here, "lim inf tree" means that splitting becomes recurrent from some point, like trees in Laver's poset. We remark here that the construction of the suitable decisive creatures in these latter two references is very complex, which makes it rather difficult to understand the intuition behind the main forcing construction.

Concerning Yorioka ideals, Kamo and Osuga [KO14] discovered connections between these ideals and anti-localisation cardinals – in particular, $\mathrm{cov}(\mathcal{I}_f)$ lies between two cardinals of the form $\mathfrak{c}^{\exists}_{b,h}$ (and dually, $\mathrm{non}(\mathcal{I}_f)$ lies between two cardinals of the form $\mathfrak{v}^{\exists}_{b,h}$). This is their starting point to force, via a FS (finite support) iteration of ccc posets, infinitely many pairwise different cardinal characteristics of each of the types $\mathrm{cov}(\mathcal{I}_f)$ and $\mathfrak{c}^{\exists}_{b,h}$, and even continuum many assuming the existence of a weakly inaccessible cardinal.

These techniques inspired Brendle and the third author [BM14] to force, via FS iterations of ccc posets, infinitely many pairwise different cardinal characteristics of the type $\mathfrak{v}^{\forall}_{b,h}$, even continuum many assuming the existence of a weakly inaccessible cardinal.[4] Much later, the first and third author [CM19] combined the methods of [KO14] and [BM14] to force infinitely many pairwise different cardinal characteristics of each of the types $\mathrm{cov}(\mathcal{I}_f)$, $\mathfrak{c}^{\exists}_{b,h}$ and $\mathfrak{v}^{\forall}_{b,h}$, even continuum many assuming the existence of a weakly inaccessible cardinal. In this framework, it is not possible to force continuum many different cardinals of these three types without using a weakly inaccessible cardinal: Usually, FS iterations of ccc posets are constructed with a length of uncountable cofinality, and this cofinality is forced between $\mathrm{non}(\mathcal{M})$ and $\mathrm{cov}(\mathcal{M})$; on the other hand, all these cardinal characteristics are below $\mathrm{non}(\mathcal{M})$ (see Figure 1 and Figure 2), so if continuum many pairwise different cardinals are forced, then the cofinality of the length of the iteration must be weakly inaccessible (and so it is forced that $\mathrm{non}(\mathcal{M}) = \mathfrak{c}$ is weakly inaccessible).

Very recently, the second and third authors [KM22] refined Kamo's and Osuga's connections between Yorioka ideals and anti-localisation cardinals (see e. g. Lemma 2.8) to prove that, consistently, there are $\aleph_1$ many cardinal characteristics of each of the types $\mathrm{non}(\mathcal{I}_f)$, $\mathfrak{v}^{\exists}_{b,h}$ and $\mathfrak{c}^{\forall}_{b,h}$. This leaves open the question on whether there could be continuum many different cardinals of each of the forms $\mathrm{non}(\mathcal{I}_f)$ and $\mathfrak{v}^{\exists}_{b,h}$, and also of the form $\mathrm{cov}(\mathcal{I}_f)$ and $\mathfrak{v}^{\forall}_{b,h}$ without using inaccessible cardinals (see [KM22, Questions C and E]). We answer all these questions in the positive – even more, as we claim in the Main Theorem, we construct a single model where there are continuum many different cardinals for each of the six types discussed.

**Method.** We mainly work with a creature forcing construction as in [GK21] mixed with lim inf creatures for parameters as in [KS12]. We review the evolution of these techniques for motivation. To force continuum many $\mathfrak{c}^{\forall}_{b,h}$ and $\mathfrak{c}^{\exists}_{b,h}$, Kellner and Shelah [KS12] proposed a mixed-limit creature forcing construction (i. e. one that mixes lim sup and lim inf creatures) that takes care of the continuum many parameters $\langle (b_\xi, h_\xi) \mid \xi < \mu \rangle$ while forcing $\mathfrak{c} = \mu$ and (pairwise different) values for $\mathfrak{c}^{\exists}_{b_\xi,h_\xi} = \mathfrak{c}^{\forall}_{b_\xi,h_\xi}$ (although they actually only use lim inf for their main theorem). Two very important features of lim inf creatures that come from [KS09] are halving and decisiveness: *Halving*, along with the more common notion of *bigness* (in the context of creatures), is essential to prove properness and continuous reading, while *decisiveness* is used to prove that the construction indeed increases $\mathfrak{c}^{\exists}_{b_\xi,h_\xi}$ to the desired value. As mentioned earlier, the construction of decisive creatures is quite elaborated, making it hard to keep track of the intuition connecting the construction to what is to be forced.

A bit later, A. Fischer, Goldstern, Kellner and Shelah [FGKS17] considerably simplified the mixed-limit creature forcing framework by allowing *subatomic creatures* in the construction of lim inf atomic creatures, ensuring decisiveness (implicitly, because the notion is not directly used) and a very intuitive construction of the creatures. This framework was used to force a constellation in Cichoń's diagram where five cardinal characteristics on the right side are pairwise different, namely $\mathrm{cov}(\mathcal{N}) = \mathfrak{d} = \aleph_1$, $\mathrm{non}(\mathcal{M})$, $\mathrm{non}(\mathcal{N})$, $\mathrm{cof}(\mathcal{N})$ and $\mathfrak{c}$ (see Theorem 9.1 and Figure 6).

[4] The original motivation in [BM14] was to force infinitely many pairwise different cardinal characteristics associated with Rothberger gaps in $F_\sigma$ ideals on $\omega$, but these turned out to be connected with the localisation cardinals of the form $\mathfrak{v}^{\forall}_{b,h}$.

Goldstern and the second author [GK21] improved [FGKS17] to force the same constellation of Cichoń's diagram plus $\aleph_1$ many different localisation cardinals of the form $\mathfrak{c}^{\forall}_{b,h}$. They succeeded to separate the components of the creatures and present the forcing construction as a CS product of one (large) lim inf part to increase $\mathrm{non}(\mathcal{M})$, one large lim sup part to increase the continuum (similar to a CS product of Sacks posets, original from [FGKS17]) and several single lim sup posets to increase $\mathrm{non}(\mathcal{N})$, $\mathrm{cof}(\mathcal{N})$ and the localisation cardinals.[5]

In this work, we combine the methods of [GK21] with a lim inf creature forcing as in [KS12] for parameters, and construct a large creature forcing to prove the Main Theorem. This construction is not a CS product in the strict sense, but it looks like one.

The first component, called the *pr-part*, is a large lim inf creature forcing (similar to the lim inf part of the forcing in [GK21]) that adds *blocks* of parameters.[6] A block is a finite sequence $\vec{\nu} = \langle d, h, g, b, f, a\rangle$ of increasing functions satisfying the requirements in Definition 2.11, whose important feature is that it satisfies the diagram of inequalities in Figure 3 (see Lemma 2.12).

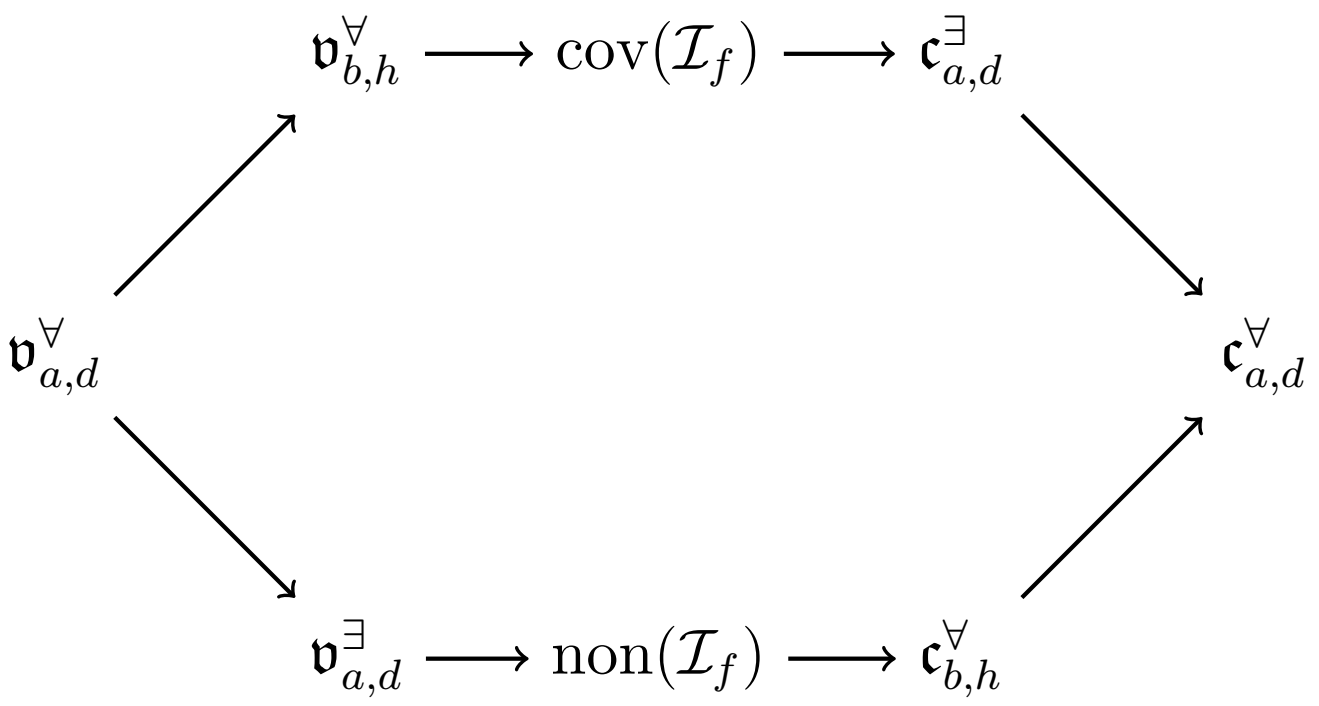

FIGURE 3. Diagram of inequalities for a block $\vec{\nu} = \langle d, h, g, b, f, a\rangle$.

With the notion of blocks, we can restate the Main Theorem more precisely.

**Theorem 1.4.** *Assume* CH. *Let $\mu = \mu^{\aleph_0}$ be an infinite cardinal, $S^{\mathrm{pr}}$ a set of size $\mu$ and, for $i \in S^{\mathrm{pr}}$, let $\kappa_i = \kappa_i^{\aleph_0}$ be an infinite cardinal.*

*Then there is a proper $\omega^\omega$-bounding $\aleph_2$-cc poset $\mathbb{Q}$ forcing, for all $i \in S^{\mathrm{pr}}$,*

[5] In fact, the subatomic creatures corresponding to $\mathrm{non}(\mathcal{M})$ in [FGKS17, GK21] come from those defined in [KS09, KS12] to increase $\mathfrak{c}^{\exists}_{b,h}$ (which we use in the current paper to add parameters); and the atomic creatures corresponding to $\mathrm{non}(\mathcal{N})$ actually increase some anti-localisation cardinal $\mathfrak{v}^{\exists}_{b,h}$ (that lies below $\mathrm{non}(\mathcal{N})$). We use an equivalent formulation of these atomic creatures in our framework.

[6] Although this can be done with lim sup forcing, it may not (easily) ensure a property called "separated support" (from [KS09, KS12]), which we require to separate the continuum many cardinal characteristics.

*(I)* $\mathfrak{c} = \mu$,
*(II)* *blocks* $\vec{\nu}_i^{\mathrm{lc}} = \langle d_i^{\mathrm{lc}}, h_i^{\mathrm{lc}}, g_i^{\mathrm{lc}}, b_i^{\mathrm{lc}}, f_i^{\mathrm{lc}}, a_i^{\mathrm{lc}}\rangle$ *and* $\vec{\nu}_i^{\mathrm{al}} = \langle h_i^{\mathrm{al}}, h_i^{\mathrm{al}}, g_i^{\mathrm{al}}, b_i^{\mathrm{al}}, f_i^{\mathrm{al}}, a_i^{\mathrm{al}}\rangle$,
*(III)* $\mathfrak{v}^{\forall}_{b_i^{\mathrm{lc}},h_i^{\mathrm{lc}}} = \mathrm{cov}(\mathcal{I}_{f_i^{\mathrm{lc}}}) = \mathfrak{c}^{\exists}_{a_i^{\mathrm{lc}},d_i^{\mathrm{lc}}} = \mathfrak{c}^{\forall}_{a_i^{\mathrm{lc}},d_i^{\mathrm{lc}}} = \kappa_i$ *(the upper part of Figure 3), and*
*(IV)* $\mathfrak{v}^{\exists}_{a_i^{\mathrm{al}},h_i^{\mathrm{al}}} = \mathrm{non}(\mathcal{I}_{f_i^{\mathrm{al}}}) = \mathfrak{c}^{\forall}_{a_i^{\mathrm{al}},h_i^{\mathrm{al}}} = \kappa_i$ *(the lower part of Figure 3).*

To prove this theorem, we use the pr-part to add the blocks in (II) (while forcing (I)), and we attach to the pr-part one (large) lim inf creature forcing to increase $\mathfrak{v}^{\forall}_{b_i^{\mathrm{lc}},h_i^{\mathrm{lc}}}$ for all $i$ as well as several individual lim sup creature forcings to increase $\mathfrak{v}^{\exists}_{a_i^{\mathrm{al}},h_i^{\mathrm{al}}}$. Although the pr-part can be isolated, the other components depend on the "possibilities" determined by the pr-part. The subatomic creatures used for these other components are the same ones used in [KM22] to manipulate $\mathfrak{v}^{\exists}_{b,h}$ (which also works for $\mathfrak{v}^{\forall}_{b,h}$ in the lim inf context).[7]

To force $\mathfrak{c}^{\forall}_{a_i^{\mathrm{lc}},d_i^{\mathrm{lc}}} \leq \kappa_i$ and $\mathfrak{c}^{\forall}_{a_i^{\mathrm{al}},h_i^{\mathrm{al}}} \leq \kappa_i$, we use similar arguments as in the cited works on creature forcing: We force that the set of slaloms that can be continuously read only using a fixed set of indices of size $\kappa_i$ is a witness of the cardinal characteristic. This is derived from a property we call $(\dot{a}_i, \dot{d}_i)$*-bounding over* $\{i\} \cup S_i^{\mathrm{lc}} \cup S_i^{\mathrm{al}}$, which is presented in Lemma 8.12.

Since our forcing construction is still quite complex, we need to motivate and explain more detailed features while building its different components.

**Structure of the paper.** We give a short overview of the following sections' contents.

Section 2: We review the definition of *Yorioka ideals* and present the notion of *blocks*.

Section 3: We build the frame in which the forcing lives, define the parameters of the forcing and which conditions they must fulfil so that we can prove Theorem 1.4.

Section 4: Inspired by the terminology of [FGKS17], we define subatomic creatures, atomic creatures and compound creatures, which are the building blocks of a condition in the forcing we want to construct.

Section 5: We present our forcing construction in detail and prove its basic properties. This construction follows the presentation in [GK21] and is divided into three components: the pr-part, the lc-part and the al-part. As in [KS12], the atomic creatures of the lc-part and al-part depend on the possibilities in the pr-part.

Section 6: The necessary bigness properties of the subatomic and compound creatures are developed. These will be essential to prove the main features of the forcing, like *pure decision*, *continuous reading of names* and *rapid reading*, and ultimately for the proof of Theorem 1.4.

[7] To prove Theorem 1.4, we expected to add only one block $\vec{\nu}_i = \langle d_i, h_i, g_i, b_i, f_i, a_i\rangle$ for each $i \in S^{\mathrm{pr}}$, and force that all cardinals in Figure 3 corresponding to $\vec{\nu}_i$ equal $\kappa_i$, by using (besides the pr-part) one large lim inf creature construction to increase $\mathfrak{v}^{\forall}_{a_i,d_i}$ for all $i$ (i. e. without lim sup forcings). However, we could not find the right construction of the parameters that allows this. Nonetheless, this may very well be possible to achieve.

Section 7: We present continuous and rapid reading for our forcing. We prove that continuous reading follows by pure decision as well as properness and $\omega^\omega$-bounding (in fact, our forcing satisfies a variant of the finite version of Baumgartner's axiom A). Rapid reading follows by continuous reading and the bigness results from Section 6. The presentation of this section is very similar to [GK21].

Section 8: We can finally prove Theorem 1.4.

Section 9: We present further discussion and open problems.

## 2. Blocks of Parameters

The main purpose in this section is to introduce the notion of a *block of parameters* and to prove the inequalities in Figure 3 (Lemma 2.12). They are consequences of connections between Yorioka ideals, localisation cardinals and anti-localisation cardinals that we present in this section.

We first recall the definition of a Yorioka ideal.

**Notation 2.1.** Let $\mathrm{pow}_k\colon \omega \to \omega$ denote the function defined by $\mathrm{pow}_k(i) := i^k$.

Define the relation $\ll$ on $\omega^\omega$ by

$$f \ll g \quad \text{iff} \quad \forall\, k < \omega\colon f \circ \mathrm{pow}_k \leq^* g.$$

**Definition 2.2** ([Yor02])**.** For $\sigma \in (2^{<\omega})^\omega$, define

$$[\sigma]_\infty := \{x \in 2^\omega \mid \exists^\infty\, n\colon \sigma(n) \subseteq x)\} = \bigcap_{n<\omega} \bigcup_{m \geqslant n} [\sigma(m)],$$

where $[s] := \{x \in 2^\omega \mid s \subseteq x\}$ for any $s \in 2^{<\omega}$. Let $\mathrm{ht}_\sigma \in \omega^\omega$ be the function defined by $\mathrm{ht}_\sigma(i) := |\sigma(i)|$ for each $i < \omega$.

Given $f \in \omega^\omega$, define

$$\mathcal{I}_f := \{X \subseteq 2^\omega \mid \exists\, \sigma \in (2^{<\omega})^\omega\colon X \subseteq [\sigma]_\infty \text{ and } \mathrm{ht}_\sigma \gg f\}.$$

Any family of the form $\mathcal{I}_f$ with $f$ increasing is called a *Yorioka ideal.*

(The proof that $\mathcal{I}_f$ is a $\sigma$-ideal for each increasing $f$ is far from trivial, see [Yor02, Lemma 3.4] and [CM19, Theorem 3.12].)

Classical cardinal characteristics of the continuum, as well as those we deal with in this work, can be expressed by relational systems.

**Definition 2.3.** A *relational system* is a triple $\mathbf{R} = \langle X, Y, R\rangle$ consisting of two non-empty sets $X, Y$ and a relation $R$.

The *dual of* $\mathbf{R}$ is the relational system $\mathbf{R}^\perp = \langle Y, X, R^\perp\rangle$ where $y\, R^\perp\, x$ is defined by $\neg(x\, R\, y)$.

The following are the *cardinal characteristics associated with* $\mathbf{R}$:

$$\mathfrak{b}(\mathbf{R}) := \min\{|F| \mid F \subseteq X,\ \neg\exists\, y \in Y\ \forall\, x \in F\colon x\, R\, y\},$$
$$\mathfrak{d}(\mathbf{R}) := \min\{|C| \mid C \subseteq X,\ \forall\, x \in X\ \exists\, y \in C\colon x\, R\, y\}.$$

It is clear that $\mathfrak{b}(\mathbf{R}^\perp) = \mathfrak{d}(\mathbf{R})$ and $\mathfrak{d}(\mathbf{R}^\perp) = \mathfrak{b}(\mathbf{R})$.

The localisation and anti-localisation cardinals are cardinals associated with relational systems as follows.

**Example 2.4.** Let $b$ be a function with domain $D$ such that $b(\ell) \neq \varnothing$ for all $\ell \in D$, and let $h\colon D \to \omega$.

(1) Let $\mathbf{Lc}(b,h) := \langle \prod b, \mathcal{S}(b,h), \in^* \rangle$ be a relational system. Then $\mathfrak{v}^{\forall}_{b,h} = \mathfrak{b}(\mathbf{Lc}(b,h))$ and $\mathfrak{c}^{\forall}_{b,h} = \mathfrak{d}(\mathbf{Lc}(b,h))$.
(2) Let $\mathbf{aLc}(b,h) := \langle \mathcal{S}(b,h), \prod b, \not\ni^{\infty} \rangle$ be a relational system. Then $\mathfrak{c}^{\exists}_{b,h} = \mathfrak{b}(\mathbf{aLc}(b,h))$ and $\mathfrak{v}^{\exists}_{b,h} = \mathfrak{d}(\mathbf{aLc}(b,h))$.

The cardinal characteristics associated with an ideal are also associated with relational systems.

**Example 2.5.** Let $X$ be a set and let $\mathcal{I}$ be an ideal on $X$ containing all finite subsets of $X$.

(1) Consider $\mathcal{I}$ as the relational system $\langle \mathcal{I}, \mathcal{I}, \subseteq \rangle$. Then $\mathfrak{b}(\mathcal{I}) = \mathrm{add}(\mathcal{I})$ and $\mathfrak{d}(\mathcal{I}) = \mathrm{cof}(\mathcal{I})$.
(2) Let $\mathbf{Cv}(\mathcal{I}) := \langle X, \mathcal{I}, \in \rangle$ be a relational system. Then $\mathfrak{b}(\mathbf{Cv}(\mathcal{I})) = \mathrm{cov}(\mathcal{I})$ and $\mathfrak{d}(\mathbf{Cv}(\mathcal{I})) = \mathrm{non}(\mathcal{I})$.

We review the concept of Tukey connections, which are useful to determine inequalities between cardinal characteristics.

**Definition 2.6** ([Voj93])**.** Let $\mathbf{R} = \langle X, Y, R \rangle$ and $\mathbf{R}' = \langle X', Y', R' \rangle$ be relational systems. We say that $\mathbf{R}$ *is Tukey below* $\mathbf{R}'$, denoted by $\mathbf{R} \preceq_{\mathrm{T}} \mathbf{R}'$, if there are maps $F\colon X \to X'$ and $G\colon Y' \to Y$ such that, for any $x \in X$ and $y' \in Y'$, if $F(x)\ R'\ y'$ then $x\ R\ G(y')$. In that case, we say that the pair $(F, G)$ is a *Tukey connection witnessing* $\mathbf{R} \preceq_{\mathrm{T}} \mathbf{R}'$.

We say that $\mathbf{R}$ and $\mathbf{R}'$ are *Tukey equivalent*, denoted by $\mathbf{R} \cong_{\mathrm{T}} \mathbf{R}'$, if $\mathbf{R} \preceq_{\mathrm{T}} \mathbf{R}'$ and $\mathbf{R}' \preceq_{\mathrm{T}} \mathbf{R}$.

**Lemma 2.7.** *With the notation in the previous definition, if* $\mathbf{R} \preceq_{\mathrm{T}} \mathbf{R}'$*, then* $(\mathbf{R}')^{\perp} \leq \mathbf{R}^{\perp}$ *as well as* $\mathfrak{d}(\mathbf{R}) \leq \mathfrak{d}(\mathbf{R}')$ *and* $\mathfrak{b}(\mathbf{R}') \leq \mathfrak{b}(\mathbf{R})$*. Also,* $(\mathbf{R}^{\perp})^{\perp} = \mathbf{R}$*.*

The Tukey order between the relational systems in Example 2.4 and $\mathbf{Cv}(\mathcal{I}_f)$ are key to our definition of blocks.

**Lemma 2.8** ([KM22, Lemma 2.4])**.** *Let* $a, d \in \omega^\omega$ *and let* $\langle I_n \mid n < \omega \rangle$ *be an interval partition of* $\omega$ *such that* $|I_n| = d(n)$*. Define* $g_{a,d} \in \omega^\omega$ *by* $g_{a,d}(k) := \lfloor \log_2 a(n) \rfloor$ *for each* $k \in I_n$*. If* $a \geq^* 2$*,* $d \geq^* 1$*,* $f$ *is an increasing function and* $g_{a,d} \gg f$*, then* $\mathbf{Cv}(\mathcal{I}_f) \preceq_{\mathrm{T}} \mathbf{aLc}(a,d)^{\perp}$*. In particular,* $\mathrm{cov}(\mathcal{I}_f) \leq \mathfrak{c}^{\exists}_{a,d}$ *and* $\mathfrak{v}^{\exists}_{a,d} \leq \mathrm{non}(\mathcal{I}_f)$*.*

The following is a variation of [KM22, Lemma 2.5] where we use a relational system for localisation instead of anti-localisation.

**Lemma 2.9.** *Let* $b, g, h \in \omega^\omega$ *and let* $\langle J_n \mid n < \omega \rangle$ *be an interval partition of* $\omega$ *such that* $|J_n| = g(n)$*. Define* $f_{b,g,h} \in \omega^\omega$ *by* $f_{b,g,h}(k) := \sum_{\ell \leq n} h(\ell) \lceil \log_2 b(\ell) \rceil$ *for each* $k \in J_n$*. If* $b \geq^* 2$*,* $g \geq^* 1$*,* $h \geq^* 1$*,* $f \in \omega^\omega$ *is an increasing function and there is some* $1 \leq m < \omega$ *such that* $f_{b,g,h}(k) \leq f(k^m)$ *for all but finitely many* $k < \omega$*, and* $g(n)h(n) < b(n)$ *for all but finitely many* $n < \omega$*, then* $\mathbf{Lc}(b,h)^{\perp} \preceq_{\mathrm{T}} \mathbf{Cv}(\mathcal{I}_f)$*. In particular,* $\mathfrak{v}^{\forall}_{b,h} \leq \mathrm{cov}(\mathcal{I}_f)$ *and* $\mathrm{non}(\mathcal{I}_f) \leq \mathfrak{c}^{\forall}_{b,h}$*.*

*Proof.* It is enough to construct two functions $F\colon \mathcal{S}(b,h)\to 2^\omega$ and $G\colon \mathcal{I}_f \to \prod b$ such that, for all $S\in\mathcal{S}(b,h)$ and $X\in\mathcal{I}_f$, $F(S)\in X$ implies $G(X)\notin^* S$.

For each $n<\omega$, fix a one-to-one function $\iota_n\colon b(n)\to 2^{\lceil\log_2 b(n)\rceil}$. Here, we interpret $2^{\lceil\log_2 b(n)\rceil}$ as the set of functions from $\lceil\log_2 b(n)\rceil$ into $2=\{0,1\}$.

Let $S\in\mathcal{S}(b,h)$ and choose some $S'\in\mathcal{S}(b,h)$ such that $S(n)\subseteq S'(n)\neq\varnothing$ whenever $h(n)\neq 0$. For each $n$, let $S'(n)=\{m_{n,j}\mid j<h(n)\}$, and define $F(S)$ as the concatenation of $\{\iota_n(m_{i,j})\mid i<\omega,\ j<h(n)\}$, where the indices $(i,j)$ are ordered lexicographically.

To define $G$, let $X\in\mathcal{I}_f$, so choose $\sigma_X\in(2^{<\omega})^\omega$ such that $\mathrm{ht}_{\sigma_X}\gg f$ and $X\subseteq[\sigma_X]_\infty$. The hypothesis of the lemma implies that $\mathrm{ht}_{\sigma_X}\geq^* f_{b,g,h}$. Consider the interval partition $\langle I_{n,\ell}\mid n<\omega,\ \ell<h(\ell)\rangle$ of $\omega$, ordered lexicographically, such that $|I_{n,\ell}|=\lceil\log_2 b(n)\rceil$. Then $I_{n,\ell}\subseteq\mathrm{ht}_{\sigma_X}(k)$ for all $k\in J_n$ and all but finitely many $n$ (because $f_{b,g,h}(k)=\sum_{\ell\leq n}h(l)\lceil\log_2 b(\ell)\rceil\leq\mathrm{ht}_{\sigma_X}(k)$), so we can define

$$H(n):=\{{\iota_n}^{-1}(\sigma_X(k){\restriction}_{I_{n,\ell}})\mid \ell<h(n),\ k\in J_n,\ I_{n,\ell}\subseteq\mathrm{ht}_{\sigma_X}(k),\ \sigma_X(k){\restriction}_{I_{n,\ell}}\in\operatorname{ran}\iota_n\}.$$

Clearly $H(n)\subseteq b(n)$ and, eventually, $|H(n)|\leq g(n)h(n)<b(n)$, so we can choose some $G(X)(n)\in b(n)\smallsetminus H(n)$. For the other finitely many $n$, choose any $G(X)(n)\in b(n)$.

Now assume $S\in\mathcal{S}(b,h)$, $X\in\mathcal{I}_f$ and $F(S)\in X$; we will show that $G(X)(n)\notin S(n)$ for infinitely many $n$. Since $F(S)\in X$, there are infinitely many $n$ such that $\sigma_X(k_n)\subseteq F(S)$ for some $k_n\in J_n$. For such a sufficiently large $n$, $\iota_n(m_{n,j})=\sigma_X(k_n){\restriction}_{I_{n,j}}$ for all $j<h(n)$, so $m_{n,j}\in H(n)$, which implies $G(X)(n)\neq m_{n,j}$. Thus $G(X)(n)\notin S(n)$. □

We also consider the following easy fact.

**Lemma 2.10.** *Let $a,b,d,h\in\omega^\omega$. Then the following statements hold:*

*(a)* $\mathbf{aLc}(a,d)^\perp\preceq_{\mathrm{T}}\mathbf{Lc}(a,d)$*; in particular,* $\mathfrak{c}^\exists_{a,d}\leq\mathfrak{c}^\forall_{a,d}$ *and* $\mathfrak{v}^\forall_{a,d}\leq\mathfrak{v}^\exists_{a,d}$.

*(b) If* $b\leq^* a$ *and* $d\leq^* h$*, then* $\mathbf{Lc}(b,h)\preceq_{\mathrm{T}}\mathbf{Lc}(a,d)$*; in particular,* $\mathfrak{v}^\forall_{a,d}\leq\mathfrak{v}^\forall_{b,h}$ *and* $\mathfrak{c}^\forall_{b,h}\leq\mathfrak{c}^\forall_{a,d}$.

Using the hypotheses of Lemma 2.8–2.10, we are finally ready to define the notion of a block of parameters. Properties (i)–(v) are chosen this way precisely to get the promised Figure 3 as a consequence of the following lemmata.

**Definition 2.11.** A sequence of increasing functions $\vec{\nu}=\langle d,h,g,b,f,a\rangle$ in $\omega^\omega$ is a *block (of parameters)* if it fulfils, for all $n<\omega$,

(i) if $k\in I_n$, then $k^n\leq\max J_n$, where $\langle I_n\mid n<\omega\rangle$ and $\langle J_n\mid n<\omega\rangle$ are interval partitions of $\omega$ such that $|I_n|=d(n)$ and $|J_n|=g(n)$;
(ii) $g(n)h(n)<b(n)$;
(iii) if $k\in J_n$, then $f(k)\geq\sum_{j\leq n}h(j)\lceil\log_2 b(j)\rceil$;
(iv) if $k\in I_n$, then $f(k^n)\leq\lfloor\log_2 a(n)\rfloor$; and
(v) $b(n)\leq a(n)$ and $d(n)\leq h(n)$.

**Lemma 2.12.** *If $\vec{\nu}=\langle d,h,g,b,f,a\rangle$ is a block, then* $\mathfrak{v}^\forall_{a,d}\leq\mathfrak{v}^\forall_{b,h}\leq\operatorname{cov}(\mathcal{I}_f)\leq\mathfrak{c}^\exists_{a,d}\leq\mathfrak{c}^\forall_{a,d}$ *and* $\mathfrak{v}^\forall_{a,d}\leq\mathfrak{v}^\exists_{a,d}\leq\operatorname{non}(\mathcal{I}_f)\leq\mathfrak{c}^\forall_{b,h}\leq\mathfrak{c}^\forall_{a,d}$ *(cf. Figure 3).*

*Proof.* We show that $\mathbf{Lc}(a,d)^\perp \preceq_{\mathrm{T}} \mathbf{aLc}(a,d) \preceq_{\mathrm{T}} \mathbf{Cv}(\mathcal{I}_f)^\perp \preceq_{\mathrm{T}} \mathbf{Lc}(b,h) \preceq_{\mathrm{T}} \mathbf{Lc}(a,d)$. By Definition 2.11 (v) and Lemma 2.10, we have that $\mathbf{Lc}(a,d)^\perp \preceq_{\mathrm{T}} \mathbf{aLc}(a,d)$ and $\mathbf{Lc}(b,h) \preceq_{\mathrm{T}} \mathbf{Lc}(a,d)$ and, by Definition 2.11 (ii)–(iii) and Lemma 2.9, $\mathbf{Cv}(\mathcal{I}_f)^\perp \preceq_{\mathrm{T}} \mathbf{Lc}(b,h)$ follows. To show $\mathbf{aLc}(a,d) \preceq_{\mathrm{T}} \mathbf{Cv}(\mathcal{I}_f)^\perp$ by application of Lemma 2.8, it remains to prove that $g_{a,d} \gg f$. Fix $m<\omega$; for $n>m$, if $k\in I_n$, then $f(k^m) \leq f(k^n) \leq \lfloor \log_2 a(n) \rfloor = g_{a,d}(k)$. □

Note that the requirement "if $k\in I_n$ then $k^n \leq \max J_n$" was not used in this proof, but we will need it to construct blocks in Lemma 3.4 (specifically, to guarantee (iv) of Definition 2.11).

As mentioned in the introduction, we plan to use lim sup forcing to increase the $\mathfrak{d}$ cardinals of (the relational systems determined by) a block, and lim inf forcing for the $\mathfrak{b}$ cardinals (specifically, to increase $\mathfrak{v}^{\forall}_{b,h}$). For the latter, we consider a variation of the relational system $\mathbf{Lc}(b,h)$.

**Definition 2.13.** Let $\bar{I} := \langle I_n \mid n<\omega\rangle$ be a sequence of pairwise disjoint non-empty finite sets and let $D := \bigcup_{n<\omega} I_n$.

(1) For two functions $x$ and $\varphi$ with domain $D$, we write

$$x \in^*_{\bar{I}} \varphi \text{ iff } \forall^\infty n\ \exists \ell \in I_n\colon x(\ell)\in\varphi(\ell).$$

(2) Let $b=\langle b(\ell)\mid \ell\in D\rangle$ be a sequence of non-empty sets and let $h\colon D\to\omega$. Let $\mathbf{Lc}_{\bar{I}}(b,h) := \langle \prod b, \mathcal{S}(b,h), \in^*_{\bar{I}}\rangle$ be a relational system.

(3) We define $\mathfrak{v}^{\forall,\bar{I}}_{b,h} := \mathfrak{b}(\mathbf{Lc}_{\bar{I}}(b,h))$ and $\mathfrak{c}^{\forall,\bar{I}}_{b,h} := \mathfrak{d}(\mathbf{Lc}_{\bar{I}}(b,h))$.

In our forcing construction, we intend to increase cardinals of the form $\mathfrak{v}^{\forall}_{b,h}$, via lim inf forcing, by adding a generic slalom $\varphi$. In this lim inf part of the forcing $\varphi(n)$ is decided by information at level $n$ of a condition of the forcing. However, in our construction, this level $n$ is divided into a finite set $I_n$ of sublevels, so we do not obtain a single $\varphi(n)$ in principle, but a sequence $\langle\varphi(\ell)\mid \ell\in I_n\rangle$, and so we end up with a slalom with domain $D$ as above. On the other hand, for any $x\in\prod_{\ell\in D} b(\ell)$, we cannot guarantee that $\forall^\infty n<\omega\ \forall \ell\in I_n\colon x(\ell)\in\varphi(\ell)$, but the best we can say is that for all but finitely many $n<\omega$, $x(\ell)\in\varphi(\ell)$ for some (even many, but not *all*) $\ell\in I_n$. This is the main motivation to introduce Definition 2.13, because our forcing increases $\mathfrak{v}^{\forall,\bar{I}}_{b,h}$ directly. This has the desired effect on cardinals of the form $\mathfrak{v}^{\forall}_{b,h}$ as well, according to the following result.

**Lemma 2.14.** *With the notation from the previous definition, let $b^* := \langle b^*(n)\mid n<\omega\rangle$, $b^*(n) := \prod_{\ell\in I_n} b(\ell)$, and define $h^*\in\omega^\omega$ by $h^*(n) := b^*(n) - \prod_{\ell\in I_n}(b(\ell)-h(\ell))$. Then $\mathbf{Lc}(b^*,h^*) \preceq_{\mathrm{T}} \mathbf{Lc}_{\bar{I}}(b,h)$. In particular, $\mathfrak{v}^{\forall,\bar{I}}_{b,h} \leq \mathfrak{v}^{\forall}_{b^*,h^*}$ and $\mathfrak{c}^{\forall}_{b^*,h^*} \leq \mathfrak{c}^{\forall,\bar{I}}_{b,h}$.*

*Proof.* Define $F\colon \prod b^* \to \prod b$ by $F(x) := \langle x(n)(\ell) \mid \ell\in I_n,\ n<\omega\rangle$, and $G\colon \mathcal{S}(b,h)\to\mathcal{S}(b^*,h^*)$ by $G(\varphi) := \varphi^*$ where $\varphi^*(n) = b^*(n)\smallsetminus \prod_{\ell\in I_n}(b(\ell)\smallsetminus\varphi(\ell))$ (which clearly has size $\leq h^*(n)$). For $x\in\prod b^*$ and $\varphi\in\mathcal{S}(b,h)$, $x(n)\in\varphi^*(n)$ iff $\exists \ell\in I_n\colon x(n)(\ell)\in\varphi(\ell)$ for all $n<\omega$. Therefore, $F(x)\in^*_{\bar{I}}\varphi$ iff $x\in^*\varphi^*$. □

**Observation 2.15.** The slaloms we add in our forcing construction are not quite in $\mathcal{S}(b,h)$, but they are modulo finite modifications. For $b=\langle b(\ell)\mid \ell\in D\rangle$ and $h\colon D\to\omega$, define

$$\mathcal{S}^*(b,h) := \big\{\varphi : D\to[\omega]^{<\aleph_0} \mid \{\ell\in D \mid \varphi(\ell)\not\subseteq b(\ell) \text{ or } |\varphi(\ell)|>h(\ell)\} \text{ is finite}\big\}.$$

Our forcing will add slaloms in $\mathcal{S}^*(b,h)$. If we replace $\mathcal{S}(b,h)$ with $\mathcal{S}^*(b,h)$ in the relational systems $\mathbf{Lc}(b,h)$ and $\mathbf{aLc}(b,h)$, then we obtain Tukey equivalent relational systems (and the same localisation and anti-localisation cardinals); likewise for the relational system $\mathbf{Lc}_{\bar{I}}(b,h)$.

## 3. The Frame and Parameters of the Forcing

In this section we present the frame in which the forcing lives and the parameters we need for the forcing construction.

3.1. **The frame.** Our forcing construction has a support $S^*$; for each $\alpha \in S^*$, we add blocks or slaloms that increase cardinals of the form $\mathfrak{v}^{\forall,\bar{I}}_{b,h}$ and $\mathfrak{v}^{\exists}_{a,h}$. Concretely, we fix:

- a set of indices $S^{\mathrm{pr}}$ (where pr stands for "parameter");
- pairwise disjoint families $\langle S^{\mathrm{lc}}_i \mid i \in S^{\mathrm{pr}}\rangle$ and $\langle S^{\mathrm{al}}_i \mid i \in S^{\mathrm{pr}}\rangle$ (where lc stands for "localisation", and al for "anti-localisation") such that $S^{\mathrm{lc}}_i \cap S^{\mathrm{pr}} = S^{\mathrm{al}}_i \cap S^{\mathrm{pr}} = S^{\mathrm{lc}}_i \cap S^{\mathrm{al}}_j = \varnothing$ for $i,j \in S^{\mathrm{pr}}$; and
- $S^* := \bigcup_{\mathrm{t}\in\{\mathrm{pr},\mathrm{lc},\mathrm{al}\}} S^{\mathrm{t}}$, where $S^{\mathrm{t}} := \bigcup_{i\in S^{\mathrm{pr}}} S^{\mathrm{t}}_i$ for $\mathrm{t} \in \{\mathrm{lc},\mathrm{al}\}$.

For each $i \in S^{\mathrm{pr}}$, we intend to force several blocks $\vec{\nu}^{\mathrm{lc}}_i = \langle d^{\mathrm{lc}}_i, h^{\mathrm{lc}}_i, g^{\mathrm{lc}}_i, b^{\mathrm{lc}}_i, f^{\mathrm{lc}}_i, a^{\mathrm{lc}}_i\rangle$ and $\vec{\nu}^{\mathrm{al}}_i = \langle h^{\mathrm{al}}_i, h^{\mathrm{al}}_i, g^{\mathrm{al}}_i, b^{\mathrm{al}}_i, f^{\mathrm{al}}_i, a^{\mathrm{al}}_i\rangle$; in fact, we add an $\omega$-sequence $y_i = \langle y_i(n) \mid n < \omega\rangle$ that defines both blocks. Definition 3.3 and Lemma 3.4 give details on how blocks are defined from reals.

For each $i \in S^{\mathrm{lc}}$ and each $\alpha \in S^{\mathrm{lc}}_i$, we add a slalom $\varphi_\alpha \in \mathcal{S}^*(b^{\mathrm{lc}}_i, h^{\mathrm{lc}}_i)$ that localises all reals in $\prod b^{\mathrm{lc}}_i$ "not depending on $\alpha$" (so $\mathfrak{v}^{\forall}_{b^{\mathrm{lc}}_i,h^{\mathrm{lc}}_i}$ is increased); and for each $\alpha \in S^{\mathrm{al}}_i$ we add a slalom $\varphi_\alpha \in \mathcal{S}^*(a^{\mathrm{al}}_i, h^{\mathrm{al}}_i)$ that anti-localises all reals in $\prod a^{\mathrm{al}}_i$ "not depending on $\alpha$" (so $\mathfrak{v}^{\exists}_{a^{\mathrm{al}}_i,h^{\mathrm{al}}_i}$ is increased). Details are presented in Lemma 8.2 and Lemma 8.4. Here, "not depending on $\alpha$" means that continuous reading (which is discussed in detail in Section 7) takes place without using the index $\alpha$.

For the moment, we are not going to add additional assumptions, but when it comes to actually prove Theorem 1.4, we will assume CH, $|S^*| = \mu = \mu^{\aleph_0}$ and $|S^{\mathrm{lc}}_i| = |S^{\mathrm{al}}_i| = \kappa_i = \kappa_i^{\aleph_0}$ for all $i \in S^{\mathrm{pr}}$. With these conditions, the collection of slaloms $\{\varphi_\alpha \mid \alpha \in S^{\mathrm{lc}}_i\}$ is used to force $\kappa_i \leq \mathfrak{v}^{\forall}_{b^{\mathrm{lc}}_i,h^{\mathrm{lc}}_i}$ and the collection of slaloms $\{\varphi_\alpha \mid \alpha \in S^{\mathrm{al}}_i\}$ is used to force $\kappa_i \leq \mathfrak{v}^{\forall}_{a^{\mathrm{al}}_i,h^{\mathrm{al}}_i}$.

In any case, we are adding an $\omega$-sequence for each $\alpha \in S^*$. As in [GK21], we are going to use different levels for each type $\mathrm{t} \in \{\mathrm{pr},\mathrm{lc},\mathrm{al}\}$ of generic we add to make the construction more intuitive and less complex. For $S^{\mathrm{al}}$, we use a lim sup creature forcing construction, while for $\mathrm{t} \in \{\mathrm{pr},\mathrm{lc}\}$ we use lim inf creature forcing. For this reason, the levels corresponding to pr and lc are divided into sublevels.

**Definition 3.1.** The *levels* (or *heights)* of the forcing construction are defined as follows.

(1) **pr-levels:** A *pr-level* is a natural number $n$, i. e. the set of pr-levels is $\omega$. For each $n < \omega$, we will determine some $0 < \iota^{\mathrm{pr}}_n < \omega$. Each pr-level $n$ is divided into **pr-sublevels** $\mathrm{ht}^{\mathrm{pr}}_n := \{(n,u) \mid u < \iota^{\mathrm{pr}}_n\}$. We refer to the set of pr-sublevels by $\mathrm{ht}^{\mathrm{pr}}$, i. e. $\mathrm{ht}^{\mathrm{pr}} := \bigcup_{n<\omega} \mathrm{ht}^{\mathrm{pr}}_n = \{(n,u) \mid n<\omega,\ u<\iota^{\mathrm{pr}}_n\}$. We often identify $(n,0)$ with $n$.

(2) **lc-levels:** An *lc-level* is $L^{\mathrm{lc}}_n := (n, \iota^{\mathrm{pr}}_n)$ for some $n<\omega$. Each lc-level is divided into **lc-sublevels** $\mathrm{ht}^{\mathrm{lc}}_n := I^*_n = \{(L^{\mathrm{lc}}_n, v) \mid v < \iota^*_n\}$ for some natural number $\iota^*_n$. We refer to the set of lc-sublevels by $\mathrm{ht}^{\mathrm{lc}} := \bigcup_{n<\omega} I^*_n = \{(L^{\mathrm{lc}}_n, u) \mid n<\omega,\ v<\iota^*_n\}$. We often identify $L^{\mathrm{lc}}_n$ with $(L^{\mathrm{lc}}_n, 0)$.

(3) **al-levels:** An *al-level* is a pair $L^{\mathrm{al}}_n := (L^{\mathrm{lc}}_n, \iota^*_n)$ for $n<\omega$. We refer to the set of al-levels by $\mathrm{ht}^{\mathrm{al}} := \{L^{\mathrm{al}}_n \mid n<\omega\}$.

Define the *set of heights* by $\mathrm{ht} := \mathrm{ht}^{\mathrm{pr}} \cup \mathrm{ht}^{\mathrm{lc}} \cup \mathrm{ht}^{\mathrm{al}}$ ordered lexicographically, i. e.:

$$\begin{aligned} \ldots < (n,0) < (n,1) < \ldots < (n, \iota^{\mathrm{pr}}_n - 1) < (L^{\mathrm{lc}}_n, 0) < \ldots < (L^{\mathrm{lc}}_n, \iota^*_n - 1) < L^{\mathrm{al}}_n \\ < (n+1,0) < (n+1,1) < \ldots < (n+1, \iota^{\mathrm{pr}}_{n+1} - 1) < (L^{\mathrm{lc}}_{n+1}, 0) < \ldots \end{aligned}$$

For $n<\omega$, let $\mathrm{ht}_n := \mathrm{ht}^{\mathrm{pr}}_n \cup I^*_n \cup \{L^{\mathrm{al}}_n\}$. Also let

$$\begin{aligned} \mathrm{Ht} &:= \{(n,0) \mid n<\omega\} \cup \{(L^{\mathrm{lc}}_n, 0) \mid n<\omega\} \cup \{L^{\mathrm{al}}_n\} \text{ and} \\ \Omega &:= \mathrm{Ht} \smallsetminus \mathrm{ht}^{\mathrm{pr}}. \end{aligned}$$

Hence, our forcing poset will "live" on the set

$$\mathrm{DOM} := (S^{\mathrm{pr}} \times \mathrm{ht}^{\mathrm{pr}}) \cup (S^{\mathrm{lc}} \times \mathrm{ht}^{\mathrm{lc}}) \cup (S^{\mathrm{al}} \times \mathrm{ht}^{\mathrm{al}})$$

and we will add a generic real $y_\alpha \colon \mathrm{ht}^{\mathrm{t}} \to \omega$ for each $\mathrm{t} \in \{\mathrm{pr}, \mathrm{lc}, \mathrm{al}\}$ and each $\alpha \in S^{\mathrm{t}}$. In the following subsection, we will explain more about the nature and purpose of these generic reals.

We also fix the sequence $\bar{I}^* := \langle I^*_n \mid n<\omega\rangle$, which will be relevant to deal with the lim inf construction increasing localisation cardinals.

The numbers $\iota^{\mathrm{pr}}_n$ and $\iota^*_n$ will be determined in Subsection 3.3 along with many other parameters we need to set for the forcing construction.

3.2. **Parameters to build blocks.** Our aim in this subsection is to define parameters following certain requirements to build the generic blocks added by the forcing. Along with the many parameters determined in Subsection 3.3, we will fix a sequence of natural numbers $\langle T^*_L \mid L \in \mathrm{ht}^{\mathrm{pr}}\rangle$ and let $T^*_n := \prod_{L\in \mathrm{ht}^{\mathrm{pr}}_n} T^*_L$ for all $n<\omega$. For every $i \in S^{\mathrm{pr}}$, we add a generic real $y_i \in \prod_{L\in\mathrm{ht}^{\mathrm{pr}}} T^*_L$, which gives us the real $y^*_i \in \prod_{n<\omega} T^*_n$ defined by $y^*_i(n) := \langle y_i(L) \mid L \in \mathrm{ht}^{\mathrm{pr}}_n\rangle$.[8]

We aim to define blocks for every real $y \in \prod_{n<\omega} T^*_n$. For each $t \in T^*_n$, we will later define functions

$$\begin{aligned} d_t &\colon I^*_n \cup \{L^{\mathrm{al}}_n\} \to \omega, \\ h_t &\colon I^*_n \cup \{L^{\mathrm{al}}_n\} \to \omega, \\ g_t &\colon \{L^{\mathrm{lc}}_n, L^{\mathrm{al}}_n\} \to \omega, \\ b_t &\colon I^*_n \cup \{L^{\mathrm{al}}_n\} \to \omega, \\ a_t &\colon \{L^{\mathrm{lc}}_n, L^{\mathrm{al}}_n\} \to \omega. \end{aligned} \tag{3.2}$$

[8] The reason why $T^*_n$ is decomposed as $\langle T^*_L \mid L \in \mathrm{ht}^{\mathrm{pr}}_L\rangle$ is due to the use of sublevels $\mathrm{ht}^{\mathrm{pr}}_L$ in the lim inf construction of the forcing.

Given these functions, we will define blocks as follows:

**Definition 3.3.** Recall that $\Omega = \mathrm{Ht} \setminus \mathrm{ht}^{\mathrm{pr}}$. For $y \in \prod_{n<\omega} T^*_n$ and $\mathrm{t} \in \{\mathrm{lc}, \mathrm{al}\}$, we make the following definitions:

(1) $d_y \colon \mathrm{ht} \setminus \mathrm{ht}^{\mathrm{pr}} \to \omega$ is defined as $d_y := \bigcup_{n<\omega} d_{y(n)}$; the functions $h_y$ and $b_y$, from $\mathrm{ht} \setminus \mathrm{ht}^{\mathrm{pr}}$ into $\omega$, are defined analogously.
(2) $g_y \colon \Omega \to \omega$ is defined as $g_y := \bigcup_{n<\omega} g_{y(n)}$, and $a_y \colon \Omega \to \omega$ is defined analogously.
(3) $d^{\mathrm{t}}_y \colon \omega \to \omega$ is defined by $d^{\mathrm{t}}_y(n) := d_y(L^{\mathrm{t}}_n)$, and define $g^{\mathrm{t}}_y, a^{\mathrm{t}}_y \in \omega^\omega$ analogously; in the case of $\mathrm{t} = \mathrm{al}$, $h^{\mathrm{al}}$ and $b^{\mathrm{al}}$ are also defined analogously.
(4) $b^{\mathrm{lc}}_y, h^{\mathrm{lc}}_y \colon \omega \to \omega$ are defined by

$$b^{\mathrm{lc}}_y(n) := \prod_{\ell \in I^*_n} b_y(\ell) \text{ and } h^{\mathrm{lc}}_y(n) := b^{\mathrm{lc}}_y(n) - \prod_{\ell \in I^*_n} (b_y(\ell) - h_y(\ell)) \text{ (cf. Lemma 2.14).}$$

(5) $\langle J^{\mathrm{t}}_{y,n} \mid n < \omega \rangle$ is the interval partition of $\omega$ such that $|J^{\mathrm{t}}_{y,n}| = g^{\mathrm{t}}_y(n)$.
(6) $f^{\mathrm{t}}_y \in \omega^\omega$ is defined by

$$f^{\mathrm{t}}_y(k) := \sum_{j \le n} h^{\mathrm{t}}_y(j) \lceil \log_2 b^{\mathrm{t}}_y(j) \rceil + k - \min J^{\mathrm{t}}_{y,n} \text{ for } k \in J^{\mathrm{t}}_{y,n} \text{ (cf. Lemma 2.9).}$$

(7) $\vec{\nu}^{\mathrm{lc}}_y := \langle d^{\mathrm{lc}}_y, h^{\mathrm{lc}}_y, g^{\mathrm{lc}}_y, b^{\mathrm{lc}}_y, f^{\mathrm{lc}}_y, a^{\mathrm{lc}}_y \rangle$ and $\vec{\nu}^{\mathrm{al}}_y := \langle h^{\mathrm{al}}_y, h^{\mathrm{al}}_y, g^{\mathrm{al}}_y, b^{\mathrm{al}}_y, f^{\mathrm{lc}}_y, a^{\mathrm{al}}_y \rangle$.
(8) Moreover, define $b^-_y := b_y{\restriction}_{\mathrm{ht}^{\mathrm{lc}}}$ and $h^-_y := h_y{\restriction}_{\mathrm{ht}^{\mathrm{lc}}}$.

In order to ensure that $\vec{\nu}^{\mathrm{lc}}_y$ and $\vec{\nu}^{\mathrm{al}}_y$ are blocks for each $y \in \prod_{n<\omega} T^*_n$, we will define the functions in Eq. (3.2) as well as new functions $b^*_t \colon \{L^{\mathrm{lc}}_n, L^{\mathrm{al}}_n\} \to \omega$ and $h^*_t \colon \{L^{\mathrm{lc}}_n, L^{\mathrm{al}}_n\} \to \omega$ satisfying the requirements below for each $n < \omega$ and $t \in T^*_n$; we will refer to the blocks added for $i \in S^{\mathrm{pr}}$ by $\vec{\nu}^{\mathrm{lc}}_i := \vec{\nu}^{\mathrm{lc}}_{y_i}$ and $\vec{\nu}^{\mathrm{al}}_i := \vec{\nu}^{\mathrm{al}}_{y_i}$.

(fp1) $d_t(L^{\mathrm{t}}_n) < h_t(L^{\mathrm{t}}_n)$ for $\mathrm{t} \in \{\mathrm{lc}, \mathrm{al}\}$;
(fp2) $g_t(L^{\mathrm{lc}}_n) \ge \left(\sum\{d_t(L) \mid L \in \Omega,\ L \le L^{\mathrm{lc}}_n\}\right)^{\max\{n,1\}}$;
(fp3) $g_t(L^{\mathrm{al}}_n) \ge \left(\sum\{h_t(L) \mid L \in \Omega,\ L \le L^{\mathrm{al}}_n\}\right)^{\max\{n,1\}}$;
(fp4) $1 - \frac{1}{g_t(L^{\mathrm{lc}}_n)} < \prod_{\ell \in I^*_n} \left(1 - \frac{h_t(\ell)}{b_t(\ell)}\right)$;
(fp5) $h_t(\ell) < b_t(\ell)$ for each $\ell \in I^*_n$;
(fp6) $b^*_t(L^{\mathrm{lc}}_n) := \prod_{\ell \in I^*_n} b_t(\ell)$ and $h^*_t(L^{\mathrm{lc}}_n) := b^*_t(L^{\mathrm{lc}}_n) - \prod_{\ell \in I^*_n}(b_t(\ell) - h_t(\ell))$;
(fp7) $b^*_t(L^{\mathrm{al}}_n) := b_t(L^{\mathrm{al}}_n)$ and $h^*_t(L^{\mathrm{al}}_n) := h_t(L^{\mathrm{al}}_n)$;
(fp8) $g_t(L^{\mathrm{al}}_n) h^*_t(L^{\mathrm{al}}_n) < b^*_t(L^{\mathrm{al}}_n)$;
(fp9) $\log_2 a_t(L^{\mathrm{t}}_n) \ge \sum\{h^*_t(L) \lceil \log_2 b^*_t(L) \rceil \mid L \in \Omega,\ L \le L^{\mathrm{t}}_n\} + g_t(L^{\mathrm{t}}_n)$ for $\mathrm{t} \in \{\mathrm{lc}, \mathrm{al}\}$;
(fp10) $2 < d_t(L^{\mathrm{lc}}_n)$, $2^{a_t(L^{\mathrm{lc}}_n)} < d_t(L^{\mathrm{al}}_n)$ and $2^{a_t(L^{\mathrm{al}}_n)} < d_{t'}(L^{\mathrm{lc}}_{n+1})$ for any $t' \in T^*_{n+1}$.

**Lemma 3.4.** *The requirements (fp1)–(fp10) imply that, for each $y \in \prod_{n<\omega} T^*_n$, the sequences $\vec{\nu}^{\mathrm{lc}}_y$ and $\vec{\nu}^{\mathrm{al}}_y$ are blocks. Moreover, $\mathfrak{v}^{\forall, \bar{I}^*}_{b^-_y, h^-_y} \le \mathfrak{v}^{\forall}_{b^{\mathrm{lc}}_y, h^{\mathrm{lc}}_y}$.*

*Proof.* For $\mathrm{t} \in \{\mathrm{lc}, \mathrm{al}\}$, let $\langle I^{\mathrm{t}}_{y,n} \mid n < \omega \rangle$ be the interval partition of $\omega$ such that $|I^{\mathrm{lc}}_{y,n}| = d^{\mathrm{lc}}_y(n)$ and $|I^{\mathrm{al}}_{y,n}| = h^{\mathrm{al}}_y(n)$. We check that all functions are increasing and that the requirements (i)–(v) of Definition 2.11 hold. Except for $f^{\mathrm{t}}_y$, it is clear that all functions are increasing (also using (fp10)). To check that $f^{\mathrm{t}}_y$ is increasing, as well, it suffices to check that $f^{\mathrm{t}}_y(k) < f^{\mathrm{t}}_y(k+1)$ when $k = \max J^{\mathrm{t}}_{y,n}$. Indeed,

$$f^{\mathrm{t}}_y(k) = \sum_{j\le n} h^{\mathrm{t}}_y(j)\lceil \log_2 b^{\mathrm{t}}_y(j)\rceil + g^{\mathrm{t}}_y(n) - 1 < \sum_{j\le n+1} h^{\mathrm{t}}_y(j)\lceil \log_2 b^{\mathrm{t}}_y(j)\rceil = f^{\mathrm{t}}_t(k+1).$$

(For the inequality, we use $2^{g^{\mathrm{t}}_y(n)} < b^{\mathrm{t}}_y(n+1)$, which follows from (fp10).)

Property (i) follows by (fp2) and (fp3). Property (ii) follows by (fp4) and (fp8), noting that (fp4) implies the inequality in (fp8) when $L^{\mathrm{al}}_n$ is replaced by $L^{\mathrm{lc}}_n$. Property (iii) is obvious by the definition of $f^{\mathrm{t}}_y$, cf. Definition 3.3 (5)–(6). To check property (iv), note that $k \in I^{\mathrm{t}}_{y,n}$ implies $k^n \le \max J^{\mathrm{t}}_{y,n}$ by (i), so $k^n \in J^{\mathrm{t}}_{y,m}$ for some $m \le n$ and

$$\begin{aligned} f^{\mathrm{t}}_y(k^n) &= \sum_{j\le m} h^{\mathrm{t}}_y(j)\lceil \log_2 b^{\mathrm{t}}_y(j)\rceil + k^n - \min J^{\mathrm{t}}_{y,m}\\ &\le \sum_{j\le m} h^{\mathrm{t}}_y(j)\lceil \log_2 b^{\mathrm{t}}_y(j)\rceil + g^{\mathrm{t}}_y(m)\\ &\le \sum_{j\le n} h^{\mathrm{t}}_y(j)\lceil \log_2 b^{\mathrm{t}}_y(j)\rceil + g^{\mathrm{t}}_y(n) \le \lfloor \log_2 a^{\mathrm{t}}_y(n)\rfloor, \end{aligned}$$

where the last inequality follows from (fp9). To check property (v), it is clear from (fp9) that $a^{\mathrm{t}}_y(n) \ge b^{\mathrm{t}}_y(n)$; by (fp1) it is obvious that $d^{\mathrm{al}}_y(n) < h^{\mathrm{al}}_y(n)$, and $d^{\mathrm{lc}}_y(n) < h^{\mathrm{lc}}_y(n)$ follows by (using (fp5))

$$\prod_{\ell\in I^*_n}\left(1 - \frac{h_y(\ell)}{b_y(\ell)}\right) \le 1 - \frac{h_y(L^{\mathrm{lc}}_n)}{b_y(L^{\mathrm{lc}}_n)} < 1 - \frac{d_y(L^{\mathrm{lc}}_n)}{b^{\mathrm{lc}}_y(L^{\mathrm{lc}}_n)} = 1 - \frac{d^{\mathrm{lc}}_y(n)}{b^{\mathrm{lc}}_y(n)}.$$

Finally, from (fp6) and Lemma 2.14 we obtain $\mathbf{Lc}(b^{\mathrm{lc}}_y, h^{\mathrm{lc}}_y) \preceq_{\mathrm{T}} \mathbf{Lc}_{\bar I^*}(b^-_y, h^-_y)$, so $\mathfrak{v}^{\forall,\bar I^*}_{b^-_y,h^-_y} \le \mathfrak{v}^{\forall}_{b^{\mathrm{lc}}_y,h^{\mathrm{lc}}_y}$ follows. □

Recall that, for each $\mathrm{t} \in \{\mathrm{pr}, \mathrm{lc}, \mathrm{al}\}$, we add a generic real $y_\alpha\colon \mathrm{ht}^{\mathrm{t}} \to \omega$. Fix $i \in S^{\mathrm{pr}}$; we have already explained the nature of $y_i$. For $\alpha \in S^{\mathrm{lc}}_i$, $y_\alpha \in \mathcal{S}^*(b^-_{y^*_i}, h^-_{y^*_i})$ will be a generic slalom such that $x \in^*_{\bar I^*} y_\alpha$ for any real $x \in \prod b^-_{y^*_i}$ in the generic extension that does not depend on $\alpha$, so $\mathfrak{v}^{\forall,\bar I^*}_{b^-_{y^*_i},h^-_{y^*_i}}$ is increased, and so is $\mathfrak{v}^{\forall}_{b^{\mathrm{lc}}_{y^*_i},h^{\mathrm{lc}}_{y^*_i}}$ by the previous lemma. In more detail, as in Lemma 2.14 and its proof, setting $\varphi_\alpha(n) := b^{\mathrm{lc}}_{y^*_i}(n) \smallsetminus \prod_{\ell\in I^*_n}(b_{y^*_i}(\ell) \smallsetminus y_\alpha(\ell))$, we get that $\varphi_\alpha \in \mathcal{S}^*(b^{\mathrm{lc}}_{y^*_i}, h^{\mathrm{lc}}_{y^*_i})$ localises all reals in $\prod b^{\mathrm{lc}}_{y^*_i}$ not depending on $\alpha$.

For $\alpha \in S^{\mathrm{al}}_i$, we will have that $y_\alpha(L) \in [a_{y^*_i}(L)]^{\le h_{y^*_i}(L)}$ for all but finitely many $L \in \mathrm{ht}^{\mathrm{al}}$. Setting $\varphi_\alpha(n) := y_\alpha(L^{\mathrm{al}}_n)$, we will get that $\varphi_\alpha \in \mathcal{S}^*(a^{\mathrm{al}}_{y^*_i}, h^{\mathrm{al}}_{y^*_i})$ and that $x \in^\infty \varphi_\alpha$ for any $x \in \prod a^{\mathrm{al}}_{y^*_i}$ not depending on $\alpha$. (A formal description of these generic reals follows in Definition 5.24 and Lemma 5.26.)

3.3. **More parameters for our forcing construction.** In addition to the functions in Eq. (3.2) and the numbers $\iota^{\mathrm{pr}}_n$, $\iota^*_n$ and $T^*_L$ (for all $n < \omega$ and $L \in \mathrm{ht}^{\mathrm{pr}}$), our forcing construction requires several further parameters to ensure properties of the forcing (like properness or continuous reading of names) and to force the desired values in Theorem 1.4.

In this subsection, we list all the parameters and the conditions they must fulfil; the reader may not clearly recognise the intuition and motivation at this point of the paper, however. Although we will briefly state the reason for each condition on the parameters, we suggest it would be more helpful for the reader to go through this subsection quite quickly at first, concentrate on the content that follows, and come back here any time the conditions on these parameters are actually used. This will likely facilitate easier comprehension of our construction's ideas.

**Notation 3.5.** We introduce the notation of all the parameters we will use for the forcing construction. They all are integers larger than 2.

(I) $n^P_{<L}$ and $n^B_L$ for $L \in \mathrm{ht}^{\mathrm{pr}} \cup \Omega$; $P$ stands for "possibilities" and $B$ for "bigness". ($n^P_{<L}$ will be an upper bound for the number of possibilities below $L$ of a modest condition in the forcing, cf. Lemma 5.16; and $n^B_L$ will be a lower bound for the bigness satisfied by a creature at level $L$, cf. Corollary 6.4, Lemma 6.2 and Lemma 6.5.)

(II) $\iota^{\mathrm{pr}}_n$ and $\iota^*_n$ for $n < \omega$.

(III) $T^*_L$ for $L \in \mathrm{ht}^{\mathrm{pr}}$. We also set $T^*_n := \prod_{L \in \mathrm{ht}^{\mathrm{pr}}_n} T^*_L$ with the lexicographic order $<_{\mathrm{lex}}$.

(IV) The functions in Equation 3.2.

(V) $n^S_{t,\ell}$ for each $t \in T^*_n$, $\ell \in I^*_n \cup \{L^{\mathrm{al}}_n\}$ and $n < \omega$. (This will be an upper bound for the number of possibilities of a subatomic creature at $(t, \ell)$, cf. Fact 4.6.)

(VI) $n^S_L$ for all $L \in \mathrm{Ht}$. (This will be an upper bound for the number of possibilities of a creature at level $L$, cf. Fact 4.23.)

We further use the following terminology: Whenever $W = \langle W, \leq_W \rangle$ is a well-ordered set, we write (for $a \in W$)

(1) $a^+$ for the (immediate) successor of $a$ in $W$, if it exists;
(2) $a^-$ for the (immediate) predecessor of $a$ in $W$, if it exists;
(3) $a{\downarrow} := \{x \in W \mid x <_W a\}$;
(4) $\mathrm{ot}(W)$ for the order type of $W$.

(Note that $\mathrm{ot}(\mathrm{ht}) = \omega$.)

In addition, we use:

(5) $t^{\mathrm{mx}}_n := \max_{<_{\mathrm{lex}}} T^*_n$, $t^{\mathrm{mn}}_n := \min_{<_{\mathrm{lex}}} T^*_n$,
(6) $\ell^{\mathrm{mx}}_n := \max I^*_n$.[9]

We also consider $T^*_n \times I^*_n$ with the lexicographic order, i. e. $(t, \ell) < (t', \ell')$ iff either $t <_{\mathrm{lex}} t'$ or $t = t'$ and $\ell < \ell'$.

We present two lists below. The first one indicates the order in which the parameters are defined and the second gives the conditions (pr1)–(pr14) on the parameters.

1. When we have $n^P_{<(n,0)}$, we define $\iota^{\mathrm{pr}}_n$ as in (pr1).
2. For $L \in \mathrm{ht}^{\mathrm{pr}}_n$, given $n^P_{<L}$, we define larger $n^B_L < T^*_L < n^P_{<L^+}$ (in this order), see (pr2)–(pr4).

[9] There is no need to explicitly define $\ell^{\mathrm{mn}}_n := \min I^*_n$, since it clearly is $(L^{\mathrm{lc}}_n, 0) = L^{\mathrm{lc}}_n$.

3. When we get to level $(L_n^{\rm lc},0)$, we define $n^S_{(n,0)}$ as in (pr5) and choose $n^P_{<L_n^{\rm lc}}$ as in (pr14) (which is stronger than (pr4) when $L$ is the predecessor of $L_n^{\rm lc}$ in ht).
4. Afterwards, we define $\iota_n^*$ as in (pr6) and define $n^B_{L_n^{\rm lc}}$ as in (pr2).
5. For each $t\in T_n^*$, we define the parameters $d_t(L_n^{\rm lc}) < g_t(L_n^{\rm lc}) < h_t(L_n^{\rm lc}) < b_t(L_n^{\rm lc}) < n^S_{t,L_n^{\rm lc}}$. For $(t,\ell) > (t_n^{\rm mn}, L_n^{\rm lc})$ in $T_n^*\times I_n^*$, we define $d_t(\ell) < h_t(\ell) < b_t(\ell) < n^S_{t,\ell}$, and get larger parameters for $(t',\ell') > (t,\ell)$ in $T_n^*\times I_n^*$. For $(t,\ell_n^{\rm mx})$, we additionally define $a_t(L_n^{\rm lc})$. See (fp1)–(fp10) and (pr7)–(pr10).
6. Define $n^S_{L_n^{\rm lc}}$ as in (pr11).
7. At $L_n^{\rm al}$, choose $n^P_{<L_n^{\rm al}}$ as in (pr14) and define $n^B_{L_n^{\rm al}} < d_t(L_n^{\rm al}) < h_t(L_n^{\rm al}) < g_t(L_n^{\rm al}) < b_t(L_n^{\rm al}) < a_t(L_n^s) < n^S_{t,L_n^s}$ for $t\in T_n^*$, while getting larger parameters for larger $t'\in T_n^*$. See (fp1)–(fp10), (pr2), (pr7)–(pr8) and (pr12).
8. Define $n^S_{L_n^{\rm al}}$ as in (pr13).
9. Define $n^P_{<(n+1,0)}$ as in (pr14), and repeat this process.

**Conditions on the parameters.** The following are the precise requirements for the parameters:

(pr1) $\iota_n^{\rm pr} \geq 3^{(n+1)\left(2^{(n+1)\cdot n^P_{<(n,0)}}+1\right)}$, which allows us to get $\mu^n(I_n^*) \geq 2^{(n+1)\cdot n^P_{<(n,0)}}+1$ (Definition 4.8) (for modesty in the pr-part and non-trivial conditions, see Lemma 4.17 and Lemma 5.9);

(pr2) $n^B_L > (n^P_{<L})^{n^P_{<L}}$ (for bigness for lim inf compound creatures and rapid reading, see Lemma 6.5 and Lemma 7.7);

(pr3) for $L\in {\rm ht}_n^{\rm pr}$, $T_L^* \geq (n_L^B)^{n_L^B\cdot 2^{(n+1)\cdot n^P_{<(n,0)}}}$, which will imply $\|T_L^*\|_L \geq 2^{(n+1)\cdot n^P_{<(n,0)}}$ (Definition 4.4) (for non-trivial conditions in the pr-part, see Lemma 4.17 and Lemma 5.9);

(pr4) for $L\in {\rm ht}_n^{\rm pr}$, $n^P_{<L^+} > n^P_{<L}\cdot T_L^*$ (for bigness for pr-compound creatures, see Lemma 6.5);

(pr5) $n^S_{(n,0)} \geq T_n^*$ (upper bound for the number of possibilities of a pr-compound creature at $n$, see Fact 4.23);

(pr6) $\iota_n^* \geq 3^{(n\cdot T_n^*+1)\left(2^{(n+1)\cdot n^P_{<L_n^{\rm lc}}}+1\right)}$, which will imply $\mu^{n\cdot T_n^*}(\iota_n^*) \geq 2^{(n+1)\cdot n^P_{<L_n^{\rm lc}}}+1$ (Definition 4.8) (for modesty in the lim inf part and non-trivial conditions, see Lemma 4.16 and Lemma 5.10.);

(pr7) for $\ell\in I_n^*$, $h_t(\ell) \geq d_t(\ell)^{d_t(\ell)\cdot 2^{(n+1)\cdot n^P_{<L_n^{\rm lc}}}}$, which will imply $\|\,{\rm POSS}^t_\ell\,\|^t_\ell \geq 2^{(n+1)\cdot n^P_{<L_n^{\rm lc}}}$ (Definition 4.5) – likewise when $\ell = L_n^{\rm al}$, but replacing $n^P_{<L_n^{\rm lc}}$ by $n^P_{<L_n^{\rm al}}$ (for non-trivial conditions, see Lemma 4.16 and Lemma 5.10);

(pr8) for $t\in T_n^*$ and $\ell\in I_n^*\cup\{L_n^{\rm al}\}$, $n^S_{t,\ell} \geq \left|[b_t(\ell)]^{\leq h_t(\ell)}\right|$ when $\ell\in{\rm ht}^{\rm lc}$, and $n^S_{t,\ell} \geq \left|[a_t(\ell)]^{\leq h_t(\ell)}\right|$ when $\ell\in{\rm ht}^{\rm al}$ – this will imply $n^S_{t,\ell}\geq |\,{\rm POSS}^t_\ell|$ in Fact 4.6 (also used for the proof of $(\dot a_i,\dot d_i)$-bounding in Lemma 8.12);

(pr9) for $t \in T^*_n$ and $\ell \in I^*_n$, if $(t,\ell)$ is not the maximum of $T^*_n \times \mathrm{ht}_{\mathrm{lc},n}$ then $b_t(\ell)^{m^S_{t,\ell}} \le d_{t'}(\ell')$, where $(t',\ell')$ is the successor of $(t,\ell)$ in $T^*_n \times I^*_n$, and

$$m^S_{t,\ell} := \begin{cases} \prod_{\substack{\ell' \in I^*_n \\ \ell' \le \ell}} n^S_{t,\ell'} \times \prod_{\substack{\ell' \in I^*_n \\ \ell' > \ell}} n^S_{t^-,\ell} & \text{if } t >_{\mathrm{lex}} t^{\mathrm{mn}}_n, \\ \prod_{\substack{\ell' \in I^*_n \\ \ell' \le \ell}} n^S_{t,\ell'} & \text{if } t = t^{\mathrm{mn}}_n; \end{cases}$$

recall that $t^-$ denotes the predecessor of $t$ in $T^*_n$ (for bigness for lc-compound creatures, see Lemma 6.3 and Lemma 8.4);

(pr10) for $t <_{\mathrm{lex}} t'$ in $T^*_n$, $n^B_{L^{\mathrm{lc}}_n} < d_t(L^{\mathrm{lc}}_n)$ and $n^P_{<L^{\mathrm{lc}}_n} \cdot a_t(L^{\mathrm{lc}}_n)^{m^S_{t,\ell^{\mathrm{mx}}_n}} < d_{t'}(L^{\mathrm{lc}}_n)$ (for bigness for lc-compound creatures and $(\dot a_i, \dot d_i)$-bounding, see Corollary 6.4 and Lemma 8.12);

(pr11) $n^S_{L^{\mathrm{lc}}_n} \ge m^S_{t^{\mathrm{mx}}_n,\ell^{\mathrm{mx}}_n}$ (upper bound for the number of possibilities of an lc-compound creature at $L^{\mathrm{lc}}_n$, see Fact 4.23);

(pr12) for $t <_{\mathrm{lex}} t'$ in $T^*_n$, $n^B_{L^{\mathrm{al}}_n} < d_t(L^{\mathrm{al}}_n)$ and $n^P_{<L^{\mathrm{al}}_n} \cdot (n^S_{t,L^{\mathrm{al}}_n} + 1) < d_{t'}(L^{\mathrm{al}}_n)$ (for $(\dot a_i, \dot d_i)$-bounding, see Lemma 8.12);

(pr13) $n^S_{L^{\mathrm{al}}_n} \ge n^S_{t^{\mathrm{mx}}_n,L^{\mathrm{al}}_n}$ (upper bound for the number of possibilities of a subatom at $L^{\mathrm{al}}_n$, see Fact 4.23); and

(pr14) for $L \in \mathrm{Ht}$, $n^P_{<L} > \prod_{L' \in \mathrm{Ht} \cap L\downarrow} n^S_{L'}$ (to ensure that $|\mathrm{poss}(p,{<}L)| < n^P_{<L}$ for any modest condition $p$ and for rapid reading, see Lemma 5.16 and Theorem 7.8).

**Theorem 3.6.** *There is a set of parameters as in Notation 3.5 satisfying all the requirements (fp1)–(fp10) and (pr1)–(pr14).*

*Proof.* Based on the properties (pr1)–(pr14), the construction of the parameters from steps 1–9 is straightforward. However, steps 5 and 7 require further explanation, in particular to guarantee properties (fp1)–(fp10).

Assume $t \in T^*_n$ and that we have taken care of step 5 for all $(t_0,\ell_0) < (t, L^{\mathrm{lc}}_n)$ in $T^*_n \times I^*_n$. Choose $d_t(L^{\mathrm{lc}}_n) > n^B_{L^{\mathrm{lc}}_n}$ such that, when $t > t^{\mathrm{mn}}_n$, $d_t(L^{\mathrm{lc}}_n) > n^P_{<L^{\mathrm{lc}}_n} \cdot a_{t^-}(L^{\mathrm{lc}}_n)^{m^S_{t^-,\ell^{\mathrm{mx}}_n}}$ as in (pr10). Define $h_t(L^{\mathrm{lc}}_n)$ as in (pr7) (so (fp1) holds for t = lc) and define $g_t(L^{\mathrm{lc}}_n)$ as in (fp2). Since $1 - \frac{1}{g_t(L^{\mathrm{lc}}_n)} < 1$, we can find some $b_t(L^{\mathrm{lc}}_n) > h_t(L^{\mathrm{lc}}_n)$ such that $1 - \frac{1}{g_t(L^{\mathrm{lc}}_n)} < 1 - \frac{h_t(L^{\mathrm{lc}}_n)}{b_t(L^{\mathrm{lc}}_n)}$. Define $n^S_{t,L^{\mathrm{lc}}_n}$ as in (pr8).

Now assume that $\ell > L^{\mathrm{lc}}_n$ in $I^*_n$ and that we have defined $d_t(\ell_0)$, $h_t(\ell_0)$, $b_t(\ell_0)$ and $n^S_{t,\ell_0}$ for all $\ell_0 < \ell$ in $I^*_n$, such that $1 - \frac{1}{g_t(L^{\mathrm{lc}}_n)} < \prod_{\ell_0 \in I^*_n \cap \ell\downarrow} \left(1 - \frac{h_t(\ell)}{b_t(\ell)}\right)$. So $m^S_{t,\ell^-}$ can be defined as in (pr9) and we can pick $d_t(\ell) \ge b_t(\ell^-)^{m^S_{t,\ell^-}}$. Now define $h_t(\ell)$ as in (pr7) and pick $b_t(\ell) > h_t(\ell)$ such that $1 - \frac{1}{g_t(L^{\mathrm{lc}}_n)} < \left(1 - \frac{h_t(\ell)}{b_t(\ell)}\right) \prod_{\ell_0 \in I^*_n \cap \ell\downarrow} \left(1 - \frac{h_t(\ell)}{b_t(\ell)}\right)$.

When we get to $\ell = \ell^{\mathrm{mx}}_n$ and $n^S_{t,\ell^{\mathrm{mx}}_n}$ has been defined, we can define $a_t(L^{\mathrm{lc}}_n)$ as in (fp9). Note that (fp4) is guaranteed.

Step 7 is easier to explain: Choose $n^B_{L^{\mathrm{al}}_n}$ as in (pr2), $d_t(L^{\mathrm{al}}_n) > n^B_{L^{\mathrm{al}}_n}$ and, whenever $t > t^{\mathrm{mn}}_n$, $d_t(L^{\mathrm{al}}_n) > n^P_{<L^{\mathrm{al}}_n} \cdot (n^S_{t^-,L^{\mathrm{al}}_n} + 1)$ as in (pr12). Define $h_t(L^{\mathrm{al}}_n)$ as in (pr7), $g_t(L^{\mathrm{al}}_n)$ as in (fp3), $b_t(L^{\mathrm{al}}_n)$ as in (fp8), $a_t(L^{\mathrm{al}}_n)$ as in (fp9), and $n^S_{t,L^{\mathrm{al}}_n}$ as in (pr8).

Note that (pr2) implies (fp10). □

## 4. Creatures Upon Creatures

Now that we have discussed the frame and the parameters, we are ready to introduce the basic components of the forcing, which we call *subatomic creatures*. Following the same idea as in [FGKS17, GK21], subatomic creatures are used to build *atomic creatures* and *compound creatures*.

### 4.1. Subatomic and atomic creatures.

The notion of subatomic creature was introduced in [FGKS17, Definition 2.1.1]. We only need a weaker version thereof.

**Definition 4.1.** Let POSS be a finite non-empty set. A *subatomic family* living on POSS is a tuple $\langle \mathrm{K}, \|\cdot\|_{\mathrm{K}}\rangle$ such that

(K1) $\mathrm{K} \subseteq \mathcal{P}(\mathrm{POSS}) \smallsetminus \{\varnothing\}$ is non-empty, and
(K2) $\|\cdot\| = \|\cdot\|_{\mathrm{K}}$ is a function from K into $[0,\infty)$, called a *norm*,

and satisfying, for any $\mathbf{c} \in \mathrm{K}$,

(K3) if $\mathbf{c} \subseteq \mathbf{d}$, then $\|\mathbf{c}\| \leq \|\mathbf{d}\|$, and
(K4) if $|\mathbf{c}| = 1$, then $\|\mathbf{c}\| = 0$.

The elements of K are called *subatomic creatures*, or *subatoms* for short.

We call $\mathfrak{c} \in \mathrm{K}$ a *trivial subatom* if $|\mathbf{c}| = 1$.

We just write K for the subatomic family $\langle \mathrm{K}, \|\cdot\|\rangle$ when the norm is clear from the context.

Our forcing construction will use two types of subatomic families. The first one will be used to add the blocks of parameters, while the second one will be used to add the generic slaloms.

**Example 4.2.** Fix the following subatomic families.

(1) Let $T \geq 1$ and $m \geq 2$ be natural numbers. We define the subatomic family $\langle \mathbf{P}^m_T, \|\cdot\|^m_T\rangle$ living on $T$ by
   (i) $\mathbf{P}^m_T := \mathcal{P}(T) \smallsetminus \{\varnothing\}$; and
   (ii) for $\mathbf{c} \in \mathbf{P}^m_T$, $\|\mathbf{c}\|^m_T := \frac{1}{m}\log_m(|\mathbf{c}|)$.
   Note that $\|\mathbf{c}\|^m_T \geq z$ iff $|\mathbf{c}| \geq m^{mz}$
(2) Given $c > l \geq 1$ and $m \geq 2$ in $\omega$, we define the subatomic family $\langle \mathbf{S}^m_{c,l}, \|\cdot\|^m_{c,l}\rangle$ living on $\mathrm{POSS}^m_{c,l} := [c]^{\leq l}$ by
   (i) $\mathbf{S}^m_{c,l} := \mathcal{P}(\mathrm{POSS}^m_{c,l}) \smallsetminus \{\varnothing\}$; and
   (ii) for $\mathbf{c} \in \mathbf{S}^m_{c,l}$,
$$\|\mathbf{c}\|^m_{c,l} := \frac{1}{m}\log_m(\|\mathbf{c}\|^{\mathrm{cov}}_c + 1),$$
   where $\|\mathbf{c}\|^{\mathrm{cov}}_c := \max\{k \leq c \mid \forall x \in [c]^{\leq k}\ \exists y \in \mathbf{c}\colon x \subseteq y\}$.
   Note that $\|\mathbf{c}\|^{\mathrm{cov}}_c \geq k$ iff every $x \subseteq c$ of size $\leq k$ is contained in some element of $\mathbf{c}$. Hence $\|\mathbf{c}\|^m_{c,l} \geq z$ iff $\|\mathbf{c}\|^{\mathrm{cov}}_c \geq m^{mz} - 1$, i. e. every $x \subseteq c$ of size $\leq\lfloor m^{mz} - 1\rfloor$ is contained in some element of $\mathbf{c}$.

**Observation 4.3.** For $\mathbf{c} \in \mathbf{S}^m_{c,l}$, note that $\|\mathbf{c}\|^{\mathrm{cov}}_c \geq 1$ iff $\bigcup \mathbf{c} = c$.

In our forcing construction, we use creatures as in Example 4.2 at each level $\ell \in \mathrm{ht}$. The generic set of our forcing construction will choose one trivial subatom for each $(\alpha, \ell) \in \mathrm{DOM}$, which determines the value $y_\alpha(\ell)$ of the generic real at $\alpha$. To determine $y_i(L) \in T^*_L$ for $i \in S^{\mathrm{pr}}$ and $L \in \mathrm{ht}^{\mathrm{pr}}$, we use:

**Definition 4.4.** Fix $L \in \mathrm{ht}^{\mathrm{pr}}$. We define the subatomic family $\langle \mathrm{K}_L, \|\cdot\|_L\rangle := \langle \mathbf{P}^{n^B_L}_{T^*_L}, \|\cdot\|^{n^B_L}_{T^*_L}\rangle$ living on $\mathrm{POSS}_L := T^*_L$ (see Example 4.2 (1)), i.e.

(i) $\mathrm{K}_L := \mathcal{P}(\mathrm{POSS}_L) \setminus \{\varnothing\}$; and
(ii) For $\mathbf{c} \in \mathrm{K}_L$, $\|\mathbf{c}\|_L := \frac{1}{n^B_L} \log_{n^B_L}(|\mathbf{c}|)$.

The subatoms in $\mathrm{K}_L$ are also called *(sub-)atomic creatures at $L$.*

For the generic slaloms, we use the following subatomic families:

**Definition 4.5.** Fix $n < \omega$ and $t \in T^*_n$. Using Example 4.2 (2), we make the following definitions:

(1) For each $\ell \in I^*_n$, we define the subatomic family $\langle \mathrm{K}^t_\ell, \|\cdot\|^t_\ell\rangle := \langle \mathbf{S}^{d_t(\ell)}_{b_t(\ell), h_t(\ell)}, \|\cdot\|^{d_t(\ell)}_{b_t(\ell), h_t(\ell)}\rangle$ living on $\mathrm{POSS}^t_\ell := [b_t(\ell)]^{\le h_t(\ell)}$, i.e.
   (i) $\mathrm{K}^t_\ell := \mathcal{P}(\mathrm{POSS}^t_\ell) \setminus \{\varnothing\}$; and
   (ii) for $\mathbf{c} \in \mathrm{K}^t_\ell$,
$$\|\mathbf{c}\|^t_\ell := \frac{1}{d_t(\ell)} \log_{d_t(\ell)}(\|\mathbf{c}\|^{\mathrm{cov}}_{\ell,t} + 1),$$
   where $\|\mathbf{c}\|^{\mathrm{cov}}_{\ell,t} = \|\mathbf{c}\|^{\mathrm{cov}}_{b_t(\ell)} := \max\{m \le b_t(\ell) \mid \forall x \in [b_t(\ell)]^{\le m}\, \exists y \in \mathbf{c} \colon x \subseteq y\}$.
(2) For $\ell = L^{\mathrm{al}}_n$, we define the subatomic family $\langle \mathrm{K}^t_\ell, \|\cdot\|^t_\ell\rangle := \langle \mathbf{S}^{d_t(\ell)}_{a_t(\ell), h_t(\ell)}, \|\cdot\|^{d_t(\ell)}_{a_t(\ell), h_t(\ell)}\rangle$ living on $\mathrm{POSS}^t_\ell := [a_t(\ell)]^{\le h_t(\ell)}$. We write $\|\mathbf{c}\|^{\mathrm{cov}}_{\ell,t} := \|\mathbf{c}\|^{\mathrm{cov}}_{a_t(\ell)}$ for $\mathbf{c} \in \mathbf{K}^t_\ell$.

The subatoms in $\mathrm{K}^t_\ell$ are also called *subatomic creatures at* $(t, \ell)$ (or just *at $\ell$).*

Our parameters bound the size of a subatom at $(t, \ell)$:

**Fact 4.6.** If $n < \omega$, $\ell \in \mathrm{ht} \setminus \mathrm{ht}^{\mathrm{pr}}$ and $t \in T^*_n$, then $|\,\mathrm{POSS}^t_\ell\,| \le n^S_{t,\ell}$.

*Proof.* Immediate by (pr8). □

Recall that, for $i \in S^{\mathrm{pr}}$ and $\alpha \in S^{\mathrm{lc}}_i$, the generic $y_\alpha$ will be in $\mathcal{S}^*(b^-_{y^*_i}, h^-_{y^*_i})$, so for all but finitely many $n < \omega$ and for all $\ell \in I^*_n$, $y_\alpha(\ell) \subseteq b^-_{y^*_i}(\ell) = b_{y^*_i(n)}(\ell)$ will have size $\le h^-_{y^*_i}(\ell) = h_{y^*_i(n)}(\ell)$. However, the value of $y^*_i(n) \in T^*_n$ is determined by the generic added at $i$. This indicates that the information a forcing contains at $(\alpha, \ell)$ must also consider the possibilities of the value of $y^*_i(n)$, i.e. a subatom in $\mathrm{K}^t_\ell$ for each possible $t \in T^*_n$. A similar situation happens when $(\alpha, \ell) \in S^{\mathrm{al}} \times \mathrm{ht}^{\mathrm{al}}$. This motivates the following notion of atomic creatures:

**Definition 4.7.** Let $n < \omega$ and $\ell \in I^*_n \cup \{L^{\mathrm{al}}_n\}$.

(1) We say that $\mathbf{x}$ is an *atomic creature at $\ell$* if $\mathbf{x} = \langle \mathbf{x}(t) \mid t \in P_{\mathbf{x}}\rangle$, where $P_{\mathbf{x}} \subseteq T^*_n$ is non-empty, and $\mathbf{x}(t)$ is a $\mathrm{K}^t_\ell$-subatom for each $t \in P_{\mathbf{x}}$. Here, $P_{\mathbf{x}}$ is called the *set of pr-indices of* $\mathbf{x}$.

(2) For such an atomic creature, we define the norm $\|\mathbf{x}\|_{\min} := \min\{\|\mathbf{x}(t)\|^t_\ell \mid t \in P_\mathbf{x}\}$.
(3) We say that an atomic creature $\mathbf{x}$ at $\ell$ is *trivial* if $x(t)$ is a trivial subatom for each $t \in P_\mathbf{x}$.
(4) We consider the following partial order of atomic creatures at $\ell$: $\mathbf{y} \leq \mathbf{x}$ iff $P_\mathbf{y} \subseteq P_\mathbf{x}$ and $\mathbf{y}(t) \subseteq \mathbf{x}(t)$ for all $t \in P_\mathbf{y}$.

4.2. **Compound creatures.** We present a method from [FGKS17] (also used in [GK21]) to build compound creatures out of (sub-)atomic creatures. This is essential for constructing the lim inf parts of our forcing. We define yet another type of atomic creature, which uses the following "measure":

**Definition 4.8** ([FGKS17, Definition 2.2.1])**.** For each $m \in \omega$, we define a "measure" $\mu^m$ on finite sets by

$$\mu^m(A) := \frac{\log_3(|A|)}{m+1}.$$

**Definition 4.9** ([FGKS17, Definition 2.2.3])**.** Fix a non-empty finite set $J$. Suppose that $\overline{\mathrm{K}} := \langle \mathrm{K}_j \mid j \in J\rangle$ is a sequence of subatomic families.

(1) We say that $\mathbf{x}$ is a *simple atomic creature along* $J$ if $\mathbf{x} = \langle \mathbf{x}(j) \mid j \in J\rangle$, where each $\mathbf{x}(j)$ is a $\mathrm{K}_j$-subatom.
(2) For $m < \omega$, the *m-norm* of a simple atomic creature $\mathbf{x}$ along $J$ is

$$\|\mathbf{x}\|^m_{\overline{\mathrm{K}}} := \max\left\{\min(\{\|\mathbf{x}(j)\|_{\mathrm{K}_j} \mid j \in A\} \cup \{\mu^m(A)\}) \;\middle|\; A \subseteq J\right\},$$

i. e. the $m$-norm is the maximal $r$ for which there is a set $A \subseteq J$ such that $\mu^m(A) \geq r$ and $\|\mathbf{x}(j)\|_{\mathrm{K}_j} \geq r$ for all $j \in A$. We say that such an $A$ *witnesses the m-norm of* $\mathbf{x}$.
(3) For simple atomic creatures $\mathbf{x}$ and $\mathbf{y}$ along $J$, we write $\mathbf{y} \leq_{\overline{\mathrm{K}}} \mathbf{x}$ if $\mathbf{y}(j) \subseteq \mathbf{x}(j)$ for all $j \in J$. This determines a partial order of simple atomic creatures along $J$. When confusion is unlikely, we just write $\mathbf{y} \leq \mathbf{x}$.
(4) We say that a simple atomic creature $\mathbf{x}$ along $J$ is *trivial* if $\mathbf{x}(j)$ is a trivial subatom for each $j \in J$.

We use these simple atomic creatures to construct compound creatures; we begin with those corresponding to the forcing at $S^{\mathrm{pr}}$. This is a variation of [FGKS17, Definition 2.5.1].

**Definition 4.10.** We say that $\mathbf{c}$ is a *pr-compound creature at level* $n < \omega$ (illustrated on the left-hand side of Figure 4) if it consists of

(1) a countable set $\operatorname{dom}\mathbf{c} \subseteq S^{\mathrm{pr}}$ (possibly empty) and a finite $\operatorname{supp}\mathbf{c} \subseteq \operatorname{dom}\mathbf{c}$;
(2) a rectangle of subatomic creatures $\langle \mathbf{c}(i, L) \mid i \in \operatorname{dom}\mathbf{c}$ and $L \in \mathrm{ht}^{\mathrm{pr}}_n\rangle$ such that each $\mathbf{c}(i) := \langle \mathbf{c}(i, L) \mid L \in \mathrm{ht}^{\mathrm{pr}}_n\rangle$ with $\mathbf{c}(i, L) \in \mathrm{K}_L$ is a simple atomic creature along $\mathrm{ht}^{\mathrm{pr}}_n$; and
(3) a real number $d_\mathbf{c} \geq 0$ which we call the *halving parameter*. We also demand that $\operatorname{supp}\mathbf{c} = \varnothing$ implies $d_\mathbf{c} = 0$.

In addition, $\mathbf{c}$ fulfils:

(4) for each $L \in \mathrm{ht}^{\mathrm{pr}}_n$, there is at most one $i \in \operatorname{supp}\mathbf{c}$ such that $\mathbf{c}(i, L)$ has size $>1$, i. e. it is a non-trivial subatom (we call this property *modesty*);

(5) if $i \in \operatorname{dom} \mathbf{c} \setminus \operatorname{supp} \mathbf{c}$, then $|\mathbf{c}(i, L)| = 1$, i. e. $\mathbf{c}(i, L)$ is a trivial subatom; and

(6) for $i \in \operatorname{supp} \mathbf{c}$, the *stacked norm* of $\mathbf{c}(i)$, denoted by $\|\mathbf{c}(i)\|^{\mathrm{pr}}_{\mathrm{stk}}$, is the $n$-norm from Definition 4.9 (2), i. e. $\|\mathbf{c}(i)\|^{\mathrm{pr}}_{\mathrm{stk}}$ is the maximal $r$ for which there is a set $A \subseteq \mathrm{ht}^{\mathrm{pr}}_n$ with $\mu^n(A) = \frac{\log_3 |A|}{n+1} \geq r$ such that $\|\mathbf{c}(i, L)\|_L \geq r$ for all $L \in A$.

Given two pr-compound creatures $\mathbf{c}$, $\mathbf{c}'$ at level $n$, we write $\mathbf{c}' \leq \mathbf{c}$ iff the following holds:

(i) $\operatorname{dom} \mathbf{c} \subseteq \operatorname{dom} \mathbf{c}'$ and $\operatorname{supp} \mathbf{c}' \cap \operatorname{dom} \mathbf{c} = \operatorname{supp} \mathbf{c}$;

(ii) $\mathbf{c}'(i, L) \subseteq \mathbf{c}(i, L)$ for all $i \in \operatorname{dom} \mathbf{c}$ and $L \in \mathrm{ht}^{\mathrm{pr}}_n$; and

(iii) $d_{\mathbf{c}'} \geq d_{\mathbf{c}}$.

This determines a partial order on the set of pr-compound creatures at $n$.

By allowing empty domains, we can define the *trivial pr-compound creature at n* as $\mathbf{1}_{\mathrm{pr}}(n)$ with empty domain (the halving parameter then must be 0 by (3)). Note that $\mathbf{c} \leq \mathbf{1}_{\mathrm{pr}}(n)$ for any pr-compound creature $\mathbf{c}$ at $n$.

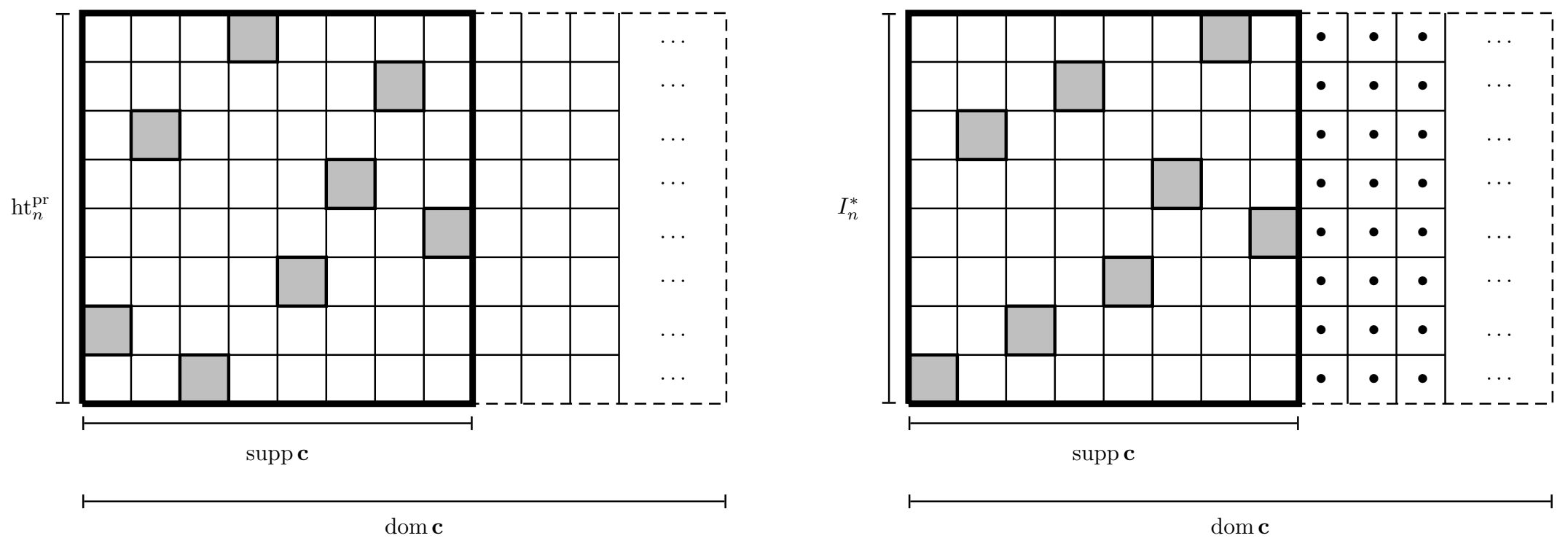

Figure 4. The structure of compound creatures. On the left-hand side, we have a pr-compound creature; the ■ cells denote subatomic creatures that may be non-trivial, while □ cells denote trivial subatomic creatures. On the right-hand side, we have an lc-compound creature; the ■ cells denote atomic creatures (of the form $\mathbf{x} = \langle \mathbf{x}(t) \mid t \in P_{\mathbf{x}} \rangle$ as in Definition 4.7) that may be non-trivial, □ cells denote trivial atomic creatures, and ⊡ cells denote trivial subatomic creatures.

**Definition 4.11.** The *norm* of a pr-compound creature $\mathbf{c}$ at level $n$, when $\sup \mathbf{c} \neq \varnothing$, is defined by

$$\|\mathbf{c}\|^{\mathrm{pr}}_n := \frac{\log_2(\max\{1, \min\{\|\mathbf{c}(i)\|^{\mathrm{pr}}_{\mathrm{stk}} \mid i \in \operatorname{supp} \mathbf{c}\} - d_{\mathbf{c}}\})}{n^P_{<(n,0)}}.$$

In the case that $\operatorname{supp}\mathbf{c} = \varnothing$, we stipulate[10]

$$\|\mathbf{c}\|_n^{\mathrm{pr}} = \begin{cases} 0 & \text{if } \operatorname{dom}\mathbf{c} \neq \varnothing, \\ n & \text{if } \operatorname{dom}\mathbf{c} = \varnothing. \end{cases}$$

In particular, $\|\mathbf{1}_{\mathrm{pr}}(n)\|_n^{\mathrm{pr}} = n$.

Since we require atomic creatures as in Definition 4.7 to add the generic slaloms, the compound creatures corresponding to the lim inf part of the forcing adding slaloms have a more complex structure:

**Definition 4.12.** For $n < \omega$, we say that $\mathbf{c}$ is an *lc-compound creature at level* $L_n^{\mathrm{lc}}$ (illustrated on the right-hand side of Figure 4) if it consists of

(1) a countable set $\operatorname{dom}\mathbf{c} \subseteq S^{\mathrm{lc}}$ (possibly empty) and a finite set $\operatorname{supp}\mathbf{c} \subseteq \operatorname{dom}\mathbf{c}$;
(2) a sequence of *pr-indices* $\bar{P}_{\mathbf{c}} := \langle P_{\mathbf{c},\alpha} \mid \alpha \in \operatorname{supp}\mathbf{c}\rangle$ where each $P_{\mathbf{c},\alpha} \subseteq T_n^*$ is non-empty;
(3) a sequence $\langle \mathbf{c}(\alpha,\ell) \mid \alpha \in \operatorname{dom}\mathbf{c},\ \ell \in I_n^*\rangle$ where
  (i) for $\alpha \in \operatorname{supp}\mathbf{c}$, $\mathbf{c}(\alpha,\ell) = \langle \mathbf{c}(\alpha,t,\ell) \mid t \in P_{\mathbf{c},\alpha}\rangle$ is an atomic creature at $\ell$ and
  (ii) for $\alpha \in \operatorname{dom}\mathbf{c} \setminus \operatorname{supp}\mathbf{c}$, $\mathbf{c}(\alpha,\ell)$ is a trivial subatom in $\mathrm{K}_\ell^{t_n^{\mathrm{mx}}}$; and
(4) a real number $d_{\mathbf{c}} \geq 0$ which we call the *halving parameter*. We also demand that $\operatorname{supp}\mathbf{c} = \varnothing$ implies $d_{\mathbf{c}} = 0$.

In addition, $\mathbf{c}$ fulfils:

(5) for each $\ell \in I_n^*$, there is at most one pair $(\alpha, t) \in \operatorname{supp}\mathbf{c} \times T_n^*$ such that $t \in P_{\mathbf{c},\alpha}$ and $\mathbf{c}(\alpha,t,\ell)$ is non-trivial (we call this property *strong modesty*); and
(6) for $\alpha \in \operatorname{supp}\mathbf{c}$ and $t_0 \in P_{\mathbf{c},\alpha}$, $\mathbf{c}(\alpha,t_0) = \langle \mathbf{c}(\alpha,t_0,\ell) \mid \ell \in I_n^*\rangle$ is considered an atom along $I_n^*$, and its norm, denoted by $\|\mathbf{c}(\alpha,t_0)\|_{\mathrm{stk}}^{t_0}$, is the $n \cdot T_n^*$-norm from Definition 4.9, i. e. $\|\mathbf{c}(\alpha,t_0)\|_{\mathrm{stk}}^{t_0}$ is the maximal $r$ for which there is a set $A \subseteq I_n^*$ with $\frac{\log_3 |A|}{n \cdot T_n^* + 1} \geq r$ such that $\|\mathbf{c}(\alpha,t_0,\ell)\|_\ell^{t_0} \geq r$ for all $\ell \in A$.

Given two lc-compound creatures $\mathbf{c}$ and $\mathbf{c}'$ at $L_n^{\mathrm{lc}}$, we write $\mathbf{c}' \leq \mathbf{c}$ iff the following holds:

(i) $\operatorname{dom}\mathbf{c} \subseteq \operatorname{dom}\mathbf{c}'$ and $\operatorname{supp}\mathbf{c}' \cap \operatorname{dom}\mathbf{c} = \operatorname{supp}\mathbf{c}$;
(ii) for all $\alpha \in \operatorname{supp}\mathbf{c}$, $P_{\mathbf{c}',\alpha} \subseteq P_{\mathbf{c},\alpha}$;
(iii) $\mathbf{c}'(\alpha,t,\ell) \subseteq \mathbf{c}(\alpha,t,\ell)$ for every $\alpha \in \operatorname{supp}\mathbf{c}$, $\ell \in I_n^*$ and $t \in P_{\mathbf{c}',\alpha}$;
(iv) $\mathbf{c}'(\alpha,\ell) = \mathbf{c}(\alpha,\ell)$ for all $\alpha \in \operatorname{dom}\mathbf{c} \setminus \operatorname{supp}\mathbf{c}$; and
(v) $d_{\mathbf{c}'} \geq d_{\mathbf{c}}$.

This determines a partial order of the lc-compound creatures at $L_n^{\mathrm{lc}}$. Items (ii) and (iii) mean that $\mathbf{c}'(\alpha,\ell) \leq \mathbf{c}(\alpha,\ell)$ as atomic creatures at $\ell$ (Definition 4.7) for any $\alpha \in \operatorname{supp}\mathbf{c}$ and $\ell \in I_n^*$.

Denote by $\mathbf{1}_{\mathrm{lc}}(L_n^{\mathrm{lc}})$ the lc-compound creature with empty domain (then $d_{\mathbf{1}_{\mathrm{lc}}(L_n^{\mathrm{lc}})} = 0$ by (4)). Clearly, $\mathbf{c} \leq \mathbf{1}_{\mathrm{lc}}(L_n^{\mathrm{lc}})$ for any lc-compound creature $\mathbf{c}$ at $L_n^{\mathrm{lc}}$.

[10] The requirement that $\|\mathbf{c}_n^{\mathrm{pr}}\| = n$ when $\operatorname{dom}\mathbf{c} = \varnothing$ is to allow the empty set as the trivial condition in the forcing.

**Definition 4.13.** The norm of an lc-compound creature $\mathbf{c}$ at level $L_n^{\mathrm{lc}}$, when $\operatorname{supp}\mathbf{c} \neq \varnothing$, is defined by

$$\|\mathbf{c}\|_n^{\mathrm{lc}} := \frac{\log_2(\max\{1, \min\{\|\mathbf{c}(\alpha, t_0)\|_{\mathrm{stk}}^{t_0} \mid \alpha \in \operatorname{supp}\mathbf{c},\ t_0 \in P_{\mathbf{c},\alpha}\} - d_{\mathbf{c}}\})}{n^P_{<L_n^{\mathrm{lc}}}}.$$

In the case that $\operatorname{supp}\mathbf{c} = \varnothing$, we stipulate

$$\|\mathbf{c}\|_n^{\mathrm{lc}} = \begin{cases} 0 & \text{if } \operatorname{dom}\mathbf{c} \neq \varnothing, \\ n & \text{if } \operatorname{dom}\mathbf{c} = \varnothing. \end{cases}$$

In particular, $\|\mathbf{1}_{\mathrm{lc}}(L_n^{\mathrm{lc}})\|_n^{\mathrm{lc}} = n$.

**Lemma 4.14.** *We state a few simple facts for calculating norms of compound creatures.*

*(a) If $x, \varepsilon \in \mathbb{R}$, $\varepsilon \geq 1$ and $m \geq 2$ in $\omega$, then*

$$\log_m(\max\{1, x - \varepsilon\}) \geq \log_m(\max\{1, x\}) - \varepsilon \log_m(2).$$

*Moreover, when $\varepsilon > 0$ and $\mathbf{e}$ is Napier's constant,*

$$\log_m(\max\{1, x - \varepsilon\}) \geq \log_m(\max\{1, x\}) - \varepsilon \log_m(\mathbf{e}).$$

*(b) Let $n < \omega$ and let $\mathbf{c}$ be an lc-compound creature at $L_n^{\mathrm{lc}}$. For each $\alpha \in \operatorname{supp}\mathbf{c}$ and $t \in P_{\mathbf{c},\alpha}$, assume that $A_{\alpha,t} \subseteq I_n^*$ witnesses $\|\mathbf{c}(\alpha, t)\|_{\mathrm{stk}}^t$. Let $\mathbf{d} \leq \mathbf{c}$ be an lc-compound creature at $L_n^{\mathrm{lc}}$ with the same domain, support, halving parameter and pr-indices, such that, for any $\alpha \in \operatorname{supp}\mathbf{c}$ and $t \in P_{\mathbf{c}_\alpha}$, there is some $B_{\alpha,t} \subseteq A_{\alpha,t}$ such that*
*(i) $|A_{\alpha,t} \smallsetminus B_{\alpha,t}| \leq n^P_{<L_n^{\mathrm{lc}}}$ and*
*(ii) for any $\ell \in B_{\alpha,t}$, $\|\mathbf{d}(\alpha, t, \ell)\|_\ell^t \geq \|\mathbf{c}(\alpha, t, \ell)\|_\ell^t - n^P_{<L_n^{\mathrm{lc}}}$.*
*Then $\|\mathbf{d}\|_n^{\mathrm{lc}} \geq \|\mathbf{c}\|_n^{\mathrm{lc}} - 1$.*

*(c) The analogous statement of (b) for pr-compound creatures also holds.*

*Proof.* We only prove (b). For any $\alpha \in \operatorname{supp}\mathbf{c}$ and $t \in P_{\mathbf{c},\alpha}$,

$$\|\mathbf{c}(\alpha, t)\|_{\mathrm{stk}}^{\mathrm{lc}} - \mu^{n \cdot T_n^*}(B_{\alpha,t}) \leq \mu^{n \cdot T_n^*}(A_{\alpha,t}) - \mu^{n \cdot T_n^*}(B_{\alpha,t}) \leq \frac{n^P_{<L_n^{\mathrm{lc}}}}{n \cdot T_n^* + 1} \leq n^P_{<L_n^{\mathrm{lc}}}$$

by (i), so $\mu^{n \cdot T_n^*}(B_{\alpha,t}) \geq \|\mathbf{c}(\alpha, t)\|_{\mathrm{stk}}^{\mathrm{lc}} - n^P_{<L_n^{\mathrm{lc}}}$. By (ii),

$$\|\mathbf{d}(\alpha, t, \ell)\|_\ell^t \geq \|\mathbf{c}(\alpha, t, \ell)\|_\ell^t - n^P_{<L_n^{\mathrm{lc}}} \geq \|\mathbf{c}(\alpha, t)\|_{\mathrm{stk}}^{\mathrm{lc}} - n^P_{<L_n^{\mathrm{lc}}}$$

for any $\ell \in B_{\alpha,t}$, so $\|\mathbf{d}(\alpha, t)\|_{\mathrm{stk}}^{\mathrm{lc}} \geq \|\mathbf{c}(\alpha, t)\|_{\mathrm{stk}}^{\mathrm{lc}} - n^P_{<L_n^{\mathrm{lc}}}$. Therefore

$$\begin{aligned}\|\mathbf{d}\|_n^{\mathrm{lc}} &\geq \frac{\log_2(\max\{1, \min\{\|\mathbf{c}(\alpha, t_0)\|_{\mathrm{stk}}^{t_0} \mid \alpha \in \operatorname{supp}\mathbf{c},\ t_0 \in P_{\mathbf{c},\alpha}\} - n^P_{<L_n^{\mathrm{lc}}} - d_{\mathbf{c}}\})}{n^P_{<L_n^{\mathrm{lc}}}} \\ &\geq \|\mathbf{c}\|_n^{\mathrm{lc}} - \frac{n^P_{<L_n^{\mathrm{lc}}}}{n^P_{<L_n^{\mathrm{lc}}}} = \|\mathbf{c}\|_n^{\mathrm{lc}} - 1,\end{aligned}$$

where the last inequality follows from (a). □

To make sense of the forcing, we need to be able to construct, at least, subatoms, atoms and compound creatures with large norm and (mostly) arbitrary domain and support. At this point, we finally start using all of the rules we have already presented in Subsection 3.3.

We need the following simple fact for the construction of compound creatures.

**Lemma 4.15** ([FGKS17, Lemma 2.2.2])**.** *Let $k \leq m+1$ and let $\langle A_i \mid i \leq k \rangle$ be a sequence of finite sets. Then there is a sequence $\langle B_i \mid i \leq k \rangle$ of pairwise disjoint sets such that $B_i \subseteq A_i$ and $\mu^m(B_i) \geq \mu^m(A_i) - 1$ for any $i \leq k$.*

**Lemma 4.16.** *Fix $B \subseteq S^{\mathrm{lc}}$ countable, $F \subseteq B$ of size $\leq n$, and a sequence $\bar{P} = \langle P_\alpha \mid \alpha \in F \rangle$ of non-empty subsets of $T^*_n$. Further assume it is not the case that $F$ is empty and $B$ is not. Then there is a lc-compound creature $\mathbf{c}$ at $L^{\mathrm{lc}}_n$ with norm $>n$ such that* $\operatorname{dom} \mathbf{c} = B$, $\operatorname{supp} \mathbf{c} = F$ *and* $\bar{P}_{\mathbf{c}} = \bar{P}$.

*Proof.* This is obvious when $F = \varnothing$ (because in this case $B = \varnothing$), so we assume that $F \neq \varnothing$. To find $\mathbf{c}$, we proceed as follows: First let $d_{\mathbf{c}} = 0$, $\operatorname{supp} \mathbf{c} = F$, $\operatorname{dom} \mathbf{c} = B$, and $P_{\mathbf{c},\alpha} = P_\alpha$ for each $\alpha \in F$. For $\alpha \in \operatorname{supp} \mathbf{c}$ and $t \in P_{\mathbf{c},\alpha}$, we could let $\mathbf{c}(\alpha, t, \ell) := \mathrm{POSS}^t_\ell$ for each $\ell \in I^*_n$, but such a $\mathbf{c}$ would not satisfy strong modesty (Definition 4.12 (5)). To get a $\mathbf{c}$ satisfying this property, we need to do a bit more.

Let $A_{\alpha,t} := I^*_n$ for each $\alpha \in F$ and $t \in P_\alpha$, and let $\mathcal{A} = \{A_{\alpha,t} \mid \alpha \in F,\ t \in P_\alpha\}$. Note that $|F| \leq n$ and $|P_\alpha| \leq T^*_n$, so $|\mathcal{A}| \leq n \cdot T^*_n$. By Lemma 4.15, we can find a pairwise disjoint sequence $\langle B_{\alpha,t} \mid \alpha \in F,\ t \in P_\alpha \rangle$ such that $B_{\alpha,t} \subseteq A_{\alpha,t}$ and $\mu^{n \cdot T^*_n}(B_{\alpha,t}) \geq \mu^{n \cdot T^*_n}(A_{\alpha,t}) - 1$ for each $\alpha \in F$ and $t \in P_\alpha$. So by (pr6), $\mu^{n \cdot T^*_n}(B_{\alpha,t}) \geq 2^{(n+1) \cdot n^P_{<L^{\mathrm{lc}}_n}}$.

Finally, we can define $\mathbf{c}$ as follows: For $\alpha \in \operatorname{supp} \mathbf{c}$ and $t \in P_{\mathbf{c},\alpha}$, we let $\mathbf{c}(\alpha, t, \ell) := \mathrm{POSS}^t_\ell$ for each $\ell \in B_{\alpha,t}$; $\mathbf{c}(\alpha, t, \ell)$ is an arbitrary singleton in $\mathrm{K}^t_\ell$ for each $\ell \in I^*_n \setminus B_{\alpha,t}$; and $\mathbf{c}(\alpha, \ell)$ is a trivial subatom in $\mathrm{K}^{t^{\mathrm{mx}}_n}_\ell$ for $\alpha \in \operatorname{dom} \mathbf{c} \setminus \operatorname{supp} \mathbf{c}$ and $\ell \in \mathrm{ht}^{\mathrm{lc}}_n$. It is clear that $\mathbf{c}$ is an lc-compound creature at $L^{\mathrm{lc}}_n$.

It remains to show that $\|\mathbf{c}\|^{\mathrm{lc}}_n > n$. To this end, note that $\|\mathbf{c}(\alpha, t, \ell)\|^t_\ell \geq 2^{(n+1) \cdot n^P_{<L^{\mathrm{lc}}_n}}$ for each $\ell \in B_{\alpha,t}$ by (pr7), so $\|\mathbf{c}(\alpha, t)\|^t_{\mathrm{stk}} \geq 2^{(n+1) \cdot n^P_{<L^{\mathrm{lc}}_n}}$ for $\alpha \in \operatorname{supp} \mathbf{c}$ and $t \in P_{\mathbf{c},\alpha}$. Therefore, $\|\mathbf{c}\|^{\mathrm{lc}}_n \geq n + 1$. □

A similar proof using (pr1) and (pr3) gives us the existence of pr-compound creatures with large norm.

**Lemma 4.17.** *Fix $B \subseteq S^{\mathrm{al}}$ countable and $F \subseteq B$ of size $\leq n$, and assume it is not the case that $F$ is empty and $B$ is not. Then there is a pr-compound creature $\mathbf{c}$ at $n$ with norm $>n$ such that* $\operatorname{dom} \mathbf{c} = B$ *and* $\operatorname{supp} \mathbf{c} = F$.

To understand the compatibility between conditions in the forcing we need to look at the compatibility between compound creatures.

**Lemma 4.18.** *Let $\mathbf{c}_0$ and $\mathbf{c}_1$ be two lc-compound creatures at level $L^{\mathrm{lc}}_n$ such that*

*(i) for all $\alpha \in \operatorname{dom} \mathbf{c}_0 \cap \operatorname{dom} \mathbf{c}_1$ and for all $\ell \in I^*_n$, $\mathbf{c}_0(\alpha, \ell) = \mathbf{c}_1(\alpha, \ell)$;*[11]

[11] Considering the nature of $\mathbf{c}_j(\alpha, \ell)$ depending on whether $\alpha$ is in $\mathbf{c}_j$'s support or not (a sequence of subatoms indexed with some subset of $T^*_n$ or a single trivial subatom, respectively), this condition (i) implies $\operatorname{supp} \mathbf{c}_{1-j} \cap \operatorname{dom} \mathbf{c}_j \subseteq \operatorname{supp} \mathbf{c}_j$ for $j < 2$.

(ii) *whenever* $\operatorname{supp}\mathbf{c}_j \not\subseteq \operatorname{supp}\mathbf{c}_{1-j}$ *for all* $j < 2$, $|\operatorname{supp}\mathbf{c}_0 \cup \operatorname{supp}\mathbf{c}_1| \leq n$ *and* $d_{\mathbf{c}_0} = d_{\mathbf{c}_1}$*; and*

(iii) *whenever* $\operatorname{supp}\mathbf{c}_i \subsetneq \operatorname{supp}\mathbf{c}_{1-j}$ *for some* $j < 2$, $d_{\mathbf{c}_j} \leq d_{\mathbf{c}_{1-j}}$.

*Then there is some lc-compound creature* $\mathbf{e}$ *at* $L_n^{\mathrm{al}}$ *such that* $\operatorname{dom}\mathbf{e} = \operatorname{dom}\mathbf{c}_0 \cup \operatorname{dom}\mathbf{c}_1$, $\operatorname{supp}\mathbf{e} = \operatorname{supp}\mathbf{c}_0 \cup \operatorname{supp}\mathbf{c}_1$, $\mathbf{e} \leq \mathbf{c}_j$ *for* $j < 2$, *and* $\|\mathbf{e}\|_n^{\mathrm{lc}} \geq \min_{j<2}\{\|\mathbf{c}_j\|_n^{\mathrm{lc}}\} - {}^1/n^P_{L_n^{\mathrm{lc}}}$.

*Proof.* To find $\mathbf{e}$, we proceed as follows: We first define $\operatorname{supp}\mathbf{e} := \operatorname{supp}\mathbf{c}_0 \cup \operatorname{supp}\mathbf{c}_1$, $\operatorname{dom}\mathbf{e} := \operatorname{dom}\mathbf{c}_0 \cup \operatorname{dom}\mathbf{c}_1$, $d_{\mathbf{e}} := \max\{d_{\mathbf{c}_0}, d_{\mathbf{c}_1}\}$, and $P_{\mathbf{e},\alpha} := P_{\mathbf{c}_j,\alpha}$ for $\alpha \in \operatorname{supp}\mathbf{c}_j$ and $j < 2$, which is fine because $\mathbf{c}_0$ and $\mathbf{c}_1$ are identical on their common domain. For the same reason, we may define $\mathbf{e}(\alpha, \ell) = \mathbf{c}_j(\alpha, \ell)$ for $j < 2$, $\alpha \in \operatorname{dom}\mathbf{c}_j$ and $\ell \in I_n^*$, which suffices when one support contains the other (in which case $\|\mathbf{e}\|_n^{\mathrm{lc}} \in \{\|\mathbf{c}_0\|_n^{\mathrm{lc}}, \|\mathbf{c}_1\|_n^{\mathrm{lc}}\}$ by (iii)), but otherwise an $\mathbf{e}$ defined this way may not satisfy strong modesty (Definition 4.12 (5)). So we need to do a bit more when neither support contains in the other, as in (ii). Hence we are assuming $|\operatorname{supp}\mathbf{c}_0 \cup \operatorname{supp}\mathbf{c}_1| \leq n$ and $d_{\mathbf{c}_0} = d_{\mathbf{c}_1}$.

Let $A_{\alpha,t} \subseteq I_n^*$ witness the norm $\|\mathbf{c}_j(\alpha, t)\|_{\mathrm{stk}}^t$ for $j < 2$, $\alpha \in \operatorname{supp}\mathbf{c}_j$ and $t \in P_{\mathbf{c}_j,\alpha}$ (which does not depend on $j$ for $\alpha$ in the common support $\operatorname{supp}\mathbf{c}_0 \cap \operatorname{supp}\mathbf{c}_1$). Let $\mathcal{A} := \{A_{\alpha,t} \mid \alpha \in \operatorname{supp}\mathbf{e},\ t \in P_{\mathbf{e},\alpha}\}$. It is clear that $|\mathcal{A}| \leq n \cdot T_n^*$.

By applying Lemma 4.15 to the family $\mathcal{A}$, we obtain a pairwise disjoint family $\mathcal{B} = \{B_{\alpha,t} \mid \alpha \in \operatorname{supp}\mathbf{e},\ t \in P_{\mathbf{e},\alpha}\}$ such that $B_{\alpha,t} \subseteq A_{\alpha,t}$ and $\mu^{n\cdot T_n^*}(B_{\alpha,t}) \geq \mu^{n\cdot T_n^*}(A_{\alpha,t}) - 1$ for $\alpha \in \operatorname{supp}\mathbf{e}$ and $t \in P_{\mathbf{e},\alpha}$.

Finally, we use $\mathcal{B}$ to define $\mathbf{e}$ as follows:

- For $j < 2$, $\alpha \in \operatorname{supp}\mathbf{c}_j$ and $t \in P_{\mathbf{e},\alpha}$, we let $\mathbf{e}(\alpha, t, \ell) := \mathbf{c}_j(\alpha, t, \ell)$ for each $\ell \in B_{\alpha,t}$; $\mathbf{e}(\alpha, t, \ell)$ is a singleton contained in $\mathbf{c}_j(\alpha, t, \ell)$ for each $\ell \in I_n^* \setminus B_{\alpha,t}$.
- $\mathbf{e}(\alpha, \ell) := \mathbf{c}_j(\alpha, \ell)$ for $\alpha \in \operatorname{dom}\mathbf{c}_j \setminus \operatorname{supp}\mathbf{c}_j$ and $\ell \in I_n^*$.

It is clear that $\mathbf{e}$ is an lc-compound creature and $\mathbf{e} \leq \mathbf{c}_j$ for $j < 2$, so it remains to show that $\|\mathbf{e}\|_n^{\mathrm{lc}} \geq \min_{j<2}\{\|\mathbf{c}_j\|_n^{\mathrm{lc}}\} - 1$. It is clear that $\|\mathbf{e}(\alpha, t)\|_{\mathrm{stk}}^t \geq \|\mathbf{c}_j(\alpha, t)\|_{\mathrm{stk}}^t - 1$ for $j < 2$, $\alpha \in \operatorname{supp}\mathbf{c}_j$ and $t \in P_{\mathbf{e},\alpha}$. Hence

$$\begin{aligned}\|\mathbf{e}\|_n^{\mathrm{lc}} &\geq \frac{\log_2(\max\{1, \min\{\|\mathbf{c}_j(\alpha, t_0)\|_{\mathrm{stk}}^{t_0} \mid j < 2,\ \alpha \in \operatorname{supp}\mathbf{c}_j,\ t_0 \in P_{\mathbf{e},\alpha}\} - 1 - d_{\mathbf{e}}\})}{n^P_{<L_n^{\mathrm{lc}}}} \\ &\geq \min_{j<2}\{\|\mathbf{c}_j\|_n^{\mathrm{lc}}\} - \frac{1}{n^P_{<L_n^{\mathrm{lc}}}} \text{ (by Lemma 4.14 (a)).}\end{aligned}$$

□

In a similar (and simpler) way, we can prove the corresponding fact for pr-compound creatures.

**Lemma 4.19.** *For* $j < 2$, *let* $\mathbf{c}_j$ *be a pr-compound creature at level* $n$ *such that:*

(i) *for all* $i \in \operatorname{dom}\mathbf{c}_0 \cap \operatorname{dom}\mathbf{c}_1$ *and for all* $L \in \mathrm{ht}_n^{\mathrm{pr}}$, $\mathbf{c}_0(i, L) = \mathbf{c}_1(i, L)$*;*

(ii) *whenever* $\operatorname{supp}\mathbf{c}_i \not\subseteq \operatorname{supp}\mathbf{c}_{1-j}$ *for all* $j < 2$, $|\operatorname{supp}\mathbf{c}_0 \cup \operatorname{supp}\mathbf{c}_1| \leq n$ *and* $d_{\mathbf{c}_0} = d_{\mathbf{c}_1}$*;*

(iii) *whenever* $\operatorname{supp}\mathbf{c}_i \subsetneq \operatorname{supp}\mathbf{c}_{1-j}$ *for some* $j < 2$, $d_{\mathbf{c}_j} \leq d_{\mathbf{c}_{1-j}}$*; and*

(iv) $\operatorname{supp}\mathbf{c}_{1-j} \cap \operatorname{dom}\mathbf{c}_j \subseteq \operatorname{supp}\mathbf{c}_j$ *for* $j < 2$.

*Then there is some pr-compound creature* $\mathbf{e}$ *at* $n$ *such that* $\operatorname{dom}\mathbf{e} = \operatorname{dom}\mathbf{c}_0 \cup \operatorname{dom}\mathbf{c}_1$, $\operatorname{supp}\mathbf{c} = \operatorname{supp}\mathbf{c}_0 \cup \operatorname{supp}\mathbf{c}_1$, $\mathbf{e} \leq \mathbf{c}_j$ *for each* $j < 2$, *and* $\|\mathbf{e}\|_n^{\mathrm{pr}} \geq \min\{\|\mathbf{c}_0\|_n^{\mathrm{pr}}, \|\mathbf{c}_1\|_n^{\mathrm{pr}}\} - {}^1\!/n^P_{<(n,0)}$.

**Observation 4.20.** In Lemma 4.18 and Lemma 4.19, we get an explicit construction of $\mathbf{e}$ when $\operatorname{supp}\mathbf{c}_j \subseteq \operatorname{supp}\mathbf{c}_{1-j}$ for some $j < 2$. We denote this $\mathbf{e}$ by $\mathbf{c}_0 \wedge \mathbf{c}_1$.

**Notation 4.21.** From now on, when considering the norms of subatoms, atoms, and compound creatures, we will simply write $\|\mathbf{c}\|$ when confusion is unlikely. The precise norm will be clearly determined by the nature of $\mathbf{c}$; for example, if $\mathbf{c}$ is an lc-compound creature at level $L_n^{\mathrm{lc}}$, we abbreviate $\|\mathbf{c}\| = \|\mathbf{c}\|_n^{\mathrm{lc}}$ and even $\|\mathbf{c}(\alpha, t_0)\| = \|\mathbf{c}(\alpha, t_0)\|_{\mathrm{stk}}^{t_0}$.

To conclude this section, we discuss the set of possibilities of a compound creature.

**Definition 4.22.** Given a compound creature $\mathbf{c}$, we define its set of possibilities according to the following cases.

(1) If $\mathbf{c}$ is a compound creature at $n < \omega$, then $\operatorname{poss}\mathbf{c} := \prod\{\mathbf{c}(i, L) \mid i \in \operatorname{dom}\mathbf{c},\ L \in \mathrm{ht}_n^{\mathrm{pr}}\}$.

(2) If $\mathbf{c}$ is a compound creature at $L_n^{\mathrm{lc}}$ and $\bar{t} = \langle t_\alpha \mid \alpha \in \operatorname{supp}\mathbf{c}\rangle \in \prod_{\alpha\in\operatorname{supp}\mathbf{c}} P_{\mathbf{c},\alpha}$, then

$$\begin{aligned}\mathbf{c}(\bar{t}) :=& \langle \mathbf{c}(\alpha, t_\alpha, \ell) \mid \alpha \in \operatorname{supp}\mathbf{c},\ \ell \in I_n^*\rangle \\ &\cup \langle \mathbf{c}(\alpha, \ell) \mid \alpha \in \operatorname{dom}\mathbf{c} \smallsetminus \operatorname{supp}\mathbf{c},\ \ell \in I_n^*\rangle, \\ \operatorname{poss}(\mathbf{c}, \bar{t}) :=& \prod\{\mathbf{c}(\bar{t})(\alpha, \ell) \mid \alpha \in \operatorname{dom}\mathbf{c},\ \ell \in I_n^*\}.\end{aligned}$$

An upper bound for the number of possibilities is given by the parameters of the forcing.

**Fact 4.23.** Let $n < \omega$.

(a) If $\mathbf{c}$ is a compound creature at $n < \omega$, then $|\operatorname{poss}\mathbf{c}| \leq T_n^* \leq n^S_{(n,0)}$.

(b) If $\mathbf{c}$ is a compound creature at $L_n^{\mathrm{lc}}$ and $\bar{t} \in \prod_{\alpha\in\operatorname{supp}\mathbf{c}} P_{\mathbf{c},\alpha}$, then $|\operatorname{poss}(\mathbf{c}, \bar{t})| \leq n^S_{L_n^{\mathrm{lc}}}$.

(c) If $t \in T_n^*$, then $|\operatorname{POSS}^t_{L_n^{\mathrm{al}}}| \leq n^S_{L_n^{\mathrm{al}}}$.

*Proof.* For (a) and (b), if $\mathbf{c}$ has empty support, then clearly the set of possibilities has only one element. Assume that $\mathbf{c}$ has non-empty support.

The modesty of $\mathbf{c}$ is essential for this proof. For (a), we can find a function $f\colon \mathrm{ht}_n^{\mathrm{pr}} \to \operatorname{supp}\mathbf{c}$ such that $f(L)$ is the unique index in $\operatorname{supp}\mathbf{c}$ such that $\mathbf{c}(f(L), L)$ is non-trivial, in case it exists. Hence, using (pr5),

$$|\operatorname{poss}(\mathbf{c})| = \left|\prod_{L\in\mathrm{ht}_n^{\mathrm{pr}}} \mathbf{c}(f(L), L)\right| \leq \prod_{L\in\mathrm{ht}_n^{\mathrm{pr}}} T_L^* = T_n^* \leq n^S_{(n,0)}.$$

For (b), we can find a function $f\colon I_n^* \to \operatorname{supp}\mathbf{c}$ such that $f(\ell)$ is the unique element of $\operatorname{supp}\mathbf{c}$ such that $\mathbf{c}(f(\ell), t_{f(\ell)}, \ell)$ is non-trivial, in case it exists. Hence, using Fact 4.6, (pr9) and (pr11),

$$|\operatorname{poss}(\mathbf{c},\bar{t})| \leq \prod_{\ell\in I_n^*} |\mathbf{c}(f(\ell), t_{f(\ell),\ell})| \leq \prod_{\ell\in I_n^*} |\operatorname{POSS}_\ell^{t_{f(\ell)}}| \leq \prod_{\ell\in I_n^*} n^S_{t_{f(\ell),\ell}}$$
$$\leq \prod_{\ell\in I_n^*} n^S_{t_n^{\mathrm{mx}},\ell} = m^S_{t_n^{\mathrm{mx}},\ell_n^{\mathrm{mx}}} \leq n^S_{L_n^{\mathrm{lc}}}.$$

Property (c) is clear by Fact 4.6 and (pr13). □

## 5. The Forcing Construction

This section is devoted to defining our forcing and proving some of its basic properties. The construction is based on [FGKS17, GK21].

5.1. **The forcing.** The forcing consists of three parts: the *pr-part* (also called *pr-forcing*), the *lc-part* and the *al-part*. As previously mentioned, the pr-part is a lim inf construction adding the blocks of parameters, the lc-part is also lim inf and it increases $\mathfrak{v}^{\forall,\bar{I}^*}_{b_{y_i^*},h_{y_i^*}}$, and the al-part is a lim sup construction that increases $\mathfrak{v}^{\exists}_{a_{y_i^*},h_{y_i^*}}$. While the pr-forcing is independent, the lc-part and al-part are not – they depend substantially on the pr-part.

We start by describing the pr-part. This is similar to the lim inf part of the forcings from [FGKS17, GK21].

**Definition 5.1.** We define the *pr-forcing* $\mathbb{Q}_{\mathrm{pr}}$ as follows:

**Conditions.** A condition $p \in \mathbb{Q}_{\mathrm{pr}}$ consists of:

(1) a trunk length $\operatorname{trl}(p) \in \omega$;
(2) a sequence $\langle p(n) \mid n < \omega\rangle$ where each $p(n)$ is a pr-compound at level $n$; and
(3) a countable set $\operatorname{supp} p \subseteq S^{\mathrm{pr}}$ (possibly empty),

and it satisfies

(4) $\operatorname{dom} p(n) = \operatorname{supp} p$ for all $n < \omega$;
(5) $\langle \operatorname{supp} p(n) \mid n < \omega\rangle$ is non-decreasing;
(6) $\operatorname{supp} p = \bigcup_{n<\omega} \operatorname{supp} p(n)$;
(7) $\operatorname{supp} p(n) = \varnothing$ for all $n < \operatorname{trl}(p)$;
(8) $\lim_{n\to\infty} \frac{|\operatorname{supp} p(n)|}{n} = 0$; and
(9) $\lim_{n\to\infty} \|p(n)\|_n^{\mathrm{pr}} = \infty$.

**Order.** For $p, q \in \mathbb{Q}_{\mathrm{pr}}$, $q \leq p$ iff $\operatorname{trl}(q) \geq \operatorname{trl}(p)$ and $q(n) \leq p(n)$ for all $n < \omega$.

It is clear that this is a partial order on $\mathbb{Q}_{\mathrm{pr}}$ and that $\mathbf{1} = \mathbf{1}_{\mathrm{pr}} := \langle \mathbf{1}_{\mathrm{pr}}(n) \mid n < \omega\rangle$ (with $\operatorname{trl}(\mathbf{1}) := 0$) is a maximal condition in $\mathbb{Q}_{\mathrm{pr}}$.[12]

This forcing adds an array $\langle \{y_i(L)\} \mid i \in S^{\mathrm{pr}},\ L \in \mathrm{ht}^{\mathrm{pr}}\rangle$ of trivial subatomic creatures, where $y_i(L) \in T^*_L$. As discussed in Subsection 3.2, we obtain a generic $y_i \in \prod_{L\in \mathrm{ht}^{\mathrm{pr}}} T^*_L$ for each $i \in S^{\mathrm{pr}}$, and also a real $y_i^* \in \prod_{n<\omega} T^*_n$ defined by $y_i^*(n) := y_i{\restriction}_{\mathrm{ht}_n^{\mathrm{pr}}}$, which defines the generic block $\bar{\nu}_i^{\mathrm{t}} := \bar{\nu}^{\mathrm{t}}_{y_i^*}$ for $\mathrm{t} \in \{\mathrm{lc}, \mathrm{al}\}$ as described in Definition 3.3. Proofs and details are presented below, beginning with Definition 5.24.

[12] We also allow conditions similar to $\mathbf{1}$, but with longer trunks.

Before describing the other components of our forcing construction, we need to fix some terminology, in particular the definition of *possibilities* of the generic reals given by a condition in $\mathbb{Q}_{\mathrm{pr}}$.

**Notation 5.2.** We fix the following terminology:

(1) For $\ell \in \mathrm{ht}$, denote by $n_*(\ell)$ the unique $n < \omega$ such that $\ell \in \mathrm{ht}_n$.
(2) For each $\alpha \in S^{\mathrm{t}}$ ($\mathrm{t} \in \{\mathrm{lc}, \mathrm{al}\}$), there is a unique $i \in S^{\mathrm{pr}}$ such that $\alpha \in S^{\mathrm{t}}_i$, which we denote by $i^*(\alpha)$. For $\alpha \in S^{\mathrm{pr}}$, let $i^*(\alpha) := \alpha$.
(3) A set $B \subseteq S^*$ is *closed* if, for any $\alpha \in B$, $i^*(\alpha) \in B$.

**Definition 5.3.** Fix a condition $p \in \mathbb{Q}_{\mathrm{pr}}$ and $i \in \operatorname{supp} p$.

(1) For each height $\ell \in \mathrm{ht}$, we let
$$\mathrm{poss}(p(i), {<}\ell) := \prod_{\substack{L<\ell \\ L\in \mathrm{ht}^{\mathrm{pr}}}} p(i, L) \quad \text{and} \quad \mathrm{poss}(p, {<}\ell) := \prod_{i\in \operatorname{supp} p} \mathrm{poss}(p(i), {<}\ell).$$
(2) For $A \subseteq \mathrm{ht}^{\mathrm{pr}}$, denote by $\operatorname{poss} p(i, A) := \prod_{L\in A} p(i, L)$ the *set of possibilities on $A$ at $i$*.
(3) Let $(\alpha, \ell) \in \mathrm{DOM} \setminus (S^{\mathrm{pr}} \times \mathrm{ht}^{\mathrm{pr}})$ such that $i^*(\alpha) \in \operatorname{supp} p$. Let
$$\mathrm{pss}(p, \alpha, \ell) := \operatorname{poss} p(i^*(\alpha), \mathrm{ht}^{\mathrm{pr}}_{n_*(\ell)}) = \prod_{L\in \mathrm{ht}^{\mathrm{pr}}_{n_*(\ell)}} p(i^*(\alpha), L).$$

We now add the lc-part to the pr-forcing: The atomic creatures at $(\alpha, \ell) \in S^{\mathrm{lc}} \times I^*_n$ must depend on the pr-part of the condition, concretely, on $\mathrm{pss}(p_{\mathrm{pr}}, \alpha, \ell)$.

**Definition 5.4.** We define the forcing $\mathbb{Q}_{\mathrm{lc}}$ as follows:

**Conditions.** A condition $p = (p_{\mathrm{pr}}, p_{\mathrm{lc}})$ in $\mathbb{Q}_{\mathrm{lc}}$ consists of:

(1) $p_{\mathrm{pr}} \in \mathbb{Q}_{\mathrm{pr}}$;
(2) a trunk length $\mathrm{trl}(p) = \mathrm{trl}(p_{\mathrm{pr}}) \in \omega$;
(3) a countable set $\operatorname{supp} p_{\mathrm{lc}} \subseteq S^{\mathrm{lc}}$ and $\operatorname{supp} p := \operatorname{supp} p_{\mathrm{pr}} \cup \operatorname{supp} p_{\mathrm{lc}} \subseteq S^*$ (which is closed by (8));
(4) a sequence $p_{\mathrm{lc}} := \langle p(L^{\mathrm{lc}}_n) \mid n < \omega\rangle$ such that
   (i) $p(L^{\mathrm{lc}}_n)$ is an lc-compound creature at level $L^{\mathrm{lc}}_n$;
   (ii) $\operatorname{dom} p(L^{\mathrm{lc}}_n) = \operatorname{supp} p_{\mathrm{lc}}$;
   (iii) $\langle \operatorname{supp} p(L^{\mathrm{lc}}_n) \mid n < \omega\rangle$ is non-decreasing;
   (iv) $\operatorname{supp} p_{\mathrm{lc}} = \bigcup_{n<\omega} \operatorname{supp} p(L^{\mathrm{lc}}_n)$; and
   (v) $P_{p(L^{\mathrm{lc}}_n),\alpha} = \mathrm{pss}(p_{\mathrm{pr}}, \alpha, L^{\mathrm{lc}}_n)$ for any $\alpha \in \operatorname{supp} p(L^{\mathrm{lc}}_n)$;

and it satisfies:

(5) $\operatorname{supp} p(L^{\mathrm{lc}}_n) = \varnothing$ and $\operatorname{supp} p_{\mathrm{pr}}(n) = \varnothing$ for all $n < \mathrm{trl}(p)$;
(6) $\lim_{n\to\infty} \frac{|\operatorname{supp} p(L^{\mathrm{lc}}_n)|}{n} = 0$;
(7) $\lim\limits_{n\to\infty} \|p(L^{\mathrm{lc}}_n)\|^{\mathrm{lc}}_n = \infty$; and
(8) if $\alpha \in \operatorname{supp} p(L^{\mathrm{lc}}_n)$, then $i^*(\alpha) \in \operatorname{supp} p_{\mathrm{pr}}(n)$.

**Order.** For $p, q \in \mathbb{Q}_{\mathrm{lc}}$, we say that $q \le p$ iff $\mathrm{trl}(q) \ge \mathrm{trl}(p)$, $q_{\mathrm{pr}} \le p_{\mathrm{pr}}$ and, for all $n < \omega$,

(i) if either $n < \mathrm{trl}(p)$ or $n \ge \mathrm{trl}(q)$, then $q(L^{\mathrm{lc}}_n) \le p(L^{\mathrm{lc}}_n)$,

(ii) if $\mathrm{trl}(p) \leq n < \mathrm{trl}(q)$, then, for all $\alpha \in \operatorname{supp} p_{\mathrm{lc}}$ and $\ell \in I^*_n$
  (a) if $\alpha \in \operatorname{supp} p(L^{\mathrm{lc}}_n)$, then $q(\alpha, \ell) \subseteq p(\alpha, t, \ell)$, where $t$ is the unique element of $\mathrm{pss}(q, \alpha, \ell)$;
  (b) else, if $\alpha \notin \operatorname{supp} p(L^{\mathrm{lc}}_n)$, then $q(\alpha, \ell) = p(\alpha, \ell)$.

It is routine to show that this defines a partial order on $\mathbb{Q}_{\mathrm{lc}}$. Moreover, $\mathbf{1}$ is a maximal condition in $\mathbb{Q}_{\mathrm{lc}}$ (with $\mathbf{1}_{\mathrm{lc}}(L^{\mathrm{lc}}_n)$ as in Definition 4.12).

We define, for $\alpha \in \operatorname{supp} p_{\mathrm{lc}}$,

$$\mathrm{trl}(p, \alpha) := \min\{n < \omega \mid \alpha \in \operatorname{supp} p(L^{\mathrm{lc}}_n)\}.$$

This forcing adds the same array of trivial subatomic creatures added by $\mathbb{Q}_{\mathrm{pr}}$ and an additional array $\langle \mathbf{c}^*_{\alpha,\ell} \mid (\alpha, \ell) \in S^{\mathrm{lc}} \times \mathrm{ht}^{\mathrm{lc}} \rangle$ of trivial atomic creatures $\mathbf{c}^*_{\alpha,\ell} := \langle \mathbf{c}^*_{\alpha,t,\ell} \mid t \in P^*_{\alpha,\ell} \rangle$ with $P^*_{\alpha,\ell} := \{y^*_{i^*(\alpha)}(n_*(\ell)\}$ and $\mathbf{c}^*_{\alpha,t,\ell} := \{y_\alpha(\ell)\}$, where $y_\alpha(\ell) \in [b_{y^*_{i_*(\alpha)}}(\ell)]^{h_{y^*_{i_*(\alpha)}}(\ell)}$ for all but finitely many $\ell \in \mathrm{ht}^{\mathrm{lc}}$. This defines a generic slalom $y_\alpha \in \mathcal{S}^*(b^-_{y^*_{i_*(\alpha)}}, h^-_{y^*_{i_*(\alpha)}})$, which is used to define the slalom $\varphi_\alpha \in \mathcal{S}^*(b^{\mathrm{lc}}_{y^*_{i_*(\alpha)}}, h^{\mathrm{lc}}_{y^*_{i_*(\alpha)}})$. The technical details are explained from Definition 5.24 onwards.

Property (8) is only relevant to force the $(\dot{a}_i, \dot{d}_i)$-bounding property in Subsection 8.2. Until that point, we only use the fact that $\operatorname{supp} p$ is closed.

We now define the forcing for the al-part (for simplicity, without the lc-part at first).

**Definition 5.5.** We define the forcing $\mathbb{Q}_{\mathrm{al}}$ as follows:

**Conditions.** A condition $p = (p_{\mathrm{pr}}, p_{\mathrm{al}}) \in \mathbb{Q}_{\mathrm{al}}$ consists of:

(1) $p_{\mathrm{pr}} \in \mathbb{Q}_{\mathrm{pr}}$;
(2) a trunk length $\mathrm{trl}(p) = \mathrm{trl}(p_{\mathrm{pr}}) \in \omega$;
(3) a countable set $\operatorname{supp} p_{\mathrm{al}} \subseteq S^{\mathrm{al}}$ such that $\operatorname{supp} p := \operatorname{supp} p_{\mathrm{pr}} \cup \operatorname{supp} p_{\mathrm{al}}$ is <u>closed</u>;
(4) each $\alpha \in \operatorname{supp} p_{\mathrm{al}}$ has a lim sup trunk length $\mathrm{trl}(p, \alpha) \geq \mathrm{trl}(p)$ in $\omega$;
(5) a sequence $p_{\mathrm{al}} := \langle p(\alpha, \ell) \mid \alpha \in \operatorname{supp} p_{\mathrm{al}}$ and $\ell \in \mathrm{ht}^{\mathrm{al}} \rangle$ such that, for $\alpha \in \operatorname{supp} p_{\mathrm{al}}$ and $\ell \in \mathrm{ht}^{\mathrm{al}}$,
  (i) if $\ell \geq \mathrm{trl}(p, \alpha)$, then $p(\alpha, \ell) = \langle p(\alpha, t, \ell) \mid t \in \mathrm{pss}(p_{\mathrm{pr}}, \alpha, \ell) \rangle$ is an atomic creature at $\ell$; and
  (ii) if $\ell < \mathrm{trl}(p, \alpha)$, then $p(\alpha, \ell)$ is a trivial subatom in $\mathrm{K}^{t^{\mathrm{mx}}_n}_\ell$;

and it satisfies:

(6) for all $\alpha \in \operatorname{supp} p_{\mathrm{al}}$,
$$\limsup_{\ell \in \mathrm{ht}_{\mathrm{al}}} \|p(\alpha, \ell)\|_{\mathrm{al}} = \infty,$$
where
$$\|p(\alpha, \ell)\|_{\mathrm{al}} := \begin{cases} \|p(\alpha, \ell)\|_{\min} & \text{(as in Definition 4.7) if } \ell \geq \mathrm{trl}(p, \alpha), \\ \|p(\alpha, \ell)\|^{t^{\mathrm{mx}}_{n_*(\ell)}}_\ell & \text{if } \ell < \mathrm{trl}(p, \alpha). \end{cases}$$
Note that, in the second case, $\|p(\alpha, \ell)\|^{t^{\mathrm{mx}}_{n_*(\ell)}}_\ell = 0$ by (5) (ii).

**Order.** For $p, q \in \mathbb{Q}_{\mathrm{al}}$, we say that $q \leq p$ iff $\operatorname{trl}(q) \geq \operatorname{trl}(p)$, $\operatorname{supp} q \supseteq \operatorname{supp} p$, $q_{\mathrm{pr}} \leq p_{\mathrm{pr}}$ and, for all $\alpha \in \operatorname{supp} p_{\mathrm{al}}$ and $\ell \in \mathrm{ht}^{\mathrm{al}}$,

(i) $\operatorname{trl}(q, \alpha) = \max\{\operatorname{trl}(p, \alpha), \operatorname{trl}(q)\}$;
(ii) if $\ell < \operatorname{trl}(p, \alpha)$, then $q(\alpha, \ell) = p(\alpha, \ell)$;
(iii) if $\operatorname{trl}(p, \alpha) \leq \ell < \operatorname{trl}(q, \alpha)$, then $q(\alpha, \ell) \subseteq p(\alpha, t, \ell)$, where $t$ is the unique element of $\operatorname{pss}(q_{\mathrm{pr}}, \alpha, \ell)$;
(iv) if $\ell \geq \operatorname{trl}(q, \alpha)$, then $q(\alpha, \ell) \leq p(\alpha, \ell)$ as atomic creatures at $\ell$ (Definition 4.7).

It is routine to check that this is a partial order on $\mathbb{Q}_{\mathrm{al}}$. Even more, $\mathbf{1}$ is a maximal condition in $\mathbb{Q}_{\mathrm{al}}$ (with $\operatorname{supp}(\mathbf{1})_{\mathrm{al}} = \varnothing$).

This forcing adds the same array of trivial subatomic creatures added by $\mathbb{Q}_{\mathrm{pr}}$ and an additional array $\langle \mathbf{c}^*_{\alpha,\ell} \mid (\alpha, \ell) \in S^{\mathrm{al}} \times \mathrm{ht}^{\mathrm{al}} \rangle$ of trivial atomic creatures $\mathbf{c}^*_{\alpha,\ell} := \langle \mathbf{c}^*_{\alpha,t,\ell} \mid t \in P^*_{\alpha,\ell} \rangle$ with $P^*_{\alpha,\ell} := \{y^*_{i^*(\alpha)}(n_*(\ell)\}$ and $\mathbf{c}^*_{\alpha,t,\ell} := \{y_\alpha(\ell)\}$, where $y_\alpha(\ell) \in [a_{y^*_{i_*(\alpha)}}(\ell)]^{h_{y^*_{i_*(\alpha)}}(\ell)}$ for all but finitely many $\ell \in \mathrm{ht}^{\mathrm{al}}$. This defines, for each $\alpha \in S^{\mathrm{al}}$, a generic slalom $\varphi_\alpha \in \mathcal{S}^*(a^{\mathrm{al}}_{y^*_{i_*(\alpha)}}, h^{\mathrm{al}}_{y^*_{i_*(\alpha)}})$ by $\varphi_\alpha(n) := y_\alpha(L^{\mathrm{al}}_n)$. The technical details are explained from Definition 5.24 onwards.

We can finally put everything together and define the entire forcing we need.

**Definition 5.6.** Let $\mathbb{Q}$ be a poset whose conditions are those $p = (p_{\mathrm{pr}}, p_{\mathrm{lc}}, p_{\mathrm{al}})$ satisfying $(p_{\mathrm{pr}}, p_{\mathrm{al}}) \in \mathbb{Q}_{\mathrm{al}}$ and $(p_{\mathrm{pr}}, p_{\mathrm{lc}}) \in \mathbb{Q}_{\mathrm{lc}}$. The order is defined by $q \leq p$ iff $(q_{\mathrm{pr}}, q_{\mathrm{al}}) \leq (p_{\mathrm{pr}}, p_{\mathrm{al}})$ and $(q_{\mathrm{pr}}, q_{\mathrm{lc}}) \leq (p_{\mathrm{pr}}, p_{\mathrm{lc}})$. This is a partial order with maximal condition $\mathbf{1}$. Figure 5 illustrates the structure of a condition $p \in \mathbb{Q}$. This forcing adds the generic objects described before for $\mathbb{Q}_{\mathrm{pr}}$, $\mathbb{Q}_{\mathrm{lc}}$ and $\mathbb{Q}_{\mathrm{al}}$.

For each $p \in \mathbb{Q}$, we define $\operatorname{supp} p := \operatorname{supp} p_{\mathrm{pr}} \cup \operatorname{supp} p_{\mathrm{lc}} \cup \operatorname{supp} p_{\mathrm{al}}$. This allows us to define restrictions of $\mathbb{Q}$ for closed $B \subseteq S^*$:

(1) $\mathbb{Q}_B := \{p \in \mathbb{Q} \mid \operatorname{supp} p \subseteq B\}$, with the same order as $\mathbb{Q}$.
(2) For $p \in \mathbb{Q}$, $p|_B$ denotes the condition in $\mathbb{Q}_B$ obtained from $p$ in the natural way, by restricting its support (as well as the supports and domains of the pr and lc-compounds) to $B$.[13]

It is unclear whether $\mathbb{Q}_B$ is a complete subposet of $\mathbb{Q}$ in general (the halving parameters are a potential problem, see Observation 5.20), but this can be guaranteed in some cases. For example, $\mathbb{Q}_{\mathrm{pr}} = \mathbb{Q}_{S^{\mathrm{pr}}}$, $\mathbb{Q}_{\mathrm{lc}} = \mathbb{Q}_{S^{\mathrm{pr}} \cup S^{\mathrm{lc}}}$ and $\mathbb{Q}_{\mathrm{al}} = \mathbb{Q}_{S^{\mathrm{pr}} \cup S^{\mathrm{al}}}$ are complete subposets of $\mathbb{Q}$ (see Corollary 5.19 for the general case).

## 5.2. **Basic properties of the forcing.**

We present some basic properties and features of our forcing construction.

**Definition 5.7.** A condition $p \in \mathbb{Q}$ is *modest* if for each $\ell \in \mathrm{ht}$, there is at most one $\alpha \in \operatorname{supp} p$ such that $p(\alpha, \ell)$ is not a trivial (sub-)atom.

By construction, the lim inf part of any condition $p \in \mathbb{Q}$ satisfies *modesty*. Although it is not necessary to demand modesty for $\ell \in \mathrm{ht}^{\mathrm{al}}$, we are nonetheless going to use this property most of the time. (Note that considering only conditions which are modest at all levels gives us an equivalent forcing notion.)

[13] So $\operatorname{trl}(p|_B) = \operatorname{trl}(p)$ even when $B = \varnothing$, and the halving parameter at a level becomes 0 when the support of the compound at that level is disjoint from $B$.

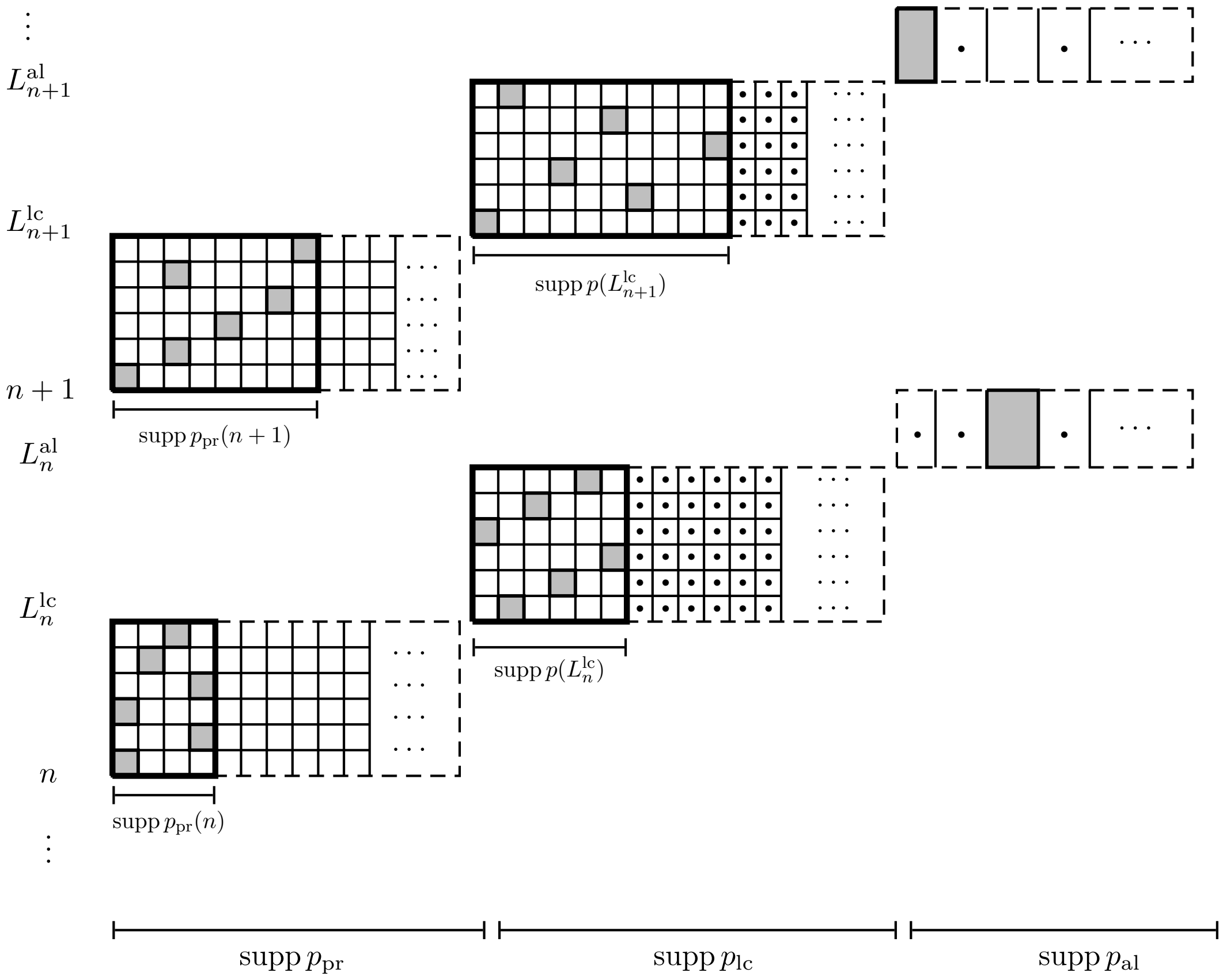

FIGURE 5. The structure of a condition $p \in \mathbb{Q}$. The lim inf compound creatures are as described in Figure 4; the al-part with support $\operatorname{supp} p_{\mathrm{al}}$ is a lim sup construction and, as for the lc-compound creatures, ■ is an atomic creature that may be non-trivial, □ is a trivial atomic creature, and ⊡ denotes a trivial subatomic creature (i. e. a singleton).

**Lemma 5.8.** *Given a condition $p \in \mathbb{Q}$, there is some modest $q \leq p$ in $\mathbb{Q}$ with the same support and trunk lengths and with identical pr- and lc-parts. In particular, the set of modest conditions is dense in $\mathbb{Q}$.*

*Proof.* Let $p$ be a condition in $\mathbb{Q}$. We only consider the case when $\operatorname{supp} p_{\mathrm{al}}$ is non-empty, so we can enumerate it as $\operatorname{supp} p_{\mathrm{al}} =: \{\alpha_m \mid m < \omega\}$, with each element of $\operatorname{supp} p_{\mathrm{al}}$ enumerated infinitely many times.

Define $q$ with the same support and trunk lengths and with identical pr- and lc-parts as $p$ and, by recursion on $m$, define $q_{\mathrm{al}}$ up to some height $\ell_m$. At step $m = 0$, according to Definition 5.5 (6), there is some $\ell_0 \in \mathrm{ht}^{\mathrm{al}}$ such that $\|p(\alpha_0, \ell_0)\|_{\mathrm{al}} > 0$, so let $q(\alpha_0, \ell_0) := p(\alpha_0, \ell_0)$ and choose a trivial atom (or subatom) $q(\alpha_j, \ell') \leq p(\alpha_j, \ell')$ for each $j < \omega$ and $\ell' < \ell_0$ in $\mathrm{ht}^{\mathrm{al}}$, and also for $\alpha_j \neq \alpha_0$ and $\ell' = \ell_0$.

At step $m$, we can find an $\ell_m > \ell_{m-1}$ in $\mathrm{ht}^{\mathrm{al}}$ such that $\|p(\alpha_m, \ell_m)\|_{\mathrm{al}} > m$. So let $q(\alpha_m, \ell_m) := p(\alpha_m, \ell_m)$ and choose a trivial atom $q(\alpha_j, \ell') \leq p(\alpha_j, \ell')$ for each $j < \omega$ and $\ell_{m-1} < \ell' < \ell_m$ in $\mathrm{ht}^{\mathrm{al}}$ and also for $\alpha_j \neq \alpha_m$ and $\ell' = \ell_m$.

The resulting $q$ is as required. □

In the following, we show that we can construct conditions with arbitrary closed support.

**Lemma 5.9.** *For any countable set of indices $B \subseteq S^{\mathrm{pr}}$, there is some $p \in \mathbb{Q}_{\mathrm{pr}}$ with* $\operatorname{supp} p = B$.

*Proof.* If $B = \varnothing$, we already have $p = \mathbf{1} \in \mathbb{Q}_{\mathrm{pr}}$, so assume $B \neq \varnothing$. Enumerate $B \cap S^{\mathrm{pr}} := \{i_j \mid j < \omega\}$ (allowing repetitions). To define the desired $p$, set the trunk length $\mathrm{trl}(p) := 0$ and choose some non-decreasing sequence $\langle x_n \mid n < \omega\rangle$ of natural numbers such that $x_n \leq n$, $\lim_{n\to\infty} x_n = \infty$ and $\lim_{n\to\infty} \frac{x_n}{n} = 0$. Let $\operatorname{supp} p(n) := \{i_j \mid j < x_n\}$. Then Definition 5.1 (8) holds, that is, $\lim_{n\to\infty} \frac{|\operatorname{supp} p(n)|}{n} = 0$. Let $n_0 := \min\{n < \omega \mid x_n > 0\}$. For each $n_0 \leq n < \omega$, by Lemma 4.17, there is a pr-compound $p(n)$ at level $n$ with norm $>n$, support $\operatorname{supp} p(n)$ and domain $B$. For $n < n_0$, let $p(n)$ be any pr-compound with empty support and domain $B$. Then $p \in \mathbb{Q}_{\mathrm{pr}}$ and $\operatorname{supp} p = \bigcup_{n<\omega} \operatorname{supp} p(n) = B$. □

**Lemma 5.10.** *Let $B$ be a closed countable subset of $S^*$ and $p_{\mathrm{pr}} \in \mathbb{Q}_{\mathrm{pr}}$ with support $B \cap S^{\mathrm{pr}}$. Then there is some condition $p$ with support $B$ whose pr-part is $p_{\mathrm{pr}}$.*

*Proof.* We assume that $B \cap S^{\mathrm{lc}} \neq \varnothing$ (otherwise there is no need to go through the "lc-part" of the proof below). Define the condition $p$ as follows:

**Trunk lengths.** Define the trunk lengths $\mathrm{trl}(p) = \mathrm{trl}(p, \alpha) := \mathrm{trl}(p_{\mathrm{pr}})$ for all $\alpha \in B \cap S^{\mathrm{al}}$.

**The pr-part.** As required, we set the pr-part of $p$ as the given $p_{\mathrm{pr}}$.

**The lc-part.** It is important to ensure (8) of Definition 5.4, i. e. for any $\alpha \in \operatorname{supp} p(L^{\mathrm{lc}}_n)$, $i^*(\alpha) \in \operatorname{supp} p_{\mathrm{pr}}(n)$. First, define $B^{\mathrm{lc}}_n := \{\alpha \in B \cap S^{\mathrm{lc}} \mid i^*(\alpha) \in \operatorname{supp} p_{\mathrm{pr}}(n)\}$. Note that $B^{\mathrm{lc}}_n \subseteq B^{\mathrm{lc}}_{n+1}$ and $\bigcup_{n<\omega} B^{\mathrm{lc}}_n = B \cap S^{\mathrm{lc}}$.

Fix $n_0 := \min\{n < \omega \mid B^{\mathrm{lc}}_n \neq \varnothing\}$ (clearly $n_0 \geq \mathrm{trl}(p)$) and a bijection $f\colon \omega \setminus n_0 \to (\omega \setminus n_0) \times \omega$ such that $f(n) = (m, k)$ implies $m \leq n$. Enumerate $B^{\mathrm{lc}}_n := \{\beta_{n,i} \mid i < \omega\}$ (allowing repetitions) for $n \geq n_0$, and define $\gamma^{\mathrm{lc}}_n := \beta_{f(n)}$. So, if $f(n) = (m, k)$, then $m \leq n$ and $\gamma^{\mathrm{lc}}_n \in B^{\mathrm{lc}}_m \subseteq B^{\mathrm{lc}}_n$. Note that $B \cap S^{\mathrm{lc}} = \{\gamma^{\mathrm{lc}}_n \mid n \geq n_0\}$.

As in the proof of Lemma 5.9, choose some non-decreasing sequence $\langle x_n \mid n < \omega\rangle$ of natural numbers such that $x_n \leq n$, $\lim_{n\to\infty} x_n = \infty$ and $\lim_{n\to\infty} \frac{x_n}{n} = 0$. Let $\operatorname{supp} p(L^{\mathrm{lc}}_n) := \varnothing$ for $n < n_0$ and $\operatorname{supp} p(L^{\mathrm{lc}}_n) := \{\gamma^{\mathrm{lc}}_{n_0}, \dots, \gamma^{\mathrm{lc}}_{x_n-1}\}$ for $n \geq n_0$. It is clear that Definition 5.4 (6) and (8) hold. Letting $n_1 := \min\{n < \omega \mid x_n > n_0\}$, by Lemma 4.16, for $n \geq n_1$, there is some lc-compound $p(L^{\mathrm{al}}_n)$ at $L^{\mathrm{lc}}_n$ with norm $>n$, support $\operatorname{supp} p(L^{\mathrm{lc}}_n)$, $\operatorname{dom} p(L^{\mathrm{lc}}_n) = B \cap S^{\mathrm{lc}}$ and $P_{p(L^{\mathrm{lc}}_n),\alpha} = \mathrm{pss}(p_{\mathrm{pr}}, \alpha, L^{\mathrm{lc}}_n)$ for all $\alpha \in \operatorname{supp} p(L^{\mathrm{lc}}_n)$. For $n < n_1$, let $p(L^{\mathrm{lc}}_n)$ be any lc-compound creature at $L^{\mathrm{lc}}_n$ with empty support and domain $B \cap S^{\mathrm{lc}}$. This determines the lc-part $p_{\mathrm{lc}}$ of $p$.

**The al-part.** For $n<\omega$ and $\alpha\in B\cap S^{\mathrm{al}}$, if $n\geq \operatorname{trl}(p)$, let $p(\alpha,t,L^{\mathrm{al}}_n):=\mathrm{POSS}^t_\ell$ for all $t\in \operatorname{pss}(p_{\mathrm{pr}},\alpha,L^{\mathrm{al}}_n)$, otherwise let $p(\alpha,L^{\mathrm{al}}_n)$ be any trivial subatom in $\mathrm{K}^{t^{\mathrm{mx}}_n}_{L^{\mathrm{al}}_n}$. This determines the al-part $p_{\mathrm{al}}$ of $p$, and it is clear that $(p_{\mathrm{pr}},p_{\mathrm{lc}})\in\mathbb{Q}_{\mathrm{lc}}$.

Let $p:=(p_{\mathrm{pr}},p_{\mathrm{lc}},p_{\mathrm{al}})$. To see that $p\in\mathbb{Q}$, it remains to check Definition 5.5 (6). By (pr7), $\|\mathrm{POSS}^t_\ell\|\geq n+1$, so $\|p(\alpha,L^{\mathrm{al}}_n)\|_{\min}\geq n+1$ for each $\alpha\in B\cap S^{\mathrm{al}}$ and $n\geq\operatorname{trl}(p)$, that is, $\|p(\alpha,L^{\mathrm{al}}_n)\|_{\mathrm{al}}\geq n+1$. Therefore, $\limsup_{\ell\in\mathrm{ht}_{\mathrm{al}}}\|p(\alpha,\ell)\|_{\mathrm{al}}=\infty$. Hence the resulting $p$ is indeed a condition in $\mathbb{Q}$ with $\operatorname{supp}p=B$. □

**Corollary 5.11.** *For any closed countable set of indices $B\subseteq S^*$, there is some $p\in\mathbb{Q}$ such that $\operatorname{supp}p=B$. In particular, given any $\alpha\in S^*$, there is some condition $p$ such that $\operatorname{supp}p=\{i^*(\alpha),\alpha\}$.*

We now look at the sets of possibilities of the generic reals determined by a condition in $\mathbb{Q}$.

**Definition 5.12.** Fix a condition $p\in\mathbb{Q}$. For each $n<\omega$ and each function $\eta$ extending an element of $\operatorname{poss}p_{\mathrm{pr}}(n)$, we define, according to Definition 4.22:

(1) for $\alpha\in\operatorname{supp}p$, $\hat\eta(\alpha,n):=\langle\eta(i^*(\alpha),L)\mid L\in\mathrm{ht}^{\mathrm{pr}}_n\rangle$, which is in $T^*_n$; in general, we write $\hat\eta(\alpha,\ell):=\hat\eta(\alpha,n_*(\ell))$ for any $\ell\in\mathrm{ht}$;
(2) $\hat\eta(L^{\mathrm{lc}}_n):=\langle\hat\eta(\alpha,n)\mid\alpha\in\operatorname{supp}p(L^{\mathrm{lc}}_n)\rangle$; and
(3) $p(L^{\mathrm{lc}}_n,\eta):=p(L^{\mathrm{lc}}_n)(\hat\eta(L^{\mathrm{lc}}_n))$.

For $L\in\mathrm{Ht}$, define $\operatorname{poss}(p,{<}L)$ as the set of functions $\eta$ with domain $\mathrm{DOM}\cap(\operatorname{supp}p\times L{\downarrow})$ such that, for any $n<L$ in $\omega$,

(i) $\eta{\restriction}_{S^{\mathrm{pr}}\times\mathrm{ht}^{\mathrm{pr}}_n}\in\operatorname{poss}p_{\mathrm{pr}}(n)$;
(ii) $\eta{\restriction}_{S^{\mathrm{lc}}\times I^*_n}\in\operatorname{poss}p(L^{\mathrm{lc}}_n,\eta)$ when $L^{\mathrm{lc}}_n<L$; and
(iii) when $L^{\mathrm{al}}_n<L$, for $\alpha\in\operatorname{supp}p_{\mathrm{al}}$,
 (a) if $n\geq\operatorname{trl}(p,\alpha)$, then $\eta(\alpha,L^{\mathrm{al}}_n)\in p(\alpha,\hat\eta(\alpha,n),L^{\mathrm{al}}_n)$;
 (b) otherwise, $\eta(\alpha,L^{\mathrm{al}}_n)$ is the unique element of $p(\alpha,L^{\mathrm{al}}_n)$.

When using possibilities $\eta\in\operatorname{poss}(p,{<}L)$, we allow the following abuse of notation:

- for $B\subseteq S^*$, $\eta|_B$ abbreviates $\eta{\restriction}_{B\times\mathrm{ht}}$; and
- for $A\subseteq\mathrm{Ht}$, $\eta{\restriction}_A$ abbreviates $\eta{\restriction}_{S^*\times A}$.

As is to be expected, $\operatorname{poss}(p,{<}L)$ denotes the set of possibilities of the generic on $\mathrm{DOM}\cap(\operatorname{supp}(p)\times L{\downarrow})$.

**Fact 5.13.** Let $L\in\mathrm{Ht}$, $p\in\mathbb{Q}$ and $\eta\in\operatorname{poss}(p,{<}L)$. Then, for each $(\alpha,\ell)\in\mathrm{DOM}$ with $\alpha\in\operatorname{supp}p$ and $\ell<L$:

(a) if $\alpha=i\in S^{\mathrm{pr}}$, then $\eta(i,\ell)\in T^*_\ell$;
(b) if $\alpha\in\operatorname{supp}p(L^{\mathrm{lc}}_{n_*(\ell)})$, then $\eta(\alpha,\ell)\subseteq b_{\hat\eta(\alpha,\ell)}(\ell)$ has size $\leq h_{\hat\eta(\alpha,\ell)}(\ell)$; else, if $\alpha\in S^{\mathrm{lc}}\smallsetminus\operatorname{supp}p(L^{\mathrm{lc}}_{n_*(\ell)})$, then $\eta(\alpha,\ell)\subseteq b_{t^{\mathrm{mx}}_{n_*(\ell)}}(\ell)$ has size $\leq h_{t^{\mathrm{mx}}_{n_*(\ell)}}(\ell)$; and
(c) if $\alpha\in S^{\mathrm{al}}$ and $\ell\geq\operatorname{trl}(p,\alpha)$, then $\eta(\alpha,\ell)\subseteq a_{\hat\eta(\alpha,\ell)}(\ell)$ has size $\leq h_{\hat\eta(\alpha,\ell)}(\ell)$; else, if $\ell<\operatorname{trl}(p,\alpha)$, then $\eta(\alpha,\ell)\subseteq a_{t^{\mathrm{mx}}_{n_*(\ell)}}(\ell)$ has size $\leq h_{t^{\mathrm{mx}}_{n_*(\ell)}}(\ell)$.

For each $\eta\in\operatorname{poss}(p,{<}L)$, we can define a condition $p\wedge\eta\leq p$ forcing the generic below $L$ to be $\eta$. This condition is basically constructed by replacing creatures below $L$ by the subatoms determined by $\eta$ (and increasing trunk lengths, if required).

**Definition 5.14.** Let $p \in \mathbb{Q}$. For every $L \in \mathrm{Ht}$ and $\eta \in \mathrm{poss}(p, {<}L)$, we define a condition $q = p \wedge \eta$ satisfying:

(i) $\operatorname{supp} q = \operatorname{supp} p$;
(ii) $\mathrm{trl}(q) := \max\{\mathrm{trl}(p), n_*(L)\}$, $\mathrm{trl}(q, \alpha) := \max\{\mathrm{trl}(p, \alpha), n_*(L)\}$ for $\alpha \in \operatorname{supp} p_{\mathrm{al}}$;
(iii) any (pr- and lc-) compound creature at a level $\geq \mathrm{trl}(q)$ has the same support, domain, halving parameter and pr-indices as the compound in $p$ at the same level;
(iv) for $i \in \operatorname{supp} p_{\mathrm{pr}}$ and $\ell \in \mathrm{ht}^{\mathrm{pr}}$,
- when $\ell < L$, let $q(i, \ell) := \{\eta(i, \ell)\}$; and
- when $\ell \geq L$, let $q(i, \ell) := p(i, \ell)$;

(v) for $(\alpha, \ell) \in \mathrm{DOM}$ with $\alpha \in \operatorname{supp} q \smallsetminus S^{\mathrm{pr}}$,
- when $\ell < \mathrm{trl}(q)$, let $q(\alpha, \ell) := \{\eta(\alpha, \ell)\}$ for all $\ell < n_*(L)$, otherwise let $q(\alpha, \ell) := p(\alpha, \ell)$;
- when $\mathrm{trl}(q) \leq \ell < \mathrm{trl}(q, \alpha)$, let $q(\alpha, \ell) := p(\alpha, \ell)$ (if there is such an $\ell$, then we must have $\mathrm{trl}(q, \alpha) = \mathrm{trl}(p, \alpha) > n_*(L)$);
- when $\mathrm{trl}(q, \alpha) \leq \ell < L$, let $q(\alpha, \ell) = \langle q(\alpha, t, \ell) \mid t \in \{\hat{\eta}(\alpha, \ell)\}\rangle$, where $q(\alpha, \hat{\eta}(\alpha, \ell), \ell) := \{\eta(\alpha, \ell)\}$ (note that $\mathrm{pss}(q, \alpha, \ell) = \{\hat{\eta}(\alpha, \ell)\}$); and
- when $\ell \geq \max\{L, \mathrm{trl}(q, \alpha)\}$ and $t \in \mathrm{pss}(q_{\mathrm{pr}}, \alpha, \ell)$, let $q(\alpha, t, \ell) := p(\alpha, t, \ell)$.

From the definition above, we get the following simple consequences:

**Fact 5.15.** For $p \in \mathbb{Q}$ and $L \in \mathrm{Ht}$, the following hold:

(a) If $\eta \in \mathrm{poss}(p, {<}L)$, then $p \wedge \eta \in \mathbb{Q}$, and it is stronger than $p$.
(b) $\{p \wedge \eta \leq p \mid \eta \in \mathrm{poss}(p, {<}L)\}$ is an antichain.
(c) For each $q \leq p$ and $\eta' \in \mathrm{poss}(q, {<}L)$, there is a unique $\eta \in \mathrm{poss}(p, {<}L)$ such that $q \wedge \eta' \leq p \wedge \eta$.
(d) $\{p \wedge \eta \leq p \mid \eta \in \mathrm{poss}(p, {<}L)\}$ is a maximal antichain below $p$.

*Proof.* Property (a) is clear by Definition 5.14. For (b), note that two different $\eta, \eta' \in \mathrm{poss}(p, {<}L)$ determine incompatible $p \wedge \eta$ and $p \wedge \eta'$ (if $\eta$ and $\eta'$ differ at some $(\alpha, \ell)$, then the trivial subatoms at this position are disjoint for both conditions).

For (c), it suffices to define $\eta := \eta'|_{\operatorname{supp} p}$. Uniqueness follows by (b).

To check (d), let $q \leq p$ in $\mathbb{Q}$ and choose some $\eta' \in \mathrm{poss}(q, {<}L)$. By (c), there is some $\eta \in \mathrm{poss}(p, {<}L)$ such that $q \wedge \eta' \leq p \wedge \eta$, so $q$ is compatible with $p \wedge \eta$ by (a). □

For our calculations, we use the parameter $n^P_{<L}$ to bound the number of possibilities of a modest condition.

**Lemma 5.16.** *Let $p \in \mathbb{Q}$ be modest and let $L \in \mathrm{Ht}$. Then $|\mathrm{poss}(p, {<}L)| < n^P_L$.*

*Proof.* Since $p$ is modest, for any $\ell \in \mathrm{ht}$ we can either pick the unique point $\alpha_\ell \in \operatorname{supp} p$ such that the creature $p(\alpha_\ell, \ell)$ is non-trivial (in case it exists), and otherwise, we let $\alpha_\ell \in \operatorname{supp} p$ be arbitrary such that $(\alpha_\ell, \ell) \in \mathrm{DOM}$; such an $\alpha_\ell$ may not even exist, in which case the creatures at $\ell$ can be ignored in the calculations below. In what follows, when we write $p(\alpha, t, \ell)$, we ignore $t$ when $\ell < \mathrm{trl}(p, \alpha)$ (i. e. it is interpreted as the trivial subatom $p(\alpha, \ell)$).

In the case that $L \in \mathrm{ht}^{\mathrm{pr}}$, using Fact 4.23 and (pr14),

$$\begin{aligned} |\operatorname{poss}(p, {<}L)| &\leq \Big| \prod_{n\in\omega\cap L\downarrow} \bigcup_{\eta\in\operatorname{poss} p_{\mathrm{pr}}(n)} \{\eta\} \times \operatorname{poss} p(L^{\mathrm{lc}}_n, \hat\eta(L^{\mathrm{lc}}_n)) \times p(\alpha_{L^{\mathrm{al}}_n}, \hat\eta(\alpha, n), L^{\mathrm{al}}_n)\Big| \\ &\leq \prod_{n\in\omega\cap L\downarrow} \sum_{\eta\in\operatorname{poss} p_{\mathrm{pr}}(n)} n^S_{L^{\mathrm{lc}}_n} n^S_{L^{\mathrm{al}}_n} = \prod_{n\in\omega\cap L\downarrow} |\operatorname{poss} p_{\mathrm{pr}}(n)| n^S_{L^{\mathrm{lc}}_n} n^S_{L^{\mathrm{al}}_n} \\ &\leq \prod_{n\in\omega\cap L\downarrow} n^S_{(n,0)} n^S_{L^{\mathrm{lc}}_n} n^S_{L^{\mathrm{al}}_n} = \prod\{n^S_{L'} \mid L' < L,\ L' \in \mathrm{Ht}\} < n^P_{<L}. \end{aligned}$$

In the case that $L = L^{\mathrm{lc}}_n$, we get

$$|\operatorname{poss}(p, {<}L)| = |\operatorname{poss}(p, {<}n) \times \operatorname{poss}(p_{\mathrm{pr}}(n))| \leq \prod\{n^S_{L'} \mid L' < L,\ L' \in \mathrm{Ht}\} < n^P_{<L}.$$

The inequality $|\operatorname{poss}(p, {<}L^{\mathrm{al}}_n)| < n^P_{<L^{\mathrm{al}}_n}$ is checked similarly. □

We now look at several results about compatibility between conditions. We also get several relevant consequences, such as conditions under which $\mathbb{Q}|_B \lessdot \mathbb{Q}$ (i. e. we have complete embeddability) and the fact that $\mathbb{Q}$ has $\mathfrak{c}^+$-cc.

**Lemma 5.17.** *Let $B \subseteq S^*$ be closed and $p \in \mathbb{Q}$. Assume that $r \in \mathbb{Q}_B$, $r \leq p|_B$ and that, for any compound creature in $p$, one of the following conditions holds:*

- *(i) its support is contained in $B$;*
- *(ii) the compound creature in $r$ at the same level has empty support; or*
- *(iii) it has the same halving parameter as the compound creature in $r$ at the same level.*

*Then there is some $q \leq p$ in $\mathbb{Q}$ with support $\operatorname{supp} p \cup \operatorname{supp} r$ such that $q|_B \leq r$. Moreover, if $\mathrm{trl}(r) = \mathrm{trl}(p)$ and (i) or (ii) always hold, then we can even get $q|_B = r$.*

*Proof.* Let $n_0 \geq \mathrm{trl}(r)$ be minimal such that, for any $n \geq n_0$,

- if (i) and (ii) do not hold at level $n$, then $|\operatorname{supp} r_{\mathrm{pr}}(n) \cup \operatorname{supp} p_{\mathrm{pr}}(n)| \leq n$, and
- if (i) and (ii) do not hold at level $L^{\mathrm{lc}}_n$, then $|\operatorname{supp} r(L^{\mathrm{lc}}_n) \cup \operatorname{supp} p(L^{\mathrm{lc}}_n)| \leq n$.

Note that in case (i) or (ii) always hold, $n_0 = \mathrm{trl}(r)$.

Choose some $\eta_0 \in \operatorname{poss}(r, {<}n_0)$ and set $r' := r \wedge \eta_0$, so $\mathrm{trl}(r') = n_0$. Since $r \leq p|_B$, we can find some $\eta_1 \in \operatorname{poss}(p, {<}n_0)$ compatible with $\eta_0$ and set $p' := p \wedge \eta_1$. Note that the hypothesis of the lemma still holds when $p$ and $r$ are replaced by $p'$ and $r'$, respectively. Note that $r' = r$ when $n_0 = \mathrm{trl}(r)$.

We define $q$ as follows.

**Trunk lengths.** $\mathrm{trl}(q) := n_0$, $\mathrm{trl}(q, \alpha) := \mathrm{trl}(r', \alpha)$ for $\alpha \in \operatorname{supp} r'_{\mathrm{al}}$, and $\mathrm{trl}(q, \alpha) := \mathrm{trl}(p', \alpha)$ for $\alpha \in \operatorname{supp} p'_{\mathrm{al}} \smallsetminus B$.

**The pr-part.** Fix $n < \omega$. Let $\mathbf{c}_n$ be the compound creature at $n$ that one gets from $p'_{\rm pr}(n)$ by restricting its domain to $\operatorname{supp} p'_{\rm pr} \smallsetminus B$. Then, by Lemma 4.19, there is some pr-compound creature $q_{\rm pr}(n)$ at $n$ such that

- $\operatorname{dom} q_{\rm pr}(n) = \operatorname{dom} \mathbf{c}_n \cup \operatorname{dom} r'_{\rm pr}(n) = \operatorname{supp} p_{\rm pr} \cup \operatorname{supp} r_{\rm pr}$;
- $\operatorname{supp} q_{\rm pr}(n) = \operatorname{supp} \mathbf{c}_n \cup \operatorname{supp} r'_{\rm pr}(n) = \operatorname{supp} p'_{\rm pr}(n) \cup \operatorname{supp} r'_{\rm pr}(n)$;
- $q_{\rm pr}(n) \leq \mathbf{c}_n$ and $q_{\rm pr}(n) \leq r'_{\rm pr}(n)$; and
- $\|q_{\rm pr}(n)\|^{\rm pr}_n \geq \min\{\|\mathbf{c}_n\|^{\rm pr}_n, \|r'_{\rm pr}(n)\|^{\rm pr}_n\} - 1$, and even
$$\|q_{\rm pr}(n)\|^{\rm pr}_n \geq \min\{\|p'_{\rm pr}(n)\|^{\rm pr}_n, \|r'_{\rm pr}(n)\|^{\rm pr}_n\} - 1$$
(which is obvious when $\mathbf{c}_n$ has non-empty support; otherwise (i) holds and $\|q_{\rm pr}(n)\|^{\rm pr}_n = \|r'_{\rm pr}(n)\|^{\rm pr}_n$).

Consequently, $q_{\rm pr}(n) \leq p_{\rm pr}(n)$. This defines a condition $q_{\rm pr} \in \mathbb{Q}_{\rm pr}$ stronger than both $p'_{\rm pr}$ and $r'_{\rm pr}|_B$. Moreover, if (i) or (ii) hold at level $n$, then $q_{\rm pr}(n) = \mathbf{c}_n \wedge r'_{\rm pr}(n)$ (see Observation 4.20).

**The lc-part.** Fix $n < \omega$. We derive $\mathbf{c}_{n,0}$ from $r'(L^{\rm lc}_n)$ by just reducing the pr-indices to $P_{\mathbf{c}_{n,0},\alpha} := \operatorname{pss}(q_{\rm pr}, \alpha, L^{\rm lc}_n)$ for each $\alpha \in \operatorname{supp} r'(L^{\rm lc}_n)$. We also derive $\mathbf{c}_{n,1}$ from $p'(L^{\rm lc}_n)$ by restricting its domain to $\operatorname{supp} p'_{\rm lc} \smallsetminus B$ and reducing the pr-indices to $P_{\mathbf{c}_{n,1},\alpha} := \operatorname{pss}(q_{\rm pr}, \alpha, L^{\rm lc}_n)$ for each $\alpha \in \operatorname{supp} p'(L^{\rm lc}_n) \smallsetminus B$. Hence, by Lemma 4.18, there is an lc-compound creature $q_{\rm lc}(L^{\rm lc}_n)$ at $L^{\rm lc}_n$ such that

- $\operatorname{dom} q_{\rm lc}(L^{\rm lc}_n) = \operatorname{dom} \mathbf{c}_{n,0} \cup \operatorname{dom} \mathbf{c}_{n,1} = \operatorname{supp} p_{\rm lc} \cup \operatorname{supp} r_{\rm lc}$;
- $\operatorname{supp} q_{\rm lc}(L^{\rm lc}_n) = \operatorname{supp} \mathbf{c}_{n,0} \cup \operatorname{supp} \mathbf{c}_{n,1} = \operatorname{supp} p'_{\rm lc}(L^{\rm lc}_n) \cup \operatorname{supp} r'_{\rm lc}(L^{\rm lc}_n)$;
- $q(L^{\rm lc}_n) \leq \mathbf{c}_{n,0}$ and $q(L^{\rm lc}_n) \leq \mathbf{c}_{n,1}$; and
- $\|q_{\rm lc}(L^{\rm lc}_n)\|^{\rm lc}_n \geq \min\{\|\mathbf{c}_{n,0}\|^{\rm lc}_n, \|\mathbf{c}_{n,1}\|^{\rm lc}_n\} - 1$, and even
$$\|q_{\rm lc}(L^{\rm lc}_n)\|^{\rm lc}_n \geq \min\{\|p'(L^{\rm lc}_n)\|^{\rm lc}_n, \|r'(L^{\rm lc}_n)\|^{\rm lc}_n\} - 1.$$

This clearly defines a condition $(q_{\rm pr}, q_{\rm lc}) \in \mathbb{Q}_{\rm lc}$ stronger than both $(p'_{\rm pr}, p'_{\rm lc})$ and $(r'_{\rm pr}, r'_{\rm lc})$. Moreover, if (i) or (ii) hold at $L^{\rm lc}_n$, we can set $q(L^{\rm al}_n) := \mathbf{c}_{n,0} \wedge \mathbf{c}_{n,1}$.

**The al-part.** Let $\operatorname{supp} q_{\rm al} = \operatorname{supp} p'_{\rm al} \cup \operatorname{supp} r'_{\rm al}$. We define $q(\alpha, \ell)$ for $(\alpha, \ell) \in \operatorname{supp} q_{\rm al} \times \operatorname{ht}^{\rm al}$ as follows: For $\ell \geq \operatorname{trl}(q, \alpha)$,

- if $\alpha \in \operatorname{supp} r'_{\rm al}$, then
$$q(\alpha, \ell) := \langle r'(\alpha, t, \ell) \mid t \in \operatorname{pss}(q_{\rm pr}, \alpha, \ell) \rangle,$$
- and if $\alpha \in \operatorname{supp} p'_{\rm al} \smallsetminus B$, then
$$q(\alpha, \ell) = \langle p'(\alpha, t, \ell) \mid t \in \operatorname{pss}(q_{\rm pr}, \alpha, \ell) \rangle;$$

for $\ell < \operatorname{trl}(q, \alpha)$, we define $q(\alpha, \ell) := r'(\alpha, \ell)$ if $\alpha \in \operatorname{supp} r'_{\rm al}$, and $q(\alpha, \ell) := p'(\alpha, \ell)$ if $\alpha \in \operatorname{supp} p'_{\rm al} \smallsetminus B$. This defines a condition $(q_{\rm pr}, q_{\rm al}) \in \mathbb{Q}_{\rm al}$ stronger than both $(p'_{\rm pr}, p'_{\rm al})$ and $(r'_{\rm pr}, r'_{\rm al})$.

It is clear that $q := (q_{\rm pr}, q_{\rm lc}, q_{\rm al}) \in \mathbb{Q}$, $q \leq p$ and $q|_B \leq r$.

In the case that $\operatorname{trl}(r) = \operatorname{trl}(p)$ and (i) or (ii) always hold, the construction actually results in $n_0 = \operatorname{trl}(r) = \operatorname{trl}(p)$, $r' = r$ and $q|_B = r$. □

**Observation 5.18.** In the previous result, when $\operatorname{trl}(r) = \operatorname{trl}(p)$ and (i) or (ii) always hold, $q$ has a concrete construction. We denote this $q$ by $p \wedge r$, which will prove to be a very useful notation, in particular for $B = S^{\rm pr}$: If $p \in \mathbb{Q}$ and $r \leq p_{\rm pr}$ in $\mathbb{Q}_{\mathrm{pr}}$ with the same trunk length as $p$, we can define the condition $p \wedge r$ (at any level $L^{\mathrm{lc}}_n$ the lc-compound of $r$ has empty domain, so (ii) is satisfied). In fact, $p \wedge r$ is derived by simply shrinking the pr-indices of the atomic creatures at the lc-part and al-part of $p$ to those given by the possibilities on $r$.

**Corollary 5.19.** *Let $B \subseteq S^*$ be closed. If for $\mathrm{t} \in \{\mathrm{pr}, \mathrm{lc}\}$ either $S^{\mathrm{t}} \subseteq B$ or $S^{\mathrm{t}} \cap B = \varnothing$, then $\mathbb{Q}_B$ is a complete subforcing of $\mathbb{Q}$. In particular, $\mathbb{Q}_{\mathrm{pr}}$, $\mathbb{Q}_{\mathrm{lc}}$ and $\mathbb{Q}_{\mathrm{al}}$ are complete subforcings of $\mathbb{Q}$.*

**Observation 5.20.** It is unclear whether $\mathbb{Q}_B$ is a complete subforcing of $\mathbb{Q}$ when $B$ is not as in Corollary 5.19. What we know is that, in such a case, the argument above using $p|_B$ may not work. For instance, fix $i_0 \in S^{\mathrm{pr}}$, assume that $|S^{\mathrm{pr}}| \geq 2$, let $B := \{i_0\}$ and $p \in \mathbb{Q}_{\mathrm{pr}}$ with $\operatorname{supp} p = \{i_0, i_1\}$ such that $i_1 \neq i_0$. We may have that, for all but finitely many $n$, $\|p(n)(i_0)\|^{\mathrm{pr}}_{\mathrm{stk}}$ is much larger than $\|p(n)(i_1)\|^{\mathrm{pr}}_{\mathrm{stk}}$, which means that $\|p(n)\|^{\mathrm{pr}}_n$ is determined by $p(n)(i_1)$. Consider some $r \leq p|_B$ in $\mathbb{Q}_B$ obtained from $p|_B$ by just increasing the halving parameters, so $r = q|_B$ where $q$ is obtained from $p$ by replacing the halving parameters with those from $r$. These new halving parameters could affect the norms of the compound creatures of $q$ in such a way that $\|q(n)\|^{\mathrm{pr}}_n$ does not diverge to infinity, i. e. $q \notin \mathbb{Q}$. This shows that $p|_B$ is not a *reduction of $p$ in $\mathbb{Q}_B$*, but does not rule out that a different reduction may exist.

We also derive the following criteria for compatibility of conditions.

**Corollary 5.21.** *If $p, r \in \mathbb{Q}$ are identical on the intersection of their supports $A := \operatorname{supp} p \cap \operatorname{supp} r$ (i. e. $p|_A = r|_A$) and they have the same halving parameters at the levels for which the supports of the two compound creatures in both conditions are non-empty, then there is a condition $q$ stronger than both $p$ and $r$ such that $\operatorname{supp} q = \operatorname{supp} p \cup \operatorname{supp} r$.*

*Proof.* Let $B := \operatorname{supp} r$. Clearly, $r \in \mathbb{Q}_B$ and $r \leq p|_B$, because $p|_B = p|_A = r|_A$. Since the preconditions of Lemma 5.17 hold, there is some $q \leq p$ in $\mathbb{Q}$ with $\operatorname{supp} q = \operatorname{supp} p \cup \operatorname{supp} r$ such that $q|_B \leq r$. Since $q \leq q|_B$, $q$ is as required. □

As a consequence, we can show that the generic real for any $\alpha \in S^*$ is always defined, i. e. the set of conditions $p$ with $\alpha \in \operatorname{supp} p$ is dense.

**Corollary 5.22.** *Given $p \in \mathbb{Q}$ and $\alpha \in S^*$, there is some $q \leq p$ such that $\operatorname{supp} q = \operatorname{supp} p \cup \{i^*(\alpha), \alpha\}$.*

*Proof.* If $\alpha \in \operatorname{supp} p$, then we can just set $q := p$; hence we assume that $\alpha \notin \operatorname{supp} p$. In the case that $i^*(\alpha) \notin \operatorname{supp} p$, by Corollary 5.11, we can find $r_0 \in \mathbb{Q}$ with support $\{i^*(\alpha), \alpha\}$ – moreover, according to the proofs of Lemma 4.16, Lemma 4.17, Lemma 5.9 and Lemma 5.10, we can even find $r_0$ such that all its subatoms are maximal and with halving parameters equal to 0. So we can modify $r_0$ by defining the halving parameters to be the same as those of $p$ at levels for which the compounds of both conditions have non-empty support, and the resulting $r$ is still in $\mathbb{Q}$ with support $\{i^*(\alpha), \alpha\}$. Since $p$ and $r$ have disjoint support, we can find the desired $q$ by Corollary 5.21.

Now assume that $i^*(\alpha) \in \operatorname{supp} p$, so we define $r_{\mathrm{pr}} := p|\{i^*(\alpha)\}$ and use Lemma 5.10 to find some $r \in \mathbb{Q}$ with support $\{i^*(\alpha), \alpha\}$ whose pr-part is $r_{\mathrm{pr}}$. Again, the subatoms of $r$ at $(\alpha, \ell)$ for $\ell > \operatorname{trl}(r) = \operatorname{trl}(p)$ are maximal and the halving parameters are 0 if $\alpha \in S^{\mathrm{lc}}$, so in this case we can modify $r$ by changing the halving parameters as before and still obtain a condition in $\mathbb{Q}$ with the same support. Since $p$ and $r$ are identical on $\operatorname{supp} p \cap \operatorname{supp} r = \{i^*(\alpha)\}$, we can find the desired $q$ by Corollary 5.21. ☐

As a consequence of Corollary 5.21, we can prove a chain condition for our forcing.

**Lemma 5.23.** *The poset $\mathbb{Q}$ has $\mathfrak{c}^+$-cc. In particular, under* CH, *$\mathbb{Q}$ has $\aleph_2$-cc.*

*Proof.* Let $A = \{p_\xi \mid \xi < \mathfrak{c}^+\}$ be a collection of conditions in $\mathbb{Q}$. Note that the set of conditions such that their halving parameters are all rational is dense, so we may without loss of generality assume that each $p_\xi$ satisfies this property.

By the $\Delta$-system lemma applied to the set of supports of conditions in $A$, we can find $C \subseteq \mathfrak{c}^+$ of size $\mathfrak{c}^+$ such that $\{\operatorname{supp} p_\xi \mid \xi \in C\}$ forms a $\Delta$-system with root $D \subseteq S^*$. We can even shrink $C$ in such a way that

(i) for $n < \omega$, there is some $f_n^{\mathrm{pr}}\colon D \cap S^{\mathrm{pr}} \to 2$ such that $f_n^{\mathrm{pr}}$ coincides with the characteristic function of $\operatorname{supp} p_{\xi,\mathrm{pr}}(n) \cap D$ on $D \cap S^{\mathrm{pr}}$ for all $\xi \in C$ (which implies that $\{\operatorname{supp} p_{\xi,\mathrm{pr}}(n) \mid \xi \in C\}$ forms a $\Delta$-system with root $D_n^{\mathrm{pr}} := (f_n^{\mathrm{pr}})^{-1}[\{1\}]$ and that $\operatorname{supp} p_{\xi,\mathrm{pr}}(n) \setminus D_n^{\mathrm{pr}} \subseteq S^{\mathrm{pr}} \setminus D$ for all $\xi \in C$);
(ii) for $n < \omega$, there is some $f_n^{\mathrm{lc}}\colon D \cap S^{\mathrm{lc}} \to 2$ such that $f_n^{\mathrm{lc}}$ coincides with the characteristic function of $\operatorname{supp} p_\xi(L_n^{\mathrm{lc}}) \cap D$ on $D \cap S^{\mathrm{lc}}$ for all $\xi \in C$ (which implies that $\{\operatorname{supp} p_\xi(L_n^{\mathrm{lc}}) \mid \xi \in C\}$ forms a $\Delta$-system with root $D_n^{\mathrm{lc}} := (f_n^{\mathrm{lc}})^{-1}[\{1\}]$, and that $\operatorname{supp} p_\xi(L_n^{\mathrm{lc}}) \setminus D_n^{\mathrm{lc}} \subseteq S^{\mathrm{lc}} \setminus D$ for all $\xi \in C$);
(iii) the halving parameters at all levels are the same for all $p_\xi$ with $\xi \in C$; and
(iv) there is some $r \in \mathbb{Q}_D$ such that $p_\xi|_D = r$ for all $\xi \in C$.

By Corollary 5.21, $p_\xi$ and $p_\eta$ are compatible for all $\xi, \eta \in C$. Therefore $A$ is not an antichain. ☐

We now finally establish how the generic reals are defined from the generic set.

**Definition 5.24.** Let $G$ be a $\mathbb{Q}$-generic set over the ground model $V$. We denote by $\dot{y}$ the $\mathbb{Q}$-name for

$$\{((\alpha, \ell), z) \mid (\alpha, \ell) \in \mathrm{DOM} \text{ and } \exists p \in G\colon \alpha \in \operatorname{supp} p,\ \operatorname{trl}(p) > \ell \text{ and } p(\alpha, \ell) = \{z\}\}.$$

For each $\mathrm{t} \in \{\mathrm{pr}, \mathrm{al}, \mathrm{lc}\}$ and each $\alpha \in S^{\mathrm{t}}$, let $\dot{y}_\alpha$ be a $\mathbb{Q}$-name for $\{(\ell, z) \mid ((\alpha, \ell), z) \in \dot{y}\}$.

As in Subsection 3.2 and keeping in mind Lemma 5.26, we now define several $\mathbb{Q}$-names – concretely, for $i \in S^{\mathrm{pr}}$ and $\mathrm{t} \in \{\mathrm{lc}, \mathrm{al}\}$, we use Definition 3.3 to define the following:

(1) $\dot{y}_i^*$ is a $\mathbb{Q}$-name for the element of $\prod_{n<\omega} T_n^*$ defined by $\dot{y}_i^*(n) := \langle \dot{y}_i(L) \mid L \in \mathrm{ht}_n^{\mathrm{pr}} \rangle$.
(2) $\dot{d}_i$ is a $\mathbb{Q}$-name for $d_{\dot{y}_i^*}$; $\dot{h}_i$, $\dot{b}_i$, $\dot{g}_i$ and $\dot{a}_i$ are defined similarly.
(3) $\dot{b}_i^-$ is a $\mathbb{Q}$-name for $b_{\dot{y}_i^*}^-$, and $\dot{h}_i^-$ is defined similarly.
(4) $\dot{d}_i^{\mathrm{t}}$ is a $\mathbb{Q}$-name for $d_{\dot{y}_i^*}^{\mathrm{t}}$; $g_i^{\mathrm{t}}$, $h_i^{\mathrm{t}}$, $b_i^{\mathrm{t}}$, $f_i^{\mathrm{t}}$ and $a_i^{\mathrm{t}}$ are defined similarly.

(5) $\dot{\vec{\nu}}^{\mathrm{t}}_i$ is a $\mathbb{Q}$-name for $\vec{\nu}^{\mathrm{t}}_{\dot{y}^*_i}$.
(6) For $\alpha \in S^{\mathrm{lc}}_i$, $\dot{\varphi}_\alpha$ is a $\mathbb{Q}$-name for a function with domain $\omega$ such that $\dot{\varphi}_\alpha(n) := \dot{b}^{\mathrm{lc}}_i(n) \smallsetminus \prod_{i\in I^*_n}(\dot{b}_i(\ell) \smallsetminus \dot{y}_\alpha(\ell))$.
(7) For $\alpha \in S^{\mathrm{al}}_i$, $\dot{\varphi}_\alpha$ is a $\mathbb{Q}$-name for a function with domain $\omega$ such that $\dot{\varphi}_\alpha(n) := y_\alpha(L^{\mathrm{al}}_n)$.

We remove the dots whenever these names are evaluated in any generic extension.

The following results show the type of reals added by $\mathbb{Q}$.

**Fact 5.25.** If $p \in \mathbb{Q}$, $L \in \mathrm{Ht}$ and $\eta \in \mathrm{poss}(p, {<}L)$, then $p \wedge \eta$ forces that $\dot{y}$ is a function extending $\eta$.

**Lemma 5.26.** *Assume that $G$ is a $\mathbb{Q}$-generic over $V$ and fix $i \in S^{\mathrm{pr}}$. Then $\mathbb{Q}$ forces the following:*

*(a) $\dot{y}$ is a function with domain* DOM.
*(b) $\dot{y}_i$ is a function in $\prod_{L\in\mathrm{ht}^{\mathrm{pr}}} T^*_L$, hence $\dot{y}^*_i$ is in $\prod_{n<\omega} T^*_n$ and both $\dot{\vec{\nu}}^{\mathrm{lc}}_i$ and $\dot{\vec{\nu}}^{\mathrm{al}}_i$ are blocks.*
*(c) For each $\alpha \in S^{\mathrm{lc}}_i$, $\dot{y}_\alpha$ is in $\mathcal{S}^*(\dot{b}^-_i, \dot{h}^-_i)$ (and thus $\dot{\varphi}_\alpha$ is a slalom in $\mathcal{S}^*(\dot{b}^{\mathrm{al}}_i, \dot{h}^{\mathrm{al}}_i)$ by the proof of Lemma 2.14).*
*(d) For each $\alpha \in S^{\mathrm{al}}_i$, $\dot{\varphi}_\alpha$ is a slalom in $\mathcal{S}^*(\dot{a}^{\mathrm{al}}_i, \dot{h}^{\mathrm{al}}_i)$.*

*Proof.* To prove (a): Clearly, $\mathbb{Q}$ forces that $\dot{y}$ is a function with domain contained in DOM. To show equality, it is enough to prove that, for any $(\alpha, \ell) \in \mathrm{DOM}$, the set $D_{\alpha,\ell} := \{p \in \mathbb{Q} \mid \alpha \in \operatorname{supp} p,\ \mathrm{trl}(p) > \ell\}$ is dense. If $p \in \mathbb{Q}$, by Corollary 5.22 there is some $q \leq p$ in $\mathbb{Q}$ such that $\alpha \in \operatorname{supp} q$. Choosing $n > \ell$ in $\omega$ and $\eta \in \mathrm{poss}(q, {<}n)$, we get $q \wedge \eta \in D_{\alpha,\ell}$ stronger than $q$.

To see (b), let $p \in \mathbb{Q}$ and $L \in \mathrm{ht}^{\mathrm{pr}}$, and choose $q \in D_{i,L}$ stronger than $p$. So $\mathrm{poss}(q, < \mathrm{trl}(q))$ only contains one possibility $\eta$ and, by Fact 5.13 (a) and Fact 5.25, $q$ forces $\dot{y}_i(L) = \eta(i, L) \in T^*_L$. The rest is clear by Lemma 3.4.

For (c) it is enough to show that, for any $p \in \mathbb{Q}$, $\alpha \in \operatorname{supp} p_{\mathrm{lc}}$ and $\ell \in \mathrm{ht}^{\mathrm{lc}}$, $p$ forces that $\dot{y}_\alpha(\ell) \subseteq b_{t^{\mathrm{mx}}_n}(\ell)$ has size $\leq h_{t^{\mathrm{mx}}_n}(\ell)$, where $n := n_*(\ell)$, and even more, whenever $\ell \geq \mathrm{trl}(p, \alpha)$, $p$ forces that $\dot{y}_\alpha(\ell) \subseteq \dot{b}_i(\ell)$ has size $\leq \dot{h}_i(\ell)$. Pick some $n' > \ell$ and let $\eta \in \mathrm{poss}(p, {<}n')$. By Fact 5.13 (b) and Fact 5.25, $p \wedge \eta$ forces $\dot{y}_\alpha(\ell) = \eta(\alpha, \ell) \subseteq b_{t^{\mathrm{mx}}_n}(\ell)$ of size $\leq h_{t^{\mathrm{mx}}_n}(\ell)$ and, whenever $\ell \geq \mathrm{trl}(p, \alpha)$, $\dot{y}_\alpha(\ell) = \eta(\alpha, \ell) \subseteq b_{\hat{\eta}(\alpha,\ell)}(\ell)$ has size $\leq h_{\hat{\eta}(\alpha,\ell)}(\ell)$. But since $p \wedge \eta$ forces $\hat{\eta}(\alpha, \ell) = y^*_i(n)$, we also have that $p \wedge \eta$ forces the desired conclusion. By Fact 5.15 (d), $p$ forces the same.

To prove (d), similar to (c) we can show that, for any $p \in \mathbb{Q}$, $\alpha \in \operatorname{supp} p_{\mathrm{lc}}$ and $\ell \in \mathrm{ht}^{\mathrm{al}}$, $p$ forces that $\dot{y}_\alpha(\ell) \subseteq a_{t^{\mathrm{mx}}_n}(\ell)$ has size $\leq h_{t^{\mathrm{mx}}_n}(\ell)$, where $n := n_*(\ell)$, and even more, whenever $\ell \geq \mathrm{trl}(p, \alpha)$, $p$ forces that $\dot{y}_\alpha(\ell) \subseteq \dot{a}_i(\ell)$ has size $\leq \dot{h}_i(\ell)$. □

## 6. Bigness

The notion of bigness is a tool in creature forcing which allows us to homogenise the decisions made by a condition or a creature (e. g. deciding the name for an ordinal). This is not only essential for proving that $\mathbb{Q}$ is proper, but also for the proof of the main theorem.

Bigness is defined for subatomic creatures as follows:

**Definition 6.1.** Let $\varepsilon > 0$ be a real number, let $B \in \omega$ and let K be a subatomic family.

(1) A subatom $\mathbf{c} \in \mathrm{K}$ has $(B, \varepsilon)$-*bigness* if for each function $F\colon \mathbf{c} \to B$, there is a $\mathbf{d} \subseteq \mathbf{c}$ in $\mathbf{K}$ such that $F{\restriction}_{\mathbf{d}}$ is constant and $\|\mathbf{d}\| \geq \|\mathbf{c}\| - \varepsilon$.
(2) We say that K has $(B, \varepsilon)$-*bigness* if each $\mathbf{c} \in \mathrm{K}$ has $(B, \varepsilon)$-bigness.

One of the reasons we defined the parameters in Subsection 3.3 the way we did is to be able to explicitly calculate the bigness of the subatomic creatures in our forcing.

**Lemma 6.2.** *Let $\ell \in \mathrm{ht}$.*

*(a) If $\ell \in \mathrm{ht}^{\mathrm{pr}}$, then $\mathrm{K}_\ell$ has $(n^B_\ell, 1/n^B_\ell)$-bigness.*
*(b) If $\ell \in \mathrm{ht}^{\mathrm{lc}} \cup \mathrm{ht}^{\mathrm{al}}$ and $t \in T^*_{n_*(\ell)}$, then $\mathrm{K}^t_\ell$ has $(d_t(\ell), 1/d_t(\ell))$-bigness.*

*Proof.* To see (a), let $\mathbf{c} \in K_\ell$ and $F\colon \mathbf{c} \to n^B_\ell$. We find $\mathbf{d} \subseteq \mathbf{c}$ with $|\mathbf{d}| \geq \frac{|\mathbf{c}|}{n^B_\ell}$ such that $F{\restriction}_{\mathbf{d}}$ is constant (in fact, $\mathbf{d}$ is the maximal $F^{-1}[\{k\}]$). Hence

$$\|\mathbf{d}\|_\ell = \frac{1}{n^B_\ell}\log_{n^B_\ell}|\mathbf{d}| \geq \frac{1}{n^B_\ell}\left(\log_{n^B_\ell}\frac{|\mathbf{c}|}{n^B_\ell}\right) = \frac{1}{n^B_\ell}(\log_{n^B_\ell}|\mathbf{c}| - \log_{n^B_\ell} n^B_\ell) = \|\mathbf{c}\|^{\mathrm{pr}}_n - \frac{1}{n^B_\ell}.$$

For (b), see [KM22, Lemma 3.10]. □

We also describe and calculate the bigness of compound creatures. Bigness here homogenises functions whose domain is the set of possibilities of the compound creature; however, note that the possibilities for the lc-compounds also depend on the pr-indices. Inspired by the proof of the main result in [KS12] (where continuum many $\mathfrak{c}^\forall_{b,h}$ and $\mathfrak{c}^\exists_{b,h}$ are separated), we state our bigness results for compound creatures, which are essential to separate the many cardinal characteristics in the proof of the main theorem. We start with lc-compounds, for which we use the lexicographic order of $T^*_n \times I^*_n$ presented in Subsection 3.3.

**Lemma 6.3.** *Let $n < \omega$, $\mathbf{c}$ a lc-compound creature at $L^{\mathrm{lc}}_n$, $\bar{t} = \langle t_\alpha \mid \alpha \in \operatorname{supp}\mathbf{c}\rangle \in \prod_{\alpha\in\operatorname{supp}\mathbf{c}} P_{\mathbf{c},\alpha}$, $t_0 \in T^*_n$, and $\ell_0 \in I^*_n$. Assume that $M \in \omega$, $f\colon \operatorname{poss}(\mathbf{c}, \bar{t}) \to M$ and $M^{m^S_{t_0,\ell_0}} \leq d_{t_0^+}(\ell_0^+)$, the latter whenever $(t_0^+, \ell_0^+)$ is the successor of $(t_0, \ell_0)$ in $T^*_n \times I^*_n$ with regard to the lexicographic order. Then there is a lc-compound creature $\mathbf{c}' \leq \mathbf{c}$ with the same domain, the same halving parameter and the same pr-indices such that, for any $\alpha \in \operatorname{supp}\mathbf{c}$, $\ell \in I^*_n$ and $t \in P_{\mathbf{c},\alpha}$,*

*(i) if $t = t_\alpha$ and $(t_\alpha, \ell) > (t_0, \ell_0)$, then $\|\mathbf{c}'(\alpha, t, \ell)\|^t_\ell \geq \|\mathbf{c}(\alpha, t, \ell)\|^t_\ell - \frac{1}{d_t(\ell)}$;*
*(ii) otherwise, $\mathbf{c}'(\alpha, t, \ell) = \mathbf{c}(\alpha, t, \ell)$; and*
*(iii) $f{\restriction}_{\operatorname{poss}(\mathbf{c}',\bar{t})}$ only depends on $C := \{(\alpha, \ell) \mid (t_\alpha, \ell) \leq (t_0, \ell_0)\}$, i. e. if $\eta, \eta' \in \operatorname{poss}(\mathbf{c}', \bar{t})$ coincide on $C$, then $f(\eta) = f(\eta')$.*

*Proof.* When $(t_0, \ell_0)$ is the maximum of $T^*_n \times I^*_n$, we can just set $\mathbf{c}' := \mathbf{c}$; hence we assume that $(t_0, \ell_0)$ is not the maximum of $T^*_n \times I^*_n$. Note that $\{\eta|_{\operatorname{dom}\mathbf{c}\smallsetminus\operatorname{supp}\mathbf{c}} \mid \eta \in \operatorname{poss}(\mathbf{c}, \bar{t})\}$ has only one element, which we denote by $\eta_-$.

We define $F\colon \prod_{(\alpha,\ell)\notin C} \mathbf{c}(\alpha, t_\alpha, \ell) \to M^{\prod_{(\alpha,\ell)\in C} \mathbf{c}(\alpha,t_\alpha,\ell)}$ such that $F(\eta_0)(\eta_1) := f(\eta_0 \cup \eta_1 \cup \eta_-)$. Note that $\prod_{(\alpha,\ell)\in C} |\mathbf{c}(\alpha, t_\alpha, \ell)| \leq m^S_{t_0,\ell_0}$ by modesty, (pr9) and Fact 4.6, so $|\operatorname{ran} F| \leq M^{m^S_{t_0,\ell_0}} \leq d_{t_0^+}(\ell_0^+)$ follows from the preconditions. In this way, we just need to find $\mathbf{c}'$ satisfying (i) and (ii) and such that $F$ is constant on $\prod_{(\alpha,\ell)\notin C} \mathbf{c}'(\alpha, t_\alpha, \ell)$.

Consider the increasing enumeration $\{(t'_k, \ell'_k) \mid k < m\}$ of all $(t, \ell) \in T^*_n \times I^*_n$ larger than $(t_0, \ell_0)$ such that, for some $\alpha \in \operatorname{supp} \mathbf{c}$, $t_\alpha = t$ and $\mathbf{c}(\alpha, t, \ell)$ is non-trivial. By strong modesty, such an $\alpha$ is unique, so we denote by $\alpha_k$ the one corresponding to $(t'_k, \ell'_k)$. Moreover, the non-trivial creatures in $\{\mathbf{c}(\alpha, t_\alpha, \ell) \mid (\alpha, \ell) \notin C\}$ are precisely $\{\mathbf{c}(\alpha_k, t'_k, \ell'_k) \mid k < m\}$, hence $F$ only depends on $\{(\alpha_k, \ell'_k) \mid k < m\}$. Hence we can define a function $F'\colon \prod_{k<m} \mathbf{c}(\alpha_k, t'_k, \ell'_k) \to \operatorname{ran} F$ that determines $F$.

We define $\langle \mathbf{c}'(\alpha_k, t'_k, \ell'_k) \mid k < m\rangle$ by decreasing induction on $k$ such that

($k$-i) $\mathbf{c}'(\alpha_k, t'_k, \ell'_k) \subseteq \mathbf{c}(\alpha_k, t'_k, \ell'_k)$, $\|\mathbf{c}'(\alpha_k, t'_k, \ell'_k)\| \geq \|\mathbf{c}(\alpha, t'_k, \ell'_k)\| - \frac{1}{d_{t'_k}(\ell'_k)}$; and

($k$-ii) $F'{\restriction}_{\prod_{k'<k} \mathbf{c}(\alpha_{k'}, t'_{k'}, \ell'_{k'}) \times \prod_{k' \geq k} \mathbf{c}'(\alpha_{k'}, t'_{k'}, \ell'_{k'})}$ only depends on $\prod_{k'<k} \mathbf{c}(\alpha_{k'}, t'_{k'}, \ell'_{k'})$.

At step $k$, by the induction hypothesis, $F'{\restriction}_{\prod_{k'\leq k} \mathbf{c}(\alpha_{k'}, t'_{k'}, \ell'_{k'}) \times \prod_{k'>k} \mathbf{c}'(\alpha_{k'}, t'_{k'}, \ell'_{k'})}$ only depends on $\prod_{k'\leq k} \mathbf{c}(\alpha_{k'}, t'_{k'}, \ell'_{k'})$. Hence we can define

$$F_k\colon \mathbf{c}(\alpha_k, t'_k, \ell'_k) \to (\operatorname{ran} F)^{\prod_{k'<k} \mathbf{c}(\alpha_{k'}, t'_{k'}, \ell'_{k'})}$$

such that $F_k(i_k)(\langle i_{k'} \mid k' < k\rangle)$ is the value of $F'$ calculated from $\langle i_{k'} \mid k' \leq k\rangle$.

We claim that $\left|(\operatorname{ran} F)^{\prod_{k'<k} \mathbf{c}(\alpha_{k'}, t'_{k'}, \ell'_{k'})}\right| \leq d_{t'_k}(\ell'_k)$. In case $k = 0$, we just have to show $|\operatorname{ran} F| \leq d_{t'_0}(\ell'_0)$, which immediately follows from $|\operatorname{ran} F| \leq d_{t_0^+}(\ell_0^+) \leq d_{t'_0}(\ell'_0)$. Assume $k > 0$; then, by (pr9),

$$\left|(\operatorname{ran} F)^{\prod_{k'<k} \mathbf{c}(\alpha_{k'}, t'_{k'}, \ell'_{k'})}\right| \leq |\operatorname{ran} F|^{m^S_{t'_{k-1},\ell'_{k-1}}} \leq d_{t'_0}(\ell'_0)^{m^S_{t'_{k-1},\ell'_{k-1}}} \leq d_{t'_k}(\ell'_k).$$

Hence we can apply Lemma 6.2 (b) to find $\mathbf{c}'(\alpha_k, t'_k, \ell'_k) \leq \mathbf{c}(\alpha_k, t'_k, \ell'_k)$ in $\mathbf{K}^{t'_k}_{\ell'_k}$ such that $F_k{\restriction}_{\mathbf{c}'(\alpha_k, t'_k, \ell'_k)}$ is constant and $\|\mathbf{c}'(\alpha_k, t'_k, \ell'_k)\|^{t'_k}_{\ell'_k} \geq \|\mathbf{c}(\alpha_k, t'_k, \ell'_k)\|^{t'_k}_{\ell'_k} - \frac{1}{d_{t'_k}(\ell'_k)}$. It is clear that both ($k$-i) and ($k$-ii) hold.

We define the other components of $\mathbf{c}'$ to be identical to $\mathbf{c}$; then $\mathbf{c}'$ is exactly as required. □

**Corollary 6.4.** *Let $n < \omega$, let $\mathbf{c}$ be a* lc*-compound creature at $L^{\mathrm{lc}}_n$, let $\bar{t} = \langle t_\alpha \mid \alpha \in \operatorname{supp} \mathbf{c}\rangle \in \prod_{\alpha\in\operatorname{supp}\mathbf{c}} P_{\mathbf{c},\alpha}$ and let $f\colon \operatorname{poss}(\mathbf{c}, \bar{t}) \to n^B_{L^{\mathrm{lc}}_n}$. Then there is an* lc*-compound creature $\mathbf{c}' \leq \mathbf{c}$ with the same domain, the same halving parameter and the same pr-indices, such that, for any $\alpha \in \operatorname{supp} \mathbf{c}$, $\ell \in I^*_n$ and $t \in P_{\mathbf{c},\alpha}$,*

*(i) if $t = t_\alpha$, then $\|\mathbf{c}'(\alpha, t, \ell)\|^t_\ell \geq \|\mathbf{c}(\alpha, t, \ell)\|^t_\ell - \frac{1}{d_t(\ell)}$;*
*(ii) otherwise, $\mathbf{c}'(\alpha, t, \ell) = \mathbf{c}(\alpha, t, \ell)$; and*
*(iii) $f{\restriction}_{\operatorname{poss}(\mathbf{c}',\bar{t})}$ is constant.*

*Proof.* Let $(t_0, \ell_0) := \min(T^*_n \times I^*_n)$. Since $n^B_{L^{\mathrm{lc}}_n} \leq d_{t_0}(\ell_0) < b_{t_0}(\ell_0)$ and $b_{t_0}(\ell_0)^{m^S_{t_0,\ell_0}} \leq d_{t'_0}(\ell'_0)$ by (pr9), where $(t'_0, \ell'_0)$ is the successor of $(t_0, \ell_0)$ in $T^*_n \times I^*_n$, we can obtain $\mathbf{c}^1 \leq \mathbf{c}$ as in Lemma 6.3. Then $f{\restriction}_{\operatorname{poss}(\mathbf{c}^1,\bar{t})}$ depends only on $\{(\alpha, \ell_0) \mid \alpha \in \operatorname{supp} \mathbf{c}^1,\ t_\alpha = t_0\}$; note that $\mathbf{c}^1(\bar{t})$ and $\mathbf{c}(\bar{t})$ coincide for those coordinates. If this set is empty, we are done. Otherwise, choose $\alpha_0 \in \operatorname{supp}\mathbf{c}$ such that $t_{\alpha_0} = t_0$, and additionally such that $\mathbf{c}(\alpha_0, t_0, \ell_0)$ is non-trivial in case there is such an $\alpha_0$ (in which case it is unique by modesty). Hence $f{\restriction}_{\mathrm{poss}(\mathbf{c}^1,\bar{t})}$ depends only on $\{(\alpha_0, \ell_0)\}$, i. e. it can be reconstructed from some function $f_0\colon \mathbf{c}(\alpha_0, t_0, \ell_0) \to n^B_{L^{\mathrm{lc}}_n}$. By (pr10), $n^B_{L^{\mathrm{lc}}_n} \leq d_{t_0}(\ell_0)$, so Lemma 6.2 (b) implies that there is some $\mathbf{c}'(\alpha_0, t_0, \ell_0) \leq \mathbf{c}(\alpha_0, t_0, \ell_0)$ such that $\|\mathbf{c}'(\alpha_0, t_0, \ell_0)\|^{t_0}_{\ell_0} \geq \|\mathbf{c}(\alpha_0, t_0, \ell_0)\|^{t_0}_{\ell_0} - \frac{1}{d_{t_0}(\ell_0)}$ and $f_0{\restriction}_{\mathbf{c}'(\alpha_0,t_0,\ell_0)}$ is constant.

We define the other components of $\mathbf{c}'$ to be identical to $\mathbf{c}^1$; then $\mathbf{c}'$ is exactly as required. □

A similar (and simpler) argument using Lemma 6.2 (a), (pr2) and (pr4) allows us to calculate the bigness of pr-compound creatures.

**Lemma 6.5.** *Let $\mathbf{c}$ be a pr-compound creature at $(n, 0)$, let $L_0 \in \mathrm{ht}^{\mathrm{pr}}_n$ and $M \in \omega$, and let $f\colon \mathrm{poss}(\mathbf{c}) \to M$ be such that $M^{\prod\{T^*_L \mid L \in \mathrm{ht}^{\mathrm{pr}}_n \cap L_0{\downarrow}\}} \leq n^B_{L_0}$. Then there is a pr-compound creature $\mathbf{c}' \leq \mathbf{c}$ with the same domain and halving parameter such that, for any $i \in \operatorname{supp}\mathbf{c}$ and $L \in \mathrm{ht}^{\mathrm{pr}}_n$,*

*(i) if $L \geq L_0$, then $\|\mathbf{c}'(i, L)\| \geq \|\mathbf{c}(i, L)\| - \frac{1}{n^B_L}$;*
*(ii) otherwise, $\mathbf{c}'(i, L) = \mathbf{c}(i, L)$; and*
*(iii) for $\eta \in \mathrm{poss}(\mathbf{c}')$, $f(\eta)$ depends only on $\eta{\restriction}_{\operatorname{supp}\mathbf{c}\times L_0{\downarrow}}$.*

*In particular, if $L_0 = (n, 0)$ and $M \leq n^B_{L_0}$, then $f{\restriction}_{\mathrm{poss}(\mathbf{c}')}$ is constant.*

*Proof.* We only check that, for $L' \geq L_0$ in $\mathrm{ht}^{\mathrm{pr}}_n$, we have $M^{\prod\{T^*_L \mid L \in \mathrm{ht}^{\mathrm{pr}}_n \cap L'{\downarrow}\}} \leq n^B_{L'}$ (the rest follows analogously to the proof of Lemma 6.3). We proceed by induction on $L'$. The case $L' = L_0$ holds from the assumption. If $L''$ is the successor of $L'$ in $\mathrm{ht}^{\mathrm{pr}}_n$, then

$$M^{\prod\{T^*_L \mid L \in \mathrm{ht}^{\mathrm{pr}}_n \cap L''{\downarrow}\}} \leq (n^B_{L'})^{T^*_{L'}} < (n^P_{<L''})^{n^P_{<L''}} < n^B_{L''},$$

where the first "<" holds by (pr4) and the second by (pr2). □

## 7. Continuous and Rapid Reading

This section is dedicated to the continuous reading of names, whose proof gives us several fundamental properties of $\mathbb{Q}$ (such as $\mathbb{Q}$ being proper and $\omega^\omega$-bounding). We also define a notion of rapid reading of names.

**Definition 7.1.** Let $B \subseteq S^*$, let $\dot{\tau}$ a $\mathbb{Q}$-name for a set in the ground model $V$, and let $p \in \mathbb{Q}$.

(1) For $L \in \mathrm{Ht}$, we say that $\dot{\tau}$ is *$L$-decided by $p$* if $p \wedge \eta$ decides the values of $\dot{\tau}$ for each $\eta \in \mathrm{poss}(p, {<}L)$. In other words, there is some map $T$ in the ground model with domain $\mathrm{poss}(p, {<}L)$ such that $p \wedge \eta \Vdash \dot{\tau} = T(\eta)$ for all $\eta \in \mathrm{poss}(p, {<}L)$.
(2) We moreover say that $\dot{\tau}$ is $L$-decided by $p$ *only using indices in $B$* if the value $T(\eta)$ depends only on $\eta|_B$ for all $\eta \in \mathrm{poss}(p, {<}L)$, i. e. if $\eta' \in \mathrm{poss}(p, {<}L)$ coincides with $\eta$ on $B \times L{\downarrow}$, then both $p \wedge \eta$ and $p \wedge \eta'$ decide the same value of $\dot{\tau}$.
(3) We say that $p$ *essentially decides $\dot{\tau}$ (only using indices in $B$)* if $\dot{\tau}$ is $L$-decided by $p$ (only using indices in $B$) for some $L \in \mathrm{Ht}$.

(4) Assume that $p \Vdash \dot{r}\colon X \to V$, where $X \in V$. We say $p$ *continuously reads* $\dot{r}$ *(only using indices in* $B$*)* if $p$ essentially decides each $\dot{r}(z)$ (only using indices in $B$) for all $z \in X$.
(5) Assume that $p \Vdash \dot{r}\colon \mathrm{Ht} \to V$. We say that $p$ *rapidly reads* $\dot{r}$ *(only using indices in* $B$*)* if, for each $L \in \mathrm{Ht}$, $\dot{r}{\restriction}_{L\downarrow}$ is $L$-decided (only using indices in $B$).
(6) In the previous notions, we can alternatively define both deciding and reading *without using indices in* $B$, which is to be read as *only using indices in* $S^* \smallsetminus B$.

One of the main results of this section is:

**Theorem 7.2.** *The forcing* $\mathbb{Q}$ *is proper,* $\omega^\omega$*-bounding and it has continuous reading of names, i. e. if* $X$ *is a countable set (in the ground model) and* $\dot{r}$ *is a* $\mathbb{Q}$*-name for a function from* $X$ *into the ground model, then the set of conditions continuously reading* $\dot{r}$ *is dense.*

The extensive proof of this theorem is presented in Subsection 7.2. However, we first discuss some consequences and related properties, e. g. that for some $\dot{r}$ we actually get that the set of conditions rapidly reading $\dot{r}$ is dense (Corollary 7.9). Preservation of cardinals also follows under CH by Lemma 5.23.

**Corollary 7.3.** *Under* CH*,* $\mathbb{Q}$ *preserves all cofinalities and cardinalities.*

For $p, q \in \mathbb{Q}$, the notation $q \leq^* p$ usually means that any condition in $\mathbb{Q}$ stronger than $q$ is compatible with $p$. This is equivalent to the fact that $q$ forces $p$ to be inside the generic set. We write $p =^* q$ when $p \leq^* q$ and $q \leq^* p$.

**Lemma 7.4.** *Let* $p, q \in \mathbb{Q}$*, let* $\dot{\tau}$ *be a* $\mathbb{Q}$*-name for a set in the ground model* $V$*, let* $X \in V$ *and let* $\dot{r}$ *be a* $\mathbb{Q}$*-name for a function from* $X$ *into* $V$*. If* $q \leq^* p$ *and* $\operatorname{supp} p \subseteq \operatorname{supp} q$*, then*

*(a) if* $L \in \mathrm{Ht}$ *and* $\dot{\tau}$ *is* $L$*-decided by* $p$*, then it is also* $L$*-decided by* $q$*;*
*(b) if* $p$ *continuously reads* $\dot{r}$*, then* $q$ *continuously reads* $\dot{r}$*; and*
*(c) if* $X = \mathrm{Ht}$ *and* $p$ *rapidly reads* $\dot{r}$*, then* $q$ *rapidly reads* $\dot{r}$*.*

*The same results hold when adding "only (or without) using indices in* $B$*" for any* $B \subseteq S^*$*.*

*Proof.* To prove (a): Since $\dot{\tau}$ is $L$-decided by $p$, there is a map $T$ as in Definition 7.1 (1), i.e. such that $p \wedge \eta \Vdash \dot{\tau} = T(\eta)$ for any $\eta \in \operatorname{poss}(p, {<}L)$. Note that $\operatorname{poss}(q|_{\operatorname{supp} p}, {<}L) \subseteq \operatorname{poss}(p, {<}L)$ because $q \leq^* p$ and $\operatorname{supp} p \subseteq \operatorname{supp} q$. If $\eta \in \operatorname{poss}(q, {<}L)$, then $\eta' := \eta|_{\operatorname{supp} p} \in \operatorname{poss}(p, {<}L)$, so $q \wedge \eta \leq^* p \wedge \eta'$ and it follows that $q \wedge \eta \Vdash \dot{\tau} = T(\eta')$.

(b) and (c) are immediate consequences of (a). □

A note about the previous lemma: [FGKS17, Lemma 5.1.3] and [GK21, Lemma 6.4] claim that whenever $q \leq^* p$ and $\dot{\tau}$ is essentially decided by $p$, then it is also essentially decided by $q$, without explicitly assuming $\operatorname{supp} p \subseteq \operatorname{supp} q$. However, in our construction, this is *not* true without explicitly assuming $\operatorname{supp} p \subseteq \operatorname{supp} q$. For example, let $i_0 \in S^{\mathrm{pr}}$ and let $p_0 \in \mathbb{Q}$ be the condition with $\operatorname{supp} p = \{i_0\}$ and $\operatorname{trl}(p_0) = 0$, where the subatom at $L \in \mathrm{ht}^{\mathrm{pr}}$ is the full $\mathrm{POSS}_L$ and all halving parameters are 0. It is clear that every $q \in \mathbb{Q}$ is compatible with $p_0$, i.e. $\mathbf{1} \leq^* p_0$. Hence $q \leq^* p_0$ for any $q \in \mathbb{Q}$, but if $i_0 \notin \operatorname{supp} q$, then we could easily construct some name $\dot{\tau}$ for a natural number that could be essentially read by $p_0$ (e.g. $\dot{\tau}$ could be (2,0)-decided by $p_0$), but not by $q$.

Continuous reading allows to estimate the size of the continuum in any generic extension.

**Lemma 7.5.** *In $V$, let $B \subseteq S^*$ and let $\kappa := \max\{\aleph_0, |B|\}^{\aleph_0}$. Then, in any $\mathbb{Q}$-generic extension, there are at most $\kappa$ many reals which are continuously read by some condition in the generic set and only using indices in $B$. To be more precise, define the following $\mathbb{Q}$-name for a subset of $2^\omega$:*

$$\dot{R}(B) := \{(\dot{r}, p) \mid \dot{r} \textit{ is a (nice) } \mathbb{Q}\textit{-name for a real in } 2^\omega,\ p \in \mathbb{Q} \textit{ and } \\ \dot{r} \textit{ is continuously read by } p \textit{ only using indices in } B\}.$$

*Then $\Vdash |\dot{R}(B)| \leq \kappa$, i.e. if $G$ is $\mathbb{Q}$-generic over $V$, then, in $V[G]$, there are at most $|\kappa|$ many reals $x \in 2^\omega$ such that there are $p \in G$ and a $\mathbb{Q}$-name $\dot{r} \in V$ such that $x = \dot{r}[G]$ and $V \vDash$"$p$ continuously reads $\dot{r}$ only using indices in $B$".*[14]

*Proof.* Let $P$ be the set of countable partial functions from $\mathrm{DOM} \cap (B \times \mathrm{ht})$ into $H_{\aleph_0}$ (the collection of hereditarily finite sets). Note that $|P| \leq \kappa^{\aleph_0} = \kappa$ (with equality in the case that $B \neq \varnothing$). Let $E$ be the set of countable partial functions from $P$ into 2; hence $|E| \leq |P|^{\aleph_0} \leq \kappa$ and $|E^\omega| \leq \kappa$.

Let $G$ be $\mathbb{Q}$-generic over $V$. In $V[G]$, if $x \in R := \dot{R}(B)[G]$, then we can choose some pair $(\dot{r}_x, p_x) \in \dot{R}(B)$ with $p_x \in G$ and $x = \dot{r}_x[G]$, and define $\bar{f}^x := \langle f^x_n \mid n < \omega\rangle \in E^\omega \cap V$ such that, in $V$, $f^x_n$ is the map that witnesses essential decision of $\dot{r}_x(n)$ by $p_x$ only using indices in $B$. More precisely, $f^x_n$ is a function with domain $D^x_n := \{\eta|_B \mid \eta \in \operatorname{poss}(p_x, {<}m_x(n))\}$ for some $m_x(n) \in \omega$ such that, for any $\eta \in \operatorname{poss}(p_x, {<}m_x(n))$, $p_x \wedge \eta \Vdash \dot{r}_x(n) = f^x_n(\eta|_B)$.

In $V[G]$, for any $n < \omega$ there is a unique $\eta^x_n \in \operatorname{poss}(p_x, {<}m_x(n))$ such that $p_x \wedge \eta^x_n \in G$, so $x(n) = \dot{r}_x[G](n) = f^x_n(\eta^x_n|_B)$. This implies that $x \mapsto \bar{f}^x$ is a one-to-one map from $R$ into $E^\omega \cap V$, so $|R| \leq |\kappa|$. In more detail, assume that $x, y \in R$ and $\bar{f}^x = \bar{f}^y$. For $n < \omega$, $D^x_n = D^y_n$, so $\eta^x_n|_B$ and $\eta^y_n|_B$ must have the same domain. On the other hand, $p_x \wedge \eta^x_n, p_y \wedge \eta^y_n \in G$, so $\eta^x_n$ and $\eta^y_n$ are compatible; hence $\eta^x_n|_B = \eta^y_n|_B$. Therefore $x(n) = f^x_n(\eta^x_n|_B) = f^y_n(\eta^y_n|_B) = y(n)$. ☐

As a direct consequence of Theorem 7.2, we get the following:

**Corollary 7.6.** *Let $\kappa := \max\{|S^*|, \aleph_0\}^{\aleph_0}$. Then $\mathbb{Q}$ forces that $\mathfrak{c} \leq |\kappa|$.*

7.1. **Rapid reading.** We show that some $\mathbb{Q}$-names for reals allow rapid reading, meaning that (as a consequence of Theorem 7.2) the set of conditions rapidly reading it is dense. We start with the following general result simplifying the proof of rapid reading, which also helps in the proof of continuous reading of names (hence Theorem 7.2 is not used in its proof).

[14] Note that $\kappa$ may be an ordinal in $V[G]$, but according to Corollary 7.3 it is still a cardinal if CH is assumed in the ground model.

**Lemma 7.7.** *Let $p \in \mathbb{Q}$ be a modest condition and $L_0 \in \omega$.*[15] *Assume that $\langle A_L \mid L \in \mathrm{Ht},\ L \leq L_0\rangle$ is a sequence of sets and that $\{\pi_{L',L} \mid L, L' \in \mathrm{Ht},\ L < L' \leq L_0\}$ is a collection of functions such that*

*(i) $|A_L| \leq n^P_{<L}$; and*
*(ii) for $L < L'$ in $\mathrm{Ht}$, $\pi_{L',L} \colon A_{L'} \to A_L$ such that, for any $L < L' < L''$ in $\mathrm{Ht}$, $\pi_{L'',L} = \pi_{L',L} \circ \pi_{L'',L'}$.*

*If $h_{L_0} \colon \mathrm{poss}(p, {<}L_0) \to A_{L_0}$, then there is some $q \leq p$ in $\mathbb{Q}$ and, for each $L < L_0$ in $\mathrm{Ht}$, there is some function $h_L \colon \mathrm{poss}(p, {<}L) \to A_L$ such that*

*(a) $p$ and $q$ have the same support, trunk lengths and halving parameters, and they are identical at levels ${\geq}L_0$;*
*(b) below $L_0$, the norm of each subatomic creature of $q$ decreases by at most 1 when compared to the subatomic creature in $p$ at the same position; and*
*(c) for $L < L_0$ in $\mathrm{Ht}$, if $\eta \in \mathrm{poss}(p, {<}L_0)$ and $\eta(\alpha, \ell) \in q(\alpha, \hat{\eta}(\alpha, n_*(\ell)), \ell)$ for all $\alpha \in \mathrm{supp}\, p$ and $L \leq \ell < L_0$ in $\mathrm{ht}$ (when $\alpha \in S^{\mathrm{pr}}$ or $\ell < \mathrm{trl}(p, \alpha)$, we commit abuse of notation and let $q(\alpha, \hat{\eta}(\alpha, n_*(\ell)), \ell) := q(\alpha, \ell)$), then $h_L(\eta{\restriction}_{L\downarrow}) = \pi_{L_0,L}(h_{L_0}(\eta))$.*

*Proof.* Fix a decreasing enumeration $\{L_j \mid 1 \leq j < m\}$ of $\mathrm{Ht} \cap L_0{\downarrow}$. Let $q_0 := p$ and $h_0 := h_{L_0}$, and write $A_j := A_{L_j}$ and $\pi_{j,k} := \pi_{L_j,L_k}$ for $j < k < m$. By induction on $j$, we construct $h_j = h_{L_j}$ and $q_j \leq q_{j-1}$ such that

(i) $q_{j-1}$ and $q_j$ have the same support, trunk lengths and halving parameters, and they are identical except at $L_j$ (and their sublevels); however, in the case that $L_j = (n, 0)$ for some $n < \omega$, they may also differ at levels $\ell \in I^*_n \cup \{L^{\mathrm{al}}_n\}$, but only in that the set of pr-indices of the atomic creatures in $q_j$ may be smaller than the set of pr-indices of the atomic creatures in $q_{j-1}$ (at the same coordinates);
(ii) at $L_j$, the norm of each subatomic creature of $q_j$ decreases by at most 1 when compared to the subatomic creature in $q_{j-1}$ at the same position; and
(iii) if $\eta \in \mathrm{poss}(p, {<}L_{j-1})$ and $\eta(\alpha, \ell) \in q_j(\alpha, \hat{\eta}(\alpha, n_*(\ell)), \ell))$ for all $\alpha \in \mathrm{supp}\, p$ and $L_j \leq \ell < L_{j-1}$ (with the same abuse of notation when $\alpha \in S^{\mathrm{pr}}$ or $\ell < \mathrm{trl}(p, \alpha)$), then $h_j(\eta{\restriction}_{L_j\downarrow}) = \pi_{j-1,j}(h_{j-1}(\eta))$.

By (i), notice that $q_j$ is identical to $p$ at levels ${<}L_j$ and ${\geq}L_0$.

Let $1 \leq j < m$ and assume that we have the desired $h_i$ and $q_i$ for $i < j$. We proceed by case distinctions depending on the height. When $L_j = L^{\mathrm{al}}_n$ for some $n < \omega$, we first assume that there is some $\alpha_j \in \mathrm{supp}\, p$ such that $p(\alpha_j, L_j)$ is non-trivial (which is unique by modesty). Recall that $q_{j-1}(\alpha_j, L_j) = p(\alpha_j, L_j)$. Fix $t \in \mathrm{pss}(p_{\mathrm{pr}}, \alpha_j, L_j)$ and write $B^t := \{\eta' \in \mathrm{poss}(p, {<}L_j) \mid \hat{\eta}'(\alpha_j, n) = t\}$.

Consider the function $F^t \colon p(\alpha_j, t, L_j) \to {A_j}^{B^t}$ such that $F^t(x) = F^t_x$ maps $\eta'$ to the value $\pi_{j-1,j}(h_{j-1}(\eta'^\frown x))$, where $\eta'^\frown x$ is the unique possibility $\eta \in \mathrm{poss}(p, {<}L_{j-1})$ extending $\eta'$ with $\eta(\alpha, L_j) = x$. (All other values at $L_j$ are trivially determined by modesty.)

Since both $B^t$ and $A_j$ have size ${\leq}n^P_{<L_j}$ (by Lemma 5.16), $|{A_j}^{B^t}| \leq (n^P_{<L_j})^{n^P_{<L_j}} < n^B_{L_j}$ by (pr2). Then, by bigness (Lemma 6.2 (b)), there is some subatom $q_j(\alpha_j, t, L_j) \leq$

[15] Recall that we identify $n \in \omega$ with the level $(n, 0)$.

$p(\alpha_j, t, L_j)$ such that $\|q_j(\alpha_j, t, L_j)\|^t_{L_j} \geq \|p(\alpha_j, t, L_j)\|^t_{L_j} - 1$ and $F^t{\restriction}_{q_j(\alpha_j,t,L_j)}$ is constant with value $f_t\colon B^t \to A_j$.

For $\eta' \in \mathrm{poss}(p, {<}L_j)$, we define $h_j(\eta') := f_{\hat{\eta}'(\alpha_j,n)}(\eta')$. We can define $q_j$ as in (i) such that $q_j(\alpha_j, L_j) = \langle q_j(\alpha_j, t, L_j) \mid t \in \mathrm{pss}(p_{\mathrm{pr}}, \alpha_j, L_j)\rangle$ and such that it coincides with $q_{j-1}$ (and $p$) at $(\alpha', L_j)$ for $\alpha' \neq \alpha_j$. Hence $q_j$ and $h_j$ are as required. Concretely, if $\eta \in \mathrm{poss}(p, {<}L_{j-1})$ and $\eta(\alpha_j, L_j) \in q_j(\alpha_j, \hat{\eta}(\alpha_j, n), L_j)$, then $h_j(\eta{\restriction}_{L_j\downarrow}) = \pi_{j-1,j}(h_{j-1}(\eta))$.

In the case that all creatures are trivial at $L_j$, we can just set $q_j := q_{j-1}$ and $h_j(\eta) := \pi_{j-1,j}(h_{j-1}(\eta^+))$ for $\eta \in \mathrm{poss}(p, < L_j)$, where $\eta^+$ is the unique possibility in $\mathrm{poss}(p, {<}L_{j-1})$ extending $\eta$.

Now consider the case that $L_j = L^{\mathrm{lc}}_n$. Enumerate $\mathrm{poss}(p, {<}L_j) = \{\eta'_k \mid k < m'\}$. By recursion, we define a decreasing sequence $\langle \mathbf{c}_k \mid k \leq m'\rangle$ of lc-compound creatures, all with the same support, halving parameter and pr-indices and such that $\mathbf{c}_0 := p(L^{\mathrm{lc}}_n)$. Assume we have defined $\mathbf{c}_k$. Define $f_k\colon \mathrm{poss}(\mathbf{c}_k, \bar{t}^k) \to A_j$, where $\bar{t}^k := \hat{\eta}'_k(L^{\mathrm{lc}}_n)$ and $f_k(\nu)$ is the value of $\pi_{j-1,j} \circ h_{j-1}$ evaluated for the unique possibility in $\mathrm{poss}(p, {<}L_{j-1})$ obtained from $\eta'_k$ and $\nu$.

By Corollary 6.4 applied to $f_k$, there is some $\mathbf{c}_{k+1} \leq \mathbf{c}_k$ as required and such that the norm of each subatom in $\mathbf{c}_{k+1}$ decreases by at most $1/n^B_{L^{\mathrm{lc}}_n}$ and $f_k$ has constant value $h_j(\eta'_k)$ on $\mathrm{poss}(\mathbf{c}_{k+1}, \bar{t}^k)$. This defines $h_j$, and we now define $q_j$ such that $q_j(L^{\mathrm{lc}}_n) := \mathbf{c}_{m'}$ and such that it is identical with $q_{j-1}$ at other levels. It is routine to check that they are as required.

In the case that $L_j = (n, 0)$, we construct $(q_j)_{\mathrm{pr}}(n)$ as in the previous case (which is simpler because the parameter $t$ is not required), but use bigness from Lemma 6.5. However, since $(q_j)_{\mathrm{pr}}(n)$ may be different from $(q_{j-1})_{\mathrm{pr}}(n) = p_{\mathrm{pr}}(n)$, we must also reduce the set of pr-indices of the atoms of $q_{j-1}$ at levels $\ell \in I^*_n \cup \{L^{\mathrm{al}}_n\}$ (as indicated in (i)). So define $q_{j,\mathrm{pr}} \in \mathbb{Q}_{\mathrm{pr}}$ identical to $(q_{j-1})_{\mathrm{pr}}$ except at $L_j$, where $q_{j,\mathrm{pr}}(n)$ is the constructed $(q_j)_{\mathrm{pr}}(n)$, and set $q_j := q_{j-1} \wedge q_{j,\mathrm{pr}}$. By Observation 5.18, $q_j$ is as required.

We have thus defined $h_j$ for all $j < m$. Set $q := q_{m-1}$, which is clearly as required. □

**Theorem 7.8.** *Let $\dot{r}$ be a $\mathbb{Q}$-name for an element of $\prod_{L\in\mathrm{Ht}} n^S_L$. Assume that $p \in \mathbb{Q}$ is modest and continuously reads $\dot{r}$. Then there is a $q \leq p$ with the same support, trunk lengths and halving parameters as $p$ which rapidly reads $\dot{r}$.*

*Proof.* For each $L \in \mathrm{Ht}$, set

$$(*_1) \qquad h(L) := \max\{L' \in \mathrm{Ht} \mid L' \leq L \text{ and } p \ L\text{-decides } \dot{r}{\restriction}_{L'\downarrow}\}$$

Note that $\langle h(L) \mid L \in \mathrm{Ht}\rangle$ is non-decreasing, and continuous reading implies that $h$ is an unbounded function. For $L' \leq L$, let $\dot{x}_{L,L'}$ be a $\mathbb{Q}$-name for $\dot{r}{\restriction}_{\min\{h(L),L'\}\downarrow}$. Note that $\dot{x}_{L,L'}$ is $L$-decided and that there are at most $|\prod_{L<L'} n^S_L|$ many possibilities for $\dot{x}_{L,L'}$.

Fix $L_0 \in \omega$. For all $L \leq L_0$ in Ht, we define

$$A_{L_0,L} := \prod\{n^S_{L'} \mid L' < \min\{h(L_0), L\},\ L' \in \mathrm{Ht}\}$$

and $\psi_{L_0,L_0}\colon \operatorname{poss}(p,{<}L_0) \to A_{L_0,L_0}$ such that $p \wedge \eta \Vdash \dot{x}_{L_0,L_0} = \psi_{L_0,L_0}(\eta)$. By (pr14) and Lemma 7.7 applied to these objects and to the projections $\pi_{L,L'}\colon A_{L_0,L} \to A_{L_0,L'}$ for $L' < L \leq L_0$, we can find $q_{L_0} \leq p$ and $\psi_{L_0,L}\colon \operatorname{poss}(p,{<}L) \to A_{L_0,L}$ for each $L < L_0$ in Ht such that

(i) $p$ and $q_{L_0}$ have the same support, trunk lengths and halving parameters, and they are identical at levels $\geq L_0$;
(ii) below $L_0$, the norm of each subatomic creature of $q_{L_0}$ decreases by at most 1 when compared to the subatomic creature in $p$ at the same position; and
(iii) for $L < L_0$ in Ht, if $\eta \in \operatorname{poss}(p,{<}L_0)$ and $\eta(\alpha,\ell) \in q_{L_0}(\alpha,\hat{\eta}(\alpha,n_*(\ell)),\ell)$ for all $\alpha \in \operatorname{supp} p$ and $L \leq \ell < L_0$ in ht (again allowing abuse of notation), then $p \wedge \eta \Vdash \dot{x}_{L_0,L} = \psi_{L_0,L}(\eta{\restriction}_{L\downarrow})$.

With regard to (iii), we say that $\eta{\restriction}_{L\downarrow}$ *decides* $\dot{x}_{L_0,L}$ *(*$= \psi_{L_0,L}(\eta{\restriction}_{L\downarrow})$*) modulo* $q_{L_0}{\restriction}_{[L,L_0)}$.

Note that, given $L_0 \in \omega$, there are only finitely many possibilities for $q_{L_0}{\restriction}_{L_0\downarrow}$ (restriction on ht) and $\psi_{L_0,L}$ for $L \leq L_0$ in Ht. Thus, by König's Lemma, there is some $q \leq p$ in $\mathbb{Q}$ and a sequence $\langle \psi^*_L \mid L \in \mathrm{Ht}\rangle$ such that

$$(*_2)\quad \begin{array}{l}\text{for all } L \in \mathrm{Ht}, \text{ there is some } L_0 \geq L \text{ in } \omega \text{ such that } q_{L_0} \text{ and } q \text{ are identical}\\ \qquad\text{below } L \text{ and } \psi_{L_0,L'} = \psi^*_{L'} \text{ for all } L' \leq L \text{ in Ht.}\end{array}$$

Here $q$ is constructed in the natural way, and it is clear that the norm of each subatom of $q$ decreases by at most 1 when compared with the subatom in $p$ at the same position, hence $q \in \mathbb{Q}$. Note that $q$ has the same support, trunk lengths, and halving parameters as $p$, and $q \leq p$.

To complete the proof, it remains to show that $q$ rapidly reads $\dot{r}$, that is, each $\eta \in \operatorname{poss}(q,{<}L)$ decides $\dot{r}{\restriction}_L$ for all $L \in \mathrm{Ht}$. Fix $L \in \mathrm{Ht}$ and choose an $L' \geq L$ in Ht such that $h(L') \geq L$. According to Eq. ($*_1$), we get that $\dot{r}{\restriction}_L$ is $L'$-decided. Choose $L_0 \geq L'$ as in Eq. ($*_2$) and recall (from (iii)) that $\dot{x}_{L_0,L}$ is decided by $\psi_{L_0,L}$ modulo $q_{L_0}{\restriction}_{[L,L_0)}$. Note that $L \leq h(L') \leq L' \leq L_0$, so $\min\{h(L_0), L\} = L$. Therefore, $\Vdash \dot{x}_{L_0,L} = \dot{r}{\restriction}_{L\downarrow}$ and, since $h(L') \geq L$, $\dot{x}_{L_0,L}$ is already $L'$-decided (by the original condition $p$). So we can decide $\dot{x}_{L_0,L}$ modulo $q_{L_0}{\restriction}_{[L,L')}$ using $\psi_{L_0,L}$.

On the other hand, $q$ and $q_{L_0}$ coincide below $L'$ and $\psi_{L',L} = \psi_{L_0,L} = \psi^*_L$, so $q \wedge \eta$ forces $\dot{r}{\restriction}_L = \dot{x}_{L_0,L} = \psi^*_L(\eta)$ for each $\eta \in \operatorname{poss}(q,{<}L)$. □

As a direct consequence of Lemma 5.8, Lemma 7.4 and Theorem 7.2, we obtain:

**Corollary 7.9.** *If* $\dot{r}$ *is as in Theorem 7.8 then the set of modest conditions rapidly reading* $\dot{r}$ *is dense in* $\mathbb{Q}$.

A consequence of rapid reading is that $\mathbb{Q}$ does not add random reals.

**Lemma 7.10.** *The forcing* $\mathbb{Q}$ *does not add random reals. In particular, under* CH, $\mathbb{Q}$ *forces* $\operatorname{cov}(\mathcal{N}) = \aleph_1$.

*Proof.* Let $\dot{r}$ be a $\mathbb{Q}$-name for a real in $2^\omega$ and let $p \in \mathbb{Q}$. For $k < \omega$, we set $n_k = \lfloor \log_2 n^S_{(k,0)} \rfloor$, and we define $\dot{f}\colon \mathrm{Ht} \to \omega$ by

$$\dot{f}(L) = \begin{cases} \dot{r}{\restriction}_{n_k} & \text{if } L = (k,0),\\ 0 & \text{otherwise.}\end{cases}$$

By Corollary 7.9, there is some modest $q \leq p$ rapidly reading $\dot{f}$. Then $\dot{r}\restriction_{n_k}$ is $(k,0)$-decided by $q$, that is, $\dot{r}\restriction_{n_k}$ is determined by $\eta \in \operatorname{poss}(q,{<}k)$. We now denote by $A^q_k$ the set of possible values of $\dot{r}\restriction_{n_k}$ determined by $\operatorname{poss}(q,{<}k)$; then $\frac{|A^q_k|}{2^{n_k}} \leq \frac{n^P_{<k}}{2^{n_k}}$ by Lemma 5.16. The sequence $\langle A^q_k \mid k < \omega\rangle$ allows us to define $N := \{x \in 2^\omega \mid \forall\, k < \omega \colon x\restriction_{n_k} \in A^q_k\}$ (in the ground model). It is clear from the definition of $N$ that $\operatorname{Lb}(N) = \lim_{k\to\infty} \frac{|A^q_k|}{2^{n_k}} \leq \lim_{k\to\infty} \frac{n^P_{<k}}{2^{n_k}} = 0$ (where Lb denotes the Lebesgue measure) and $q \Vdash \dot{r} \in N$. □

**Observation 7.11.** Under CH, by Theorem 7.2 and Lemma 7.10, $\mathbb{Q}$ is $\omega^\omega$-bounding and does not add random or Cohen reals. However, $\mathbb{Q}$ does not have the Laver property; if that were the case (together with $\omega^\omega$-bounding), it would have the Sacks property, which would imply $\Vdash_{\mathbb{Q}} \operatorname{cof}(\mathcal{N}) = \aleph_1$, contradicting Figure 2 and the results of Subsection 8.1 (i. e. $\mathbb{Q}$ increases localisation and anti-localisation cardinals).

7.2. **The proof of continuous reading and more.** In this subsection, we prove Theorem 7.2 and present more features of the forcing (such as fusion). This is the only part of the paper where the halving parameters are actually used. Our presentation is based on [GK21, Section 7].

Recall that the norm of an lc-compound creature $\mathbf{c}$ with non-empty support is defined by

$$\|\mathbf{c}\|^{\mathrm{lc}}_n := \frac{\log_2(\max\{1, \min\{\|\mathbf{c}(\alpha,t_0)\|^{t_0}_{\mathrm{stk}} \mid \alpha \in \operatorname{supp}\mathbf{c},\ t_0 \in P_{\mathbf{c},\alpha}\} - d_{\mathbf{c}}\})}{n^P_{<L^{\mathrm{lc}}_n}}.$$

Let $D_{\mathbf{c}} := \min\{\|\mathbf{c}(\alpha,t_0)\|^{t_0}_{\mathrm{stk}} \mid \alpha \in \operatorname{supp}\mathbf{c},\ t_0 \in P_{\mathbf{c},\alpha}\}$. If we change $d_{\mathbf{c}}$ to

$$d'_{\mathbf{c}} := d_{\mathbf{c}} + \frac{D_{\mathbf{c}} - d_{\mathbf{c}}}{2} = \frac{D_{\mathbf{c}} + d_{\mathbf{c}}}{2},$$

then the norm of the resulting lc-compound creature decreases by at most $1/n^P_{<L^{\mathrm{lc}}_n}$. We call this procedure *halving*, which can also be applied to pr-compound creatures.

**Definition 7.12** ([FGKS17, Definition 5.2.2])**.** Given a pr- or lc-compound creature $\mathbf{c}$, we define $\operatorname{hf}(\mathbf{c})$, *the half of* $\mathbf{c}$, to be the same creature as $\mathbf{c}$ except that we replace each halving parameter $d_{\mathbf{c}}$ by the $d'_{\mathbf{c}}$ described above. In the case that $\operatorname{supp}\mathbf{c} = 0$, we set $d'_{\mathbf{c}} = d_{\mathbf{c}} = 0$, i. e. $\operatorname{hf}(\mathbf{c}) = \mathbf{c}$.

By extending this procedure, we can also halve entire conditions.

**Definition 7.13.** Given a condition $q \in \mathbb{Q}$ and $m < \omega$, we define $r := \operatorname{hf}(q,m)$ to be the same condition as $q$ except that all compound creatures $q_{\mathrm{pr}}(n)$ and $q(L^{\mathrm{lc}}_n)$ for $n \geq m$ are halved, i. e. $r_{\mathrm{pr}}(n) = \operatorname{hf}(q_{\mathrm{pr}}(n))$ and $r(L^{\mathrm{lc}}_n) = \operatorname{hf}(q(L^{\mathrm{lc}}_n))$.

We also have a converse procedure to *unhalve* conditions, illustrated in the following result.

**Lemma 7.14.** *Let $M \in \mathbb{R}$, $p \in \mathbb{Q}$ and $N \in \omega$, $N \geq \operatorname{trl}(p)$ be such that all compound creatures of $p$ at levels $\geq N$ have norm $\geq M$. Assume that $r \leq \operatorname{hf}(p,N)$ in $\mathbb{Q}$ such that $\operatorname{trl}(r) = N$ and all compound creatures of $r$ at levels $\geq N$ have norm $>0$. Then there is a condition $q \leq p$ and an $N^* > N$ such that*

*(i)* $q$ *is identical to* $r$ *except for the halving parameters of the compound creatures at heights in* $[(N,0),(N^*,0))$*;*

*(ii)* $\|q_{\mathrm{pr}}(n)\|_n^{\mathrm{pr}} \geq M$ *for* $n \geq N^*$*;*

*(iii)* $\|q(L_n^{\mathrm{lc}})\|_n^{\mathrm{lc}} \geq M$ *for* $n \geq N^*$*;*

*(iv)* $\|q_{\mathrm{pr}}(n)\|_n^{\mathrm{pr}} \geq M - {}^1/{n^P_{<(n,0)}}$ *for all* $N \leq n < N^*$*; and*

*(v)* $\|q(L_n^{\mathrm{lc}})\|_n^{\mathrm{lc}} \geq M - {}^1/{n^P_{<L_n^{\mathrm{lc}}}}$ *for all* $N \leq n < N^*$*.*

*Moreover,* $q =^* r$*, so by Lemma 7.4,* $r$ *essentially decides some* $\mathbb{Q}$*-name* $\dot\tau$ *for a ground model object iff* $q$ *does.*

*Proof.* Choose $N^* > N$ such that $\|r_{\mathrm{pr}}(n)\|_n^{\mathrm{pr}} \geq M$ and $\|r(L_n^{\mathrm{lc}})\|_n^{\mathrm{lc}} \geq M$ for $n \geq N^*$. We define $q$ to be identical to $r$ except that, for all $N \leq n < N^*$, we replace the halving parameters $d_{r_{\mathrm{pr}}(n)}$ and $d_{r(L_n^{\mathrm{lc}})}$ by $d_{p_{\mathrm{pr}}(n)}$ and $d_{p(L_n^{\mathrm{lc}})}$, respectively; that is, $d_{q_{\mathrm{pr}}(n)} := d_{p_{\mathrm{pr}}(n)}$ and $d_{q(L_n^{\mathrm{lc}})} := d_{p(L_n^{\mathrm{lc}})}$. It is clear that $r \leq q$ and $q \leq^* r$, so $q =^* r$.

It is straightforward to see that $q \leq p$ and (i)–(iii) hold. We only show (v) (since (iv) is quite similar). Assume $N \leq n < N^*$. If $\operatorname{supp} q(L_n^{\mathrm{lc}}) = \varnothing$, then $q(L_n^{\mathrm{lc}}) = r(L_n^{\mathrm{lc}})$ and $p(L_n^{\mathrm{lc}})$ has empty support, too. But $\|r(L_n^{\mathrm{lc}})\| > 0$, so $\operatorname{supp} r_{\mathrm{lc}} = \varnothing$, which implies that $\operatorname{supp} p_{\mathrm{lc}} = \varnothing$ and $\|q(L_n^{\mathrm{lc}})\| = \|r(L_n^{\mathrm{lc}})\| = \|p(L_n^{\mathrm{lc}})\| = n \geq M$.

So, assuming $\operatorname{supp} q(L_n^{\mathrm{lc}}) \neq \varnothing$, it remains to show that

$$\|q(L_n^{\mathrm{lc}})\|_n^{\mathrm{lc}} = \frac{\log_2(\min\{\|q(\alpha,t_0)\|_{\mathrm{stk}}^{t_0} \mid (\alpha,t_0) \in Q_n\} - d_{q(L_n^{\mathrm{lc}})})}{n^P_{<L_n^{\mathrm{lc}}}} \geq M - \frac{1}{n^P_{<L_n^{\mathrm{lc}}}},$$

where $Q_n = \{(\alpha,t) \in \operatorname{supp} q(L_n^*) \times T_n^* \mid t \in \operatorname{pss}(q,\alpha,L_n^{\mathrm{lc}})\}$ (since $\|r(L_n^{\mathrm{lc}})\| > 0$, the term inside the logarithm is certainly $>1$).

Recall that $0 < \|r(L_n^{\mathrm{lc}})\|_n^{\mathrm{lc}} \leq \|q(L_n^{\mathrm{lc}})\|_n^{\mathrm{lc}}$. Fix some $(\alpha,t_0) \in Q_n$. Since $q$ and $r$ only differ in the halving parameters, we have

$$0 < \log_2(\|q(\alpha,t_0)\|_{\mathrm{stk}}^{t_0} - d_{r(L_n^{\mathrm{lc}})}),$$

which implies that

$$\|q(\alpha,t_0)\|_{\mathrm{stk}}^{t_0} > d_{r(L_n^{\mathrm{lc}})} + 1 \geq d_{\mathrm{hf}(p(L_n^{\mathrm{lc}}))} + 1 = d_{p(L_n^{\mathrm{lc}})} + \frac{D_{p(L_n^{\mathrm{lc}})} - d_{p(L_n^{\mathrm{lc}})}}{2} + 1.$$

Since $d_{q(L_n^{\mathrm{lc}})} = d_{p(L_n^{\mathrm{lc}})}$,

$$\|q(\alpha,t_0)\|_{\mathrm{stk}}^{t_0} - d_{q(L_n^{\mathrm{lc}})} > \frac{D_{p(L_n^{\mathrm{lc}})} - d_{p(L_n^{\mathrm{lc}})}}{2} + 1$$

for any $(\alpha,t_0) \in Q_n$, and hence

$$\frac{\log_2(\|q(\alpha,t_0)\|_{\mathrm{stk}}^{t_0} - d_{q(L_n^{\mathrm{lc}})})}{n^P_{<L_n^{\mathrm{lc}},0)}} \geq \|p(L_n^{\mathrm{lc}}))\|_n^{\mathrm{lc}} - \frac{1}{n^P_{<L_n^{\mathrm{lc}}}}.$$

Therefore

$$\|q(L_n^{\mathrm{lc}})\|_n^{\mathrm{lc}} \geq \|p(L_n^{\mathrm{lc}})\|_n^{\mathrm{lc}} - \frac{1}{n^P_{<L_n^{\mathrm{lc}}}} \geq M - \frac{1}{n^P_{<L_n^{\mathrm{lc}}}},$$

which proves (iv). □

Before engaging in the proof of continuous reading, similar to [GK21, Section 7] we first present some orders that give us a variant of the finite version of Baumgartner's axiom A for our forcing $\mathbb{Q}$,[16] from which it immediately follows that $\mathbb{Q}$ is proper and $\omega^\omega$-bounding and has continuous reading of names.

**Definition 7.15.**

(1) Let $N < \omega$. We define the order $\leq_N$ on $\mathbb{Q}$ as follows: $q \leq_N p$ iff
  (i) $\mathrm{trl}(q) = \mathrm{trl}(p) \leq N$;
  (ii) $q \leq p$;
  (iii) at each level $<N$ (excluding $\mathrm{ht}^{\mathrm{al}}$), the compound creatures of both conditions have the same support and halving parameter (note, however, that $\operatorname{supp} q$ could be larger than $\operatorname{supp} p$);
  (iv) $\mathrm{trl}(q, \alpha) \geq N$ for $\alpha \in \operatorname{supp} q_{\mathrm{al}} \setminus \operatorname{supp} p_{\mathrm{al}}$; and
  (v) $q(\alpha, \ell) = p(\alpha, \ell)$ for any $(\alpha, \ell) \in \mathrm{DOM}$ with $\alpha \in \operatorname{supp} p$ and $\ell < (N, 0)$.

  Properties (i) and (ii) imply that $\mathrm{trl}(q, \alpha) = \mathrm{trl}(p, \alpha)$ for all $\alpha \in \operatorname{supp} p_{\mathrm{al}}$ (by the definition of the order of $\mathbb{Q}$).

(2) Let $n < \omega$ and let $F \subseteq S^{\mathrm{al}}$ be finite. We define the order $\leq_{n,F}$ on $\mathbb{Q}$ as follows: $q \leq_{n,F} p$ if there is some $n' \geq n$ such that
  (i) $q \leq_{n'} p$;
  (ii) every compound creature at any level $\geq (n', 0)$ has norm $>n$; and
  (iii) for each $\alpha \in F \cap \operatorname{supp} p$, there is some $\ell < (n', 0)$ in $\mathrm{ht}^{\mathrm{al}}$ such that $\|p(\alpha, \ell)\| > n$.

(3) We say that a sequence $\bar{p} = \langle p_n \mid n < \omega \rangle$ of elements of $\mathbb{Q}$ is a *fusion sequence* if there is some sequence $\bar{F} = \langle F_n \mid n < \omega \rangle$ of finite subsets of $S^{\mathrm{al}}$ such that
  (F1) $F_n \subseteq F_{n+1}$;
  (F2) $p_{n+1} \leq_{n,F_n} p_n$; and
  (F3) $\bigcup_{n<\omega} F_n = \bigcup_{n<\omega} \operatorname{supp}(p_n)_{\mathrm{al}}$.

  We then say that $\bar{F}$ *witnesses the fusion sequence* $\bar{p}$.

It is easy to check that $\leq_{n,F}$ is a preorder on $\mathbb{Q}$.[17] Moreover, whenever $n_0 \leq n$ and $F_0 \subseteq F$, $q \leq_{n,F} p$ implies $q \leq_{n_0,F_0} p$; and whenever $\mathrm{trl}(p) \leq n_0 \leq n$, $q \leq_n p$ implies $q \leq_{n_0} p$. These properties of the preorders $\leq_{n,F}$ along with *fusion* (Lemma 7.16) and *pure decision* (Lemma 7.17) define our variation of Baumgartner's axiom A.

We can always construct a *fusion condition* from a fusion sequence.

**Lemma 7.16.** *If $\langle p_n \mid n < \omega \rangle$ is a fusion sequence witnessed by $\bar{F} = \langle F_n \mid n < \omega \rangle$, then there is a condition $q \in \mathbb{Q}$ satisfying* $\operatorname{supp} q = \bigcup_{n<\omega} \operatorname{supp} p_n$ *and $q \leq_{n,F_n} p_n$ for all $n < \omega$.*

*Proof.* For each $n$, let $h(n)$ be like $n'$ in Definition 7.15 (2) for $p_{n+1} \leq_{n,F_n} p_n$. Choose an increasing sequence $\langle n_j \mid j < \omega \rangle$ of natural numbers with $n_0 = 0$ such that $\langle h(n_j) \mid j < \omega \rangle$ is strictly increasing. Set $\mathrm{trl}(q) := \mathrm{trl}(p_0)$, $\operatorname{supp} q := \bigcup_{n<\omega} \operatorname{supp} p_n$ and, for $\alpha \in \operatorname{supp} q_{\mathrm{al}}$, $\mathrm{trl}(q, \alpha) := \mathrm{trl}(p_{n_j}, \alpha)$ for some $j$ such that $\alpha \in \operatorname{supp} p_{n_j}$. (This value does not depend on the chosen $j$.)

[16] See e. g. [MS19, Definition 1.2 and 1.3] for the definition of axiom A and its finite version (also known as "strong axiom A").

[17] Although the order $\leq_N$ is transitive, reflexivity fails for conditions with trunk length larger than $N$.

We fix $j < \omega$ and define $q$ at heights $\ell \in \mathrm{ht}_m$ for $h(n_{j-1}) \le m < h(n_j)$ as follows (with $h(n_{-1}) := 0$): We let $\operatorname{supp} q_{\mathrm{pr}}(m) := \operatorname{supp}(p_{n_j})_{\mathrm{pr}}(m)$, $\operatorname{supp} q(L^{\mathrm{lc}}_m) := \operatorname{supp} p_{n_j}(L^{\mathrm{lc}}_m)$, $d_{q_{\mathrm{pr}}(m)} := d_{(p_{n_j})_{\mathrm{pr}}(m)}$, $d_{q(L^{\mathrm{lc}}_m)} = d_{p_{n_j}(L^{\mathrm{lc}}_m)}$, and $q(\alpha,\ell) := p_{n_j}(\alpha,\ell)$ for all $\alpha \in \operatorname{supp} p_{n_j}$. Note that all of these objects are the same when $j$ is replaced by any $j' > j$. In the case that $\alpha \in \operatorname{supp} q \smallsetminus \operatorname{supp} p_{n_j}$, we let $q(\alpha,\ell) := p_{n_{j'}}(\alpha,\ell)$ for some $j'$ such that $\alpha \in \operatorname{supp} p_{n_{j'}}$. (Again, this value does not depend on the chosen $j'$.)

Routine calculations show $q \in \mathbb{Q}$ and $q \le_{n_j,F_{n_j}} p_{n_j}$ for all $j < \omega$. If $n \le n_j$, then $p_{n_j} \le_{n,F_n} p_n$, so $q \le_{n,F_n} p_n$. □

The following lemma concludes our variation of Baumgartner's axiom A.

**Lemma 7.17** (pure decision). *Assume that $\dot\tau$ is a $\mathbb{Q}$-name for a ground model object, that $p \in \mathbb{Q}$, that $N_0 \in \omega$ and that $M_0 > 1$ is a real such that any compound creature of $p$ at any level $\ge N_0$ has norm $\ge M_0 + 1$ (and hence $\mathrm{trl}(p) \le N_0$). Then there is some $q \in \mathbb{Q}$ such that*

*(i) $q \le_{N_0} p$;*
*(ii) $q$ essentially decides $\dot\tau$; and*
*(iii) any compound creature of $q$ at any level $\ge N_0$ has norm $\ge M_0$.*

*In particular, for any $n < \omega$ and any finite $F \subseteq \operatorname{supp} p_{\mathrm{al}}$, there is some $q \le_{n,F} p$ in $\mathbb{Q}$ essentially deciding $\dot\tau$.*

*Proof.* This proof consists of three parts.

**Part 1: Halving, the single step.**

Suppose that we are given $p \in \mathbb{Q}$, $N \in \omega$ and $M \in \mathbb{R}$ such that $N \ge \mathrm{trl}(p)$ and any compound creature of $p$ at any level $\ge N$ has norm $\ge M + 1$. We will show how to construct a condition $r := r(p,N,M) \in \mathbb{Q}$ satisfying

(H1) $r \le_N p$;
(H2) any compound creature of $r$ at any level $\ge N$ has norm $\ge M$; and
(H3) if $\eta \in \mathrm{poss}(r, {<}N)$ and if there is an $s \le r \wedge \eta$ such that $s$ essentially decides $\dot\tau$, $\mathrm{trl}(s) = N$ and the compound creatures of $s$ at any level $\ge N$ have norm $>0$, then $r \wedge \eta$ already essentially decides $\dot\tau$.

To see this, first enumerate $\mathrm{poss}(p, {<}N)$ as $\eta^1, \dots, \eta^m$. Let $q^0 := p$. By induction on $m$ we construct conditions $\hat q^1, q^1, \dots, \hat q^m, q^m$ such that, for each $j < m$,

(1) $\hat q^{j+1}$ is derived from $q^j$ by replacing the creature at each $(\alpha,\ell) \in \mathrm{DOM} \cap (\operatorname{supp} p \times (N,0){\downarrow})$ by the trivial creature $\{\eta^{j+1}(\alpha,\ell)\}$;
(2) $q^{j+1} \le \hat q^{j+1}$;
(3) $\mathrm{trl}(\hat q^{j+1}) = \mathrm{trl}(q^{j+1}) = N$;
(4) $\|q^{j+1}_{\mathrm{pr}}(n)\|^{\mathrm{pr}}_n \ge M + 1 - \frac{j+1}{n^P_{<(n,0)}}$ for $n \ge N$;
(5) $\|q^{j+1}(L^{\mathrm{lc}}_n)\|^{\mathrm{lc}}_n \ge M + 1 - \frac{j+1}{n^P_{<L^{\mathrm{lc}}_n}}$ for $n \ge N$; and
(6) one of the following two cases holds:
- Either $q^{j+1}$ essentially decides $\dot\tau$, or
- $q^{j+1} = \mathrm{hf}(\hat q^{j+1}, N)$.

Note that $\hat{q}^1$ is just $p \wedge \eta^1$, but $q^1$ may have support larger than $\operatorname{supp} p$, so we can not simply write $q^j \wedge \eta^{j+1}$ in (1). The definition of $\hat{q}^{j+1}$ is correct because $\operatorname{trl}(q^j) = N$, so $q^j(\alpha, \ell)$ is trivial for $(\alpha, \ell)$ in $\mathrm{DOM} \times (\operatorname{supp} q^j \times (N, 0){\downarrow})$ with $\alpha \notin \operatorname{supp} p$. The condition $q^{j+1}$ is easy to get: If there is some $q^{j+1}$ essentially deciding $\dot{\tau}$ and satisfying (2)–(5), we just choose it; otherwise we let $q^{j+1} = \operatorname{hf}(\hat{q}^{j+1}, N)$ (which clearly satisfies (2)–(6)). We thus always give priority to choosing an essentially deciding $q^{j+1}$ over halving $\hat{q}^{j+1}$.

We now define $r$ as follows: On $\mathrm{DOM} \cap (\operatorname{supp} p \times (N, 0){\downarrow})$, $r$ is identical to $p$; otherwise $r$ is identical to $q^m$. In more detail, we define

- $\operatorname{supp} r := \operatorname{supp} q^m$;
- $\operatorname{trl}(r) = \operatorname{trl}(p)$, $\operatorname{trl}(r, \alpha) := \operatorname{trl}(p, \alpha)$ for each $\alpha \in \operatorname{supp} p_{\mathrm{al}}$, and $\operatorname{trl}(r, \alpha) := \operatorname{trl}(q_m, \alpha)$ for each $\alpha \in \operatorname{supp} r_{\mathrm{al}} \setminus \operatorname{supp} p_{\mathrm{al}}$ (so the trunk lengths are $\geq N$);
- for $(\alpha, \ell) \in \mathrm{DOM}$ with $\ell < N$,
  - if $\alpha \in \operatorname{supp} p$, then $r(\alpha, \ell) := p(\alpha, \ell)$;
  - if $\alpha \in \operatorname{supp} r \setminus \operatorname{supp} p$, then $r(\alpha, \ell) := q^m(\alpha, \ell)$;
- for $(\alpha, \ell) \in \mathrm{DOM}$ with $\alpha \in \operatorname{supp} r$ and $\ell \geq N$, $r(\alpha, \ell) := q^m(\alpha, \ell)$; and
- the supports and halving parameters of compound creatures in $r$ coincide with those of compound creatures in $p$ at levels $<N$, and with those of compound creatures in $q^m$ at levels $\geq N$.

It is clear by the construction that $r$ satisfies (H1). To see (H2), note that, for $n \geq N$, $\|r(L_n^{\mathrm{lc}})\| = \|q_m(L_n^{\mathrm{lc}})\| \geq M + 1 - \frac{m}{n^P_{<L_n^{\mathrm{lc}}}} \geq M$, where the final inequality holds because $m \leq n^P_{<L_n^{\mathrm{lc}}}$ by Lemma 5.16.

It remains to prove (H3). Let $\eta \in \operatorname{poss}(r, <N)$ and $s$ as in (H3). Note that $\eta$ extends some $\eta^j \in \operatorname{poss}(p, <N)$, so $s \leq r \wedge \eta \leq q^j \leq \hat{q}^j$. It suffices to show that $q^j$ was constructed using the "decision" case. Assume instead towards a contradiction that the "halving" case was used. Then $s$ is stronger than $\operatorname{hf}(\hat{q}^j, N)$, so we can use Lemma 7.14 and unhalve $s$ to get some $s' \leq \hat{q}^j$, with $\|s'_{\mathrm{pr}}(n)\|_n^{\mathrm{pr}} \geq M + 1 - \frac{j}{n^P_{<(n,0)}}$ and $\|s'(L_n^{\mathrm{lc}})\|_n^{\mathrm{lc}} \geq M + 1 - \frac{j}{n^P_{<L_n^{\mathrm{lc}}}}$ for $n \geq N$, such that $s' =^* s$. But now, by Lemma 7.4, we could have used the "decision" case after all by choosing $q^j := s'$, which is a contradiction.

**Part 2: Iterating the single step.**

Given $p$, $N_0$ and $M_0$ as in the lemma's statement, we inductively construct conditions $p_j$ and natural numbers $N_j$ for $j < \omega$, starting with $N_0$. We define $p_0 := r(p, N_0, M_0)$. Given $p_j$ and $N_j$, we define $p_{j+1}$ and $N_{j+1}$ as follows:

- Choose $N_{j+1} > N_j$ such that
  - $\|(p_j)_{\mathrm{pr}}(n)\|_n^{\mathrm{pr}} \geq M_0 + j + 2$ for $n \geq N_{j+1}$;
  - $\|p_j(L_n^{\mathrm{lc}})\|_n^{\mathrm{lc}} \geq M_0 + j + 2$ for $n \geq N_{j+1}$; and
  - for all $m, i \leq j$, there is some $\ell \in \mathrm{ht}^{\mathrm{al}}$ with $(N_j, 0) < \ell < (N_{j+1}, 0)$ such that $\|p_j(\alpha_{m,i}, \ell)\| \geq M_0 + j + 1$, where $\operatorname{supp}_{\mathrm{al}}(p_j) := \{\alpha_{j,i} \mid i < \omega\}$.
- Let $p_{j+1} := r(p_j, N_{j+1}, M_0 + j + 1)$.

Thus $\langle p_j \mid j < \omega \rangle$ is a fusion sequence converging to a condition $q \in \mathbb{Q}$ constructed as in the proof of Lemma 7.16 (with $n_j = j$, $h(n_j) = N_j$ and $F_n = \{\alpha_{j,i} \mid j, i \leq n\}$).

By the construction, $\|q_{\mathrm{pr}}(n)\|_n^{\mathrm{pr}} \geq M_0 + j$ and $\|q(L_n^{\mathrm{lc}})\|_n^{\mathrm{lc}} \geq M_0 + j$ for $n \geq N_j$. On the other hand, $\operatorname{supp}_{\mathrm{al}} q = \bigcup_{j<\omega} \operatorname{supp}_{\mathrm{al}} p_j = \{\alpha_{m,i} \mid m, i < \omega\}$, so, for each $m, i < \omega$, if $j := \max\{m, i\}$, then there is an $\ell \in \mathrm{ht}^{\mathrm{al}}$ with $(N_j, 0) < \ell < (N_{j+1}, 0)$ such that $\|q(\alpha_{m,i}, \ell)\| = \|p_j(\alpha_{m,i}, \ell)\| \geq M_0 + j + 1$. This guarantees $p_{j+1} \leq_{N_j} p_j$ and $p_{j+1} \leq_{j,F_j} p_j$ as well as properties (i) and (iii) of the lemma.

It remains to show essential decision for $q$. The following property will be crucial for this proof:

$(*_3)$ If $j < \omega$, $\eta \in \operatorname{poss}(q, {<}N_j)$, $r \leq q \wedge \eta$ essentially decides $\dot\tau$, $\operatorname{trl}(r) = N_j$, $\|r_{\mathrm{pr}}(n)\|_n^{\mathrm{pr}} > 0$ and $\|r(L_n^{\mathrm{lc}})\|_n^{\mathrm{lc}} > 0$ for all $n \geq N_j$, then $q \wedge \eta$ already essentially decides $\dot\tau$.

Let $j$, $\eta$ and $r$ be as in Eq. ($*_3$). Then $r \leq q\wedge\eta \leq p_j \wedge \eta'$ for some $\eta' \in \operatorname{poss}(p_j, {<}N_j)$ by Fact 5.15 (c), hence by (H3), $p_j \wedge \eta'$ already essentially decides $\dot\tau$, and so does $q \wedge \eta$ by Lemma 7.4.

**Part 3: Bigness, thinning out.**

We use Lemma 7.7 to homogenise with regard to whether $q\wedge\eta$ essentially decides $\dot\tau$ or not. For $L_0 \in \omega$ and $L \leq L_0$ in Ht, let $A_{L_0,L} := 2$ and $h_{L_0,L_0} \colon \operatorname{poss}(q, {<}L_0) \to 2$ such that $h_{L_0,L_0}(\eta) = 1$ if $q \wedge \eta$ essentially decides $\dot\tau$. Let $B_{L_0,L_0}$ be the set of $\eta \in \operatorname{poss}(q, {<}L_0)$ such that $h_{L_0,L_0}(\eta) = 1$. By Lemma 7.7 applied to these objects and to the identity functions $\pi_{L,L'} \colon A_{L_0,L} \to A_{L_0,L'}$ for $L' < L \leq L_0$ in Ht, we find $q_{L_0} \leq q$ and, for $L < L_0$ in Ht, a function $h_{L_0,L} \colon \operatorname{poss}(q, {<}L) \to 2$ satisfying

(i) $q_{L_0}$ and $q$ have the same support, trunk lengths and halving parameters, and they are identical at any level $\geq L_0$;
(ii) below $L_0$, the norm of each subatomic creature from $q_{L_0}$ decreases by at most 1 when compared to the subatomic creature in $q$ at the same position; and
(iii) for $L < L_0$ in Ht, if $\eta \in \operatorname{poss}(q, {<}L_0)$ and $\eta(\alpha, \ell) \in q_{L_0}(\alpha, \hat\eta(\alpha, n_*(\ell)), \ell)$ for all $\alpha \in \operatorname{supp} q$ and $L \leq \ell < L_0$ (yet again allowing for some abuse of notation), then $h_{L_0,L}(\eta{\restriction}_{L\downarrow}) = h_{L_0,L_0}(\eta)$.

Let $B_{L_0,L} := \{\eta \in \operatorname{poss}(q, {<}L) \mid h_{L_0,L}(\eta) = 1\}$ for $L \leq L_0$ in Ht. Then (iii) states that, for any $\eta$ as in (iii), $\eta{\restriction}_{L\downarrow} \in B_{L_0,L}$ iff $q \wedge \eta$ essentially decides $\dot\tau$.

Given $L \leq L_0$ and $\eta \in \operatorname{poss}(q, {<}L)$, if $q \wedge \eta$ essentially decides $\dot\tau$ and $\eta' \in \operatorname{poss}(q, {<}L_0)$ extends $\eta$, then $q \wedge \eta'$ essentially decides $\dot\tau$, as well. Therefore, by (iii):

$(*_4)$ If $L \in \mathrm{Ht}$, $\eta \in \operatorname{poss}(q, {<}L)$ and $q \wedge \eta$ essentially decides $\dot\tau$, then $\eta \in B_{L_0,L}$ for each $L_0 \geq L$ in $\omega$.

The converse also holds (in some sense):

$(*_5)$
For $j < \omega$ and $L_0 \in \omega$, if $L_0 \geq N_j$ and $\eta \in B_{L_0,N_j}$ then $q \wedge \eta$ essentially decides $\dot\tau$.

To see this, we find a condition $r$ identical to $q_{L_0}$ at levels in $[N_j, L_0)$ and identical to $q \wedge \eta$ at the remaining levels. Any $\eta' \in \operatorname{poss}(r, {<}L_0) \subseteq \operatorname{poss}(q, {<}L_0)$ is in $B_{L_0,L_0}$ by (iii), so $q \wedge \eta'$ essentially decides $\dot\tau$, and so does $r \wedge \eta'$; hence $r$ also essentially decides $\dot\tau$. Moreover, $\|r_{\mathrm{pr}}(n)\|_n^{\mathrm{pr}} > 0$ and $\|r(L_{n^*}^{\mathrm{lc}})\|_{n^*}^{\mathrm{lc}} > 0$ for all $n \geq N_j$ (because $M_0 - 1 > 0$), so by Eq. ($*_3$), $q \wedge \eta$ essentially decides $\dot\tau$.

Note that $q$ essentially decides $\dot\tau$ iff, for some $L \in \mathrm{Ht}$, $q \wedge \eta$ essentially decides $\dot\tau$ for all $\eta \in \mathrm{poss}(q, {<}L)$. So, by Eq. ($*_5$), it suffices to prove that, for all $\eta \in \mathrm{poss}(q, {<}N_0)$, there is some $L_0 \geq N_0$ in $\omega$ such that $\eta \in B_{L_0,N_0}$.

As in Theorem 7.8, since for fixed $L_0 \in \omega$ there are only finitely many possibilities for $q_{L_0}\restriction_{L_0\downarrow}$ and $B_{L_0,L}$ for $L \leq L_0$ in Ht, by König's Lemma there is some $q^* \leq q$ and some sequence $\langle B_L^* \mid L \in \mathrm{Ht},\ L \geq N_0\rangle$ such that, for any $L \geq N_0$ in Ht, there is some $L_0 \geq L$ in $\omega$ such that $q^*$ is identical to $q$ below $N_0$ and identical to $q_{L_0}$ at levels in $[N_0, L]$, and $B_{L'}^* = B_{L_0,L'}$ for all $N_0 \leq L' \leq L$ in Ht.

Now fix $\eta \in \mathrm{poss}(q, {<}N_0)$. We find some $r \leq q^* \wedge \eta$ deciding $\dot\tau$. Without loss of generality, for some $m$, $\mathrm{trl}(r) = N_m$, $\|r(n)\|_n^{\mathrm{pr}} \geq 1$ and $\|r(L_n^{\mathrm{lc}})\|_n^{\mathrm{lc}} \geq 1$ for all $n \geq N_m$. Let $\eta'$ be the unique possibility in $\mathrm{poss}(r, {<}N_m)$ restricted to $\mathrm{supp}(q)$ which ensures $\eta' \in \mathrm{poss}(q, {<}N_m)$ and $r \leq q \wedge \eta'$. Then, by Eq. ($*_3$), $q \wedge \eta'$ already essentially decides $\dot\tau$.

Finally, pick some $L_0 \geq N_m$ in $\omega$ such that $q^*$ and $q_{L_0}$ are identical below $N_m$, and $B_L^* = B_{L_0,L}$ for all $N_0 \leq L \leq N_m$ in Ht. According to Eq. ($*_4$), $\eta' \in B_{L_0,N_m} = B_{N_m}^*$, so $\eta = \eta'\restriction_{N_0\downarrow} \in B_{N_0}^* = B_{L_0,N_0}$ because $\eta'(\alpha,\ell) \in q^*(\alpha, \hat\eta(\alpha, n_*(\ell)), \ell) = q_{L_0}(\alpha, \hat\eta(\alpha, n_*(\ell)), \ell)$ for all $\alpha \in \mathrm{supp}\, q$ and $\ell \in [N_0, N_m)$ in ht. ☐

We can now finally prove Theorem 7.2 using pure decision. To this end, the theorem is reformulated as follows:

**Theorem 7.18.** *Let $\dot r$ be a $\mathbb{Q}$-name for a function from $\omega$ into the ground model.*

*(a) For each $p \in \mathbb{Q}$, there is some fusion sequence $\langle p_n \mid n < \omega\rangle$ such that $p_0 \leq p$ and each $p_n$ essentially decides $\dot r(n)$.*
*(b) The set of conditions continuously reading $\dot r$ is dense in $\mathbb{Q}$.*
*(c) $\mathbb{Q}$ is $\omega^\omega$-bounding.*
*(d) $\mathbb{Q}$ is proper.*

*Proof.* To prove (a): By induction on $n < \omega$, we want to find a sequence $\bar F = \langle F_n \mid n < \omega\rangle$ of finite subsets of $S^{\mathrm{al}}$ such that

(i) $F_n \subseteq F_{n+1}$;
(ii) $p_{n+1} \leq_{n,F_n} p_n$;
(iii) $\bigcup_{n<\omega} F_n = \bigcup_{n<\omega} \mathrm{supp}(p_n)_{\mathrm{al}}$; and
(iv) $p_n$ essentially decides $\dot r(n)$.

We find some $p_0 \leq p$ in $\mathbb{Q}$ deciding $\dot r(0)$. We now assume that $p_n$ has already been defined and set $\mathrm{supp}(p_m)_{\mathrm{al}} = \{\alpha_{m,j} \mid j < \omega\}$ for all $m \leq n$. Next, we define $F_n := \{\alpha_{i,j} \mid i, j \leq n\}$. By applying Lemma 7.17, we can find $p_{n+1} \leq_{n,F_n} p_n$ essentially deciding $\dot r(n+1)$.

For (b) and (c): Given $p \in \mathbb{Q}$, we choose, using (a), a fusion sequence $\langle p_n \mid n < \omega\rangle$ such that $p_0 \leq p$ and such that each $p_n$ essentially decides $\dot r(n)$. Next, by Lemma 7.16, there is a $q \in \mathbb{Q}$ satisfying $\mathrm{supp}\, q = \bigcup_{n<\omega} F_n$ and $q \leq_{n,F_n} p_n$ for all $n < \omega$. This implies that $q \leq p$ and $q$ continuously read $\dot r$.

In addition, if $\dot{r}$ is a name for a function in $\omega^\omega$, then, since $\mathrm{poss}(q, {<}L)$ is finite for all $L \in \mathrm{Ht}$, we can find some $f \in \omega^\omega$ (in the ground model) such that $q \Vdash \dot{r}(n) \leq f(n)$, proving (c).

To see (d), let $\chi$ be a sufficiently large regular cardinal, let $p \in \mathbb{Q}$, and let $N \preceq H_\chi$ countable containing (as elements) $p$ and all the parameters of the forcing $\mathbb{Q}$. Let $\langle A_n \mid n \in \omega\rangle$ enumerate all maximal antichains in $\mathbb{Q}$ belonging to $N$, and choose a $\mathbb{Q}$-name $\dot{\tau}_n \in N$ for the element of $A_n$ chosen by the generic. As in (a), we can construct a fusion sequence $\langle p_n \mid n < \omega\rangle$ such that $p_n \in N$ and such that $p_n$ essentially decides $\dot{\tau}_n$. (Note, however, that the entire sequence may not be in $N$.) For each $n < \omega$, we can find (as in the argument for (c)) a finite set $E_n \in N$ with $E_n \subseteq A_n$ such that $p_n \Vdash \dot{\tau}_n \in E_n$. Hence $q \Vdash \dot{\tau}_n \in E_n \cap \dot{G} \subseteq A_n \cap N \cap \dot{G}$ (where $\dot{G}$ is the $\mathbb{Q}$-name for the generic set). □

## 8. Proving the Main Theorem

In this section, we finally prove Theorem 1.4. To this end, we require several assumptions:

**Assumption 8.1.** For what follows, we will make the following assumptions:

(1) CH.
(2) $|S^{\mathrm{pr}}| = \mu = \mu^{\aleph_0}$ is infinite.
(3) For each $i \in S^{\mathrm{pr}}$, $\kappa_i := |S^{\mathrm{lc}}_i| = |S^{\mathrm{al}}_i|$ is infinite and $\kappa_i = \kappa_i^{\aleph_0} \leq \mu$.

Note that (2) and (3) imply that $|S^*| = \mu$.

We present several results that, together, prove Theorem 1.4. We do not need to assume all of Assumption 8.1 (1)–(3) for all of them; we will present each result with the minimal necessary assumptions, i. e. we will not make use of the entirety of Assumption 8.1 all of the time.

First of all, recall from Theorem 7.2 that $\mathbb{Q}$ is proper and $\omega^\omega$-bounding and, under CH, $\mathbb{Q}$ has $\aleph_2$-cc by Lemma 5.23, so it preserves all cofinalities and cardinalities. On the other hand, under CH, $\mathbb{Q}$ forces $\mathrm{cov}(\mathcal{N}) = \aleph_1$ by Lemma 7.10.

According to Lemma 5.26 (b), at each $i \in S^{\mathrm{pr}}$ we add a real generic $y^*_i \in \prod_{n<\omega} T^*_n$ and the sequences $\vec{\nu}^{\mathrm{lc}}_i$ and $\vec{\nu}^{\mathrm{al}}_i$ are the generic blocks added at $i \in S^{\mathrm{pr}}$. This proves (II) of Theorem 1.4.

In the rest of this section, we will prove the remaining parts of Theorem 1.4, i. e. we will prove that, under Assumption 8.1, $\mathbb{Q}$ forces, for all $i \in S^{\mathrm{pr}}$,

(I) $\mathfrak{c} = \mu$;
(III) $\mathfrak{v}^\forall_{\dot{b}^{\mathrm{lc}}_i, \dot{h}^{\mathrm{lc}}_i} = \mathrm{cov}(\mathcal{I}_{\dot{f}^{\mathrm{lc}}_i}) = \mathfrak{c}^\exists_{\dot{a}^{\mathrm{lc}}_i, \dot{d}^{\mathrm{lc}}_i} = \mathfrak{c}^\forall_{\dot{a}^{\mathrm{lc}}_i, \dot{d}^{\mathrm{lc}}_i} = \kappa_i$; and
(IV) $\mathfrak{v}^\exists_{\dot{a}^{\mathrm{al}}_i, \dot{h}^{\mathrm{al}}_i} = \mathrm{non}(\mathcal{I}_{\dot{f}^{\mathrm{al}}_i}) = \mathfrak{c}^\forall_{\dot{a}^{\mathrm{al}}_i, \dot{h}^{\mathrm{al}}_i} = \kappa_i$.

To prove (III) and (IV), by Lemma 2.12 it suffices to show that $\mathbb{Q}$ forces, for all $i \in S^{\mathrm{pr}}$, $\kappa_i \leq \mathfrak{v}^\forall_{\dot{b}^{\mathrm{lc}}_i, \dot{h}^{\mathrm{lc}}_i}, \mathfrak{v}^\exists_{\dot{a}^{\mathrm{al}}_i, \dot{h}^{\mathrm{al}}_i}$, and $\mathfrak{c}^\forall_{\dot{a}^{\mathrm{lc}}_i, \dot{d}^{\mathrm{lc}}_i}, \mathfrak{c}^\forall_{\dot{a}^{\mathrm{al}}_i, \dot{h}^{\mathrm{al}}_i} \leq \kappa_i$. In Subsection 8.1, we prove that the $\kappa_i$ are lower bounds by using the generic slaloms.

In Subsection 8.2, we force (I) and $\mathfrak{c}^{\forall}_{\dot{a}^{\text{lc}}_i,\dot{d}^{\text{lc}}_i}, \mathfrak{c}^{\forall}_{\dot{a}^{\text{al}}_i,\dot{h}^{\text{al}}_i} \leq \kappa_i$, basically by forcing that the reals in $\dot{R}(\{i\} \cup S^{\text{lc}}_i \cup S^{\text{al}}_i)$ (see Lemma 7.5) give witnesses of both cardinals. A feature of the forcing called *separated support* is essential in the proofs.

8.1. **Forcing the lower bounds.** We start with the anti-localisation cardinals. Recall from Corollary 5.19 that, whenever $\alpha \in S^{\text{al}}$, $\mathbb{Q}_{S^*\smallsetminus\{\alpha\}}$ is a complete subforcing of $\mathbb{Q}$.

**Lemma 8.2.** *Let $i \in S^{\text{pr}}$, let $\alpha \in S^{\text{al}}_i$ and let $\dot{r}$ be a $\mathbb{Q}_{S^*\smallsetminus\{\alpha\}}$-name for a real in $\prod_{n<\omega} \dot{a}^{\text{al}}_i(n)$. Then $\Vdash_{\mathbb{Q}} \dot{r} \in^\infty \dot{\varphi}_\alpha$.*

*Proof.* We fix $p \in \mathbb{Q}$ and assume, without loss of generality, that $\alpha \in \operatorname{supp} p$. It is enough to show that, for any $n_0 < \omega$, there are $n \geq n_0$ and some $q \leq p$ forcing $\dot{r}(n) \in \dot{\varphi}_\alpha(n)$.

We pick some $n \geq n_0$ such that $\|p(\alpha, L^{\text{al}}_n)\|_{\text{al}} \geq 1$, and we define $\ell_n := L^{\text{al}}_n$. We now choose $t_0 \in \operatorname{pss}(p, \alpha, \ell_n)$ such that $\|p(\alpha, t_0, \ell_n)\|^{t_0}_{\ell_n} \geq 1$. We pick some $\eta \in \operatorname{poss}(p, {<}n+1)$ such that $\hat{\eta}(\alpha, n) = t_0$; then there are $k$ and $p' \leq p \wedge \eta$ forcing $\dot{r}(n) = k$. We define $q$ as being identical to $p'$ except at $(\alpha, \ell_n)$, where $q(\alpha, \ell_n) = \{w\}$ for some $w \in p(\alpha, t_0, \ell_n)$ containing $k$ (which exists by Observation 4.3 because $\|p(\alpha, t_0, \ell_n)\|^{t_0}_{\ell_n} \geq 1$). It is clear that $q \in \mathbb{Q}$, $q \leq p$ and $q \Vdash \dot{r}(n) \in w = \dot{y}_\alpha(L^{\text{al}}_n) = \dot{\varphi}_\alpha(n)$. □

As an immediate consequence, we get:

**Corollary 8.3.** *Assume CH. Then $\mathbb{Q}$ forces $|S^{\text{al}}_i| \leq \mathfrak{v}^{\exists}_{\dot{a}^{\text{al}}_i,\dot{h}^{\text{al}}_i}$ for all $i \in S^{\text{pr}}$.*

*Proof.* Let $\langle \dot{r}_\zeta \mid \zeta < \kappa\rangle$ be a sequence of $\mathbb{Q}$-names for elements of $\prod \dot{a}^{\text{al}}_i$ with $\kappa < |S^{\text{al}}_i|$ a cardinal. If $\kappa \leq \aleph_0$, it is clear that $\mathbb{Q}$ forces that, for some slalom $\varphi \in \mathcal{S}(\dot{a}^{\text{al}}_i, \dot{h}^{\text{al}}_i)$, $\dot{r}_\zeta \in^\infty \varphi$ for all $\zeta < \kappa$ (because $\mathfrak{v}^{\exists}_{a,h}$ is uncountable whenever $\lim_{n\to\infty} \frac{h(n)}{a(n)} = 0$). So assume that $\kappa$ is uncountable.

For each $\zeta < \kappa$, we can assume, without loss of generality, that $\dot{r}_\zeta$ is a nice name, i. e. that each $\dot{r}_\zeta(n)$ is determined by a maximal antichain $A_{\zeta,n} \subseteq \mathbb{Q}$. Then $B := \bigcup\{\operatorname{supp}(q) \mid q \in A_{\zeta,n},\ \zeta < \kappa,\ n < \omega\}$ has size $\leq \kappa < |S^{\text{al}}_i|$, because $\mathbb{Q}$ has $\aleph_2$-cc, and so we can pick some $\alpha \in S^{\text{al}}_i \smallsetminus B$.

It is clear that $\dot{r}_\zeta$ is a $\mathbb{Q}_{S^*\smallsetminus\{\alpha\}}$-name for all $\zeta < \kappa$. Therefore, by Lemma 8.2, $\Vdash_{\mathbb{Q}} \dot{r}_\zeta \in^\infty \dot{\varphi}_\alpha$, which finishes the proof. □

For the previous results we could, alternatively, use continuous reading. Namely, if instead of assuming in Lemma 8.2 that $\dot{r}$ is a $\mathbb{Q}_{S^*\smallsetminus\{\alpha\}}$-name, we assume that $p \in \mathbb{Q}$ continuously reads $\dot{r}$ without using the index $\alpha$, then we can show that there is some $q \leq p$ forcing $\dot{r} \in^\infty \dot{\varphi}_\alpha$. This is actually the approach we must use to increase the localisation cardinals, because we cannot say that $\mathbb{Q}_{S^*\smallsetminus\{\alpha\}}$ is a complete subforcing of $\mathbb{Q}$ whenever $\alpha \in S^{\text{lc}}$.

**Lemma 8.4.** *Let $i \in S^{\text{pr}}$, let $\alpha_0 \in S^{\text{lc}}_i$, let $\dot{r}$ be a $\mathbb{Q}$-name for a real in $\prod \dot{b}^-_i$, and let $p_0 \in \mathbb{Q}$. Assume that, for all $n < \omega$, $\dot{r}(n) := \langle \dot{r}(\ell) \mid \ell \in I^*_n\rangle$ is $L^{\text{al}}_n$-decided by $p_0$ without using the index $\alpha_0$. Then there is some $q \leq p_0$ forcing $\dot{r} \in^*_{\bar{I}^*} \dot{y}_{\alpha_0}$.*

*Proof.* We fix a modest $p \leq p_0$ in $\mathbb{Q}$ with $\alpha_0 \in \operatorname{supp} p$. We then pick some $n_0 < \omega$ such that $\alpha_0 \in \operatorname{supp} p(L^{\mathrm{lc}}_n)$ and $\|p(L^{\mathrm{lc}}_n)\|^{\mathrm{lc}}_n \geq 2$ for all $n \geq n_0$. We construct $q$ identical to $p$ except on the lc-part above $L^{\mathrm{lc}}_{n_0}$, hence we will now construct $q(L^{\mathrm{lc}}_n)$ by induction on $n \geq n_0$.

We fix some $n \geq n_0$. Assume that we have constructed $q$ up to $<L^{\mathrm{lc}}_n$ already. Since $\dot r(n)$ is $L^{\mathrm{al}}_n$-decided by $p$ without using the index $\alpha_0$, there is a function $F^n\colon \operatorname{poss}(p, <L^{\mathrm{al}}_n) \to b^*_{t^{\mathrm{mx}}_n}(L^{\mathrm{lc}}_n)$ (see (fp6)) that calculates $\dot r(n)$ and such that $F^n(\eta) = F^n(\eta')$ whenever $\eta|_{\operatorname{supp} p \setminus \{\alpha_0\}} = \eta'|_{\operatorname{supp} p \setminus \{\alpha_0\}}$. For each $\eta \in \operatorname{poss}(p, <L^{\mathrm{lc}}_n)$, we define $\bar t^\eta := \hat\eta(L^{\mathrm{lc}}_n) = \langle t^\eta_\alpha \mid \alpha \in \operatorname{supp} p(L^{\mathrm{lc}}_n)\rangle$ (so $t^\eta_\alpha := \hat\eta(\alpha, n)$) and $t^\eta_* := t^\eta_{\alpha_0}$. We also define $F^n_\eta\colon \operatorname{poss} p(L^{\mathrm{lc}}_n, \eta) \to b^*_{t^\eta_*}(L^{\mathrm{lc}}_n)$ such that $F^n_\eta(\bar x)$ is $F^n$ evaluated on the unique element of $\operatorname{poss}(p, <L^{\mathrm{al}}_n)$ constructed from $\eta$ and $\bar x$.

We enumerate $\operatorname{poss}(p, <L^{\mathrm{lc}}_n)$ as $\eta_0, \dots, \eta_{m-1}$. By induction on $k \leq m$, we construct compound lim inf creatures $\mathbf{d}_k \leq \mathbf{d}_{k-1} \leq p(L^{\mathrm{lc}}_n)$ and $\ell_k \in I^*_n$ for $k > 0$ such that the creatures have the same domain, halving parameters and pr-indices, and such that

$$\|\mathbf{d}_{k+1}(\alpha, t, \ell)\|^t_\ell \geq \|\mathbf{d}_k(\alpha, t, \ell)\|^t_\ell - \frac{1}{d_t(\ell)},$$

for each valid $(\alpha, t, \ell)$, except for $(\alpha_0, t^{\eta_k}_*, \ell_{k+1})$, where $\mathbf{d}_{k+1}(\alpha_0, t^{\eta_k}_*, \ell_{k+1})$ is a singleton. Let $\mathbf{d}_0 := p(L^{\mathrm{lc}}_n)$. Assume we have constructed $\mathbf{d}_k$ already. We choose $\ell_{k+1} \in I^*_n \setminus \{\ell_{k'} \mid 0 < k' \leq k\}$ such that $\|\mathbf{d}_k(\alpha_0, t^k_*, \ell_{k+1})\|^{t^k_*}_{\ell_{k+1}} \geq 1$ (this is possible because $\|\mathbf{d}_k\|^{\mathrm{lc}}_n > 1$). We define $\bar t^k := \bar t^{\eta_k}$, $t^k_* := t^{\eta_k}_*$ and $F^n_k := F^n_{\eta_k}$. Let $f^n_k\colon \operatorname{poss} p(L^{\mathrm{lc}}_n, \eta_k) \to b_{t^k_*}(\ell_{k+1})$, where $f^n_k(\bar x)$ is the $\ell_{k+1}$-th coordinate of $F^n_k(\bar x)$. According to Lemma 6.3 and (pr9), there is a compound creature $\mathbf{d}'_k \leq \mathbf{d}_k$ with the same domain, halving parameter and pr-indices such that $f^n_k\restriction_{\operatorname{poss}(\mathbf{d}'_k, \bar t^k)}$ only depends on

$$C_{k+1} := \{(\alpha, \ell) \in (\operatorname{supp} p(L^{\mathrm{lc}}_n)) \times I^*_n \mid (t^k_\alpha, \ell) \leq (t^k_*, \ell_{k+1})\}.$$

But recall that $\dot r(n)$ does not depend on $\alpha_0$, and that $\mathbf{d}_k(\alpha_0, \ell_{k+1})$ is not trivial, so $p(\alpha, \ell_{k+1})$ is trivial for all $\alpha \neq \alpha_0$ by strong modesty; thus the decision of $\dot r(n)$ does not depend on level $\ell_{k+1}$. Therefore, $f^n_k$ does not depend on $(\alpha_0, \ell_{k+1})$, so $f^n_k\restriction_{\operatorname{poss}(\mathbf{d}'_k, \bar t^k)}$ only depends on $C^*_{k+1} := \{(\alpha, \ell) \in (\operatorname{supp} p(L^{\mathrm{lc}}_n)) \times I^*_n \mid (t^k_\alpha, \ell) < (t^k_*, \ell_{k+1})\}$. Considering this, we can redefine $\mathbf{d}'_k(\alpha_0, t^k_*, \ell_{k+1}) := \mathbf{d}_k(\alpha_0, t^k_*, \ell_{k+1})$.

We define $H_k := \{f^n_k(\bar x) \mid \bar x \in \operatorname{poss}(\mathbf{d}'_k, \bar t^k)\}$. Then $|H_k| \leq \prod_{(\alpha,\ell)\in C^*_{k+1}} |\mathbf{d}'_k(\alpha, t^k_\alpha, \ell)| < d_{t^k_*}(\ell_{k+1})$ by (pr9). Since $\|\mathbf{d}'_k(\alpha_0, t^k_*, \ell_{k+1})\| \geq 1$, we have $\|\mathbf{d}'_k(\alpha_0, t^k_*, \ell_{k+1})\|^{\mathrm{cov}}_{\ell_{k+1}, t^k_*} > d_{t^k_*}(\ell_{k+1})$, so there is some $W_k \in \mathbf{d}'_k(\alpha_0, t^k_*, \ell_{k+1})$ such that $H_k \subseteq W_k$. We define $\mathbf{d}_{k+1}$ identical to $\mathbf{d}'_k$, except at $(\alpha_0, t^k_*, \ell_{k+1})$, where we define it as $\mathbf{d}_{k+1}(\alpha_0, t^k_*, \ell_{k+1}) := \{W_k\}$.

Define $q(L^{\mathrm{lc}}_n) := \mathbf{d}_m$. According to the construction, the norm of each subatom of $q(L^{\mathrm{lc}}_n)$ decreases by at most $n^P_{<L^{\mathrm{lc}}_n}/n^B_{L^{\mathrm{lc}}_n} \leq 1$ when compared to the norm of the subatom of $p(L^{\mathrm{lc}}_n)$ at the same coordinate, except at $(\alpha_0, t^k_*, \ell_{k+1})$ for $k < m$, where we get singletons. Hence, by Lemma 4.14 (b), $\|q(L^{\mathrm{lc}}_n)\|^{\mathrm{lc}}_n \geq \|p(L^{\mathrm{lc}}_n)\|^{\mathrm{lc}}_n - 1$.

We now prove that $p \wedge \eta \Vdash \exists \ell \in I^*_n\colon \dot r(\ell) \in \dot y_{\alpha_0}(\ell)$ for all $\eta \in \operatorname{poss}(q, <L^{\mathrm{al}}_n)$. (Although $q$ has not been fully defined at this point, we can talk about $\operatorname{poss}(q, <L^{\mathrm{al}}_n)$ because $q$ is already defined below $L^{\mathrm{al}}_n$.) Each $\eta \in \operatorname{poss}(q, <L^{\mathrm{al}}_n)$ depends on some $\eta_k \in \operatorname{poss}(p, {<}L^{\mathrm{lc}}_n)$ and $\bar{x} \in \operatorname{poss} q(L^{\mathrm{lc}}_n)(\bar{t}^{\eta_k})$, and by the construction of $W_k$, we know that $p \wedge \eta \Vdash \dot{r}(\ell_{k+1}) = f^n_k(\bar{x}) \in H_k \subseteq W_k = \dot{y}_{\alpha_0}(\ell_{k+1})$.

This finishes the construction, and it is clear that $q$ is in $\mathbb{Q}$ and that it is as required. □

Just as in Corollary 8.3, this implies:

**Corollary 8.5.** *Assume* CH. *Then, for all $i \in S^{\mathrm{pr}}$, $\mathbb{Q}$ forces $|S^{\mathrm{lc}}_i| \leq \mathfrak{v}^{\forall,\bar{I}^*}_{\dot{b}^-_i,\dot{h}^-_i} \leq \mathfrak{v}^{\forall}_{\dot{b}^{\mathrm{lc}}_i,\dot{h}^{\mathrm{lc}}_i}$.*

*Proof.* This proof is a variation of the proof of Corollary 8.3. Let $\kappa < |S^{\mathrm{lc}}_i|$ be a cardinal and let $\langle \dot{r}_\zeta \mid \zeta < \kappa \rangle$ be a sequence of names for elements of $\prod \dot{b}^-_i$. If $\kappa \leq \aleph_0$, it is clear that $\mathbb{Q}$ forces that, for some slalom $\varphi \in \mathcal{S}(\dot{b}^-_i, \dot{h}^-_i)$, $\dot{r}_\zeta \in^*_{\bar{I}^*} \varphi$ for all $\zeta < \kappa$. So assume that $\kappa$ is uncountable.

For each $\zeta < \kappa$, we define a $\mathbb{Q}$-name $\dot{r}'_\zeta$ for a function with domain Ht such that $\mathbb{Q}$ forces

$$\dot{r}'_\zeta(L) = \begin{cases} \langle \dot{r}_\zeta(\ell) \mid \ell \in I^*_n \rangle & \text{if } L = L^{\mathrm{lc}}_n, \\ 0 & \text{otherwise,} \end{cases}$$

and, by Corollary 7.9, we pick a maximal antichain $A_\zeta \subseteq \mathbb{Q}$ such that each condition in $A_\zeta$ rapidly reads $\dot{r}'_\zeta$. Then $B := \bigcup_{\zeta<\kappa} \bigcup_{p \in A_\zeta} \operatorname{supp}(p)$ has size $\leq \kappa < |S^{\mathrm{lc}}_i|$, so we can pick some $\alpha \in S^{\mathrm{lc}}_i \setminus B$.

It is enough to show that, for any $\zeta < \kappa$ and $p_0 \in A_\zeta$, $p_0 \Vdash \dot{r}_\zeta \in^*_{\bar{I}^*} \dot{y}_\alpha$. If $p \leq p_0$, then, by Lemma 7.4, $p$ rapidly reads $\dot{r}'_\zeta$ without using the index $\alpha$, so by Lemma 8.4, there is some $q \leq p$ forcing $\dot{r}_\zeta \in^*_{\bar{I}^*} \dot{y}_\alpha$. □

8.2. **Forcing the upper bounds.** To calculate the upper bounds, we need a property (called *separated support* in [KS09, KS12]) which basically says that the set of conditions $p \in \mathbb{Q}_{\mathrm{pr}}$ such that $\langle \hat{\eta}(i,n) \mid i \in \operatorname{supp} p(n) \rangle$ is a one-to-one sequence for all $n < \omega$ and $\eta \in \operatorname{poss}(p, {<}n+1)$ is dense. This is the main reason why we constructed the pr-part as a lim inf forcing.

Motivated by [KS12, Lemma 2.8], we start by proving a version of separated support for pr-compound creatures.

**Lemma 8.6.** *Let $\mathbf{c}$ be a pr-compound creature at level $n$ and let $i_0 \in \operatorname{supp} \mathbf{c}$. Then there is some pr-compound creature $\mathbf{d} \leq \mathbf{c}$ with the same domain and halving parameter such that $\|\mathbf{d}\|^{\mathrm{pr}}_n \geq \|\mathbf{c}\|^{\mathrm{pr}}_n - \frac{1}{n^B_{(n,0)}}$ and, for $i \neq i_0$ in $\operatorname{supp} \mathbf{c}$ and $L \in \mathrm{ht}^{\mathrm{pr}}_n$, if $\mathbf{d}(i,L)$ is non-trivial, then $\mathbf{d}(i,L) \cap \mathbf{d}(i_0,L) = \varnothing$.*

*Proof.* Let $d_{\mathbf{d}} := d_{\mathbf{c}}$, $\operatorname{supp} \mathbf{d} := \operatorname{supp} \mathbf{c}$, $\operatorname{dom} \mathbf{d} := \operatorname{dom} \mathbf{c}$ and, for each $i \in \operatorname{dom} \mathbf{c}$ and $L \in \mathrm{ht}^{\mathrm{pr}}_n$,

$$\mathbf{d}(i,L) = \begin{cases} \mathbf{c}(i,L) \setminus \mathbf{c}(i_0,L) & \text{if } i \neq i_0, \text{ and } \mathbf{c}(i,L) \text{ is non-trivial,} \\ \mathbf{c}(i,L) & \text{otherwise.} \end{cases}$$

In the first case, $\mathbf{c}(i_0,L)$ is a singleton by modesty, so $\|\mathbf{d}(i,L)\| \geq \|\mathbf{c}(i,L)\| - {}^{1}/_{n^B_L}$. It is clear from the definition above that $\mathbf{d}$ is as required. □

**Corollary 8.7.** *Let $\mathbf{c}$ be a pr-compound creature at level $n$. Then there is some pr-compound $\mathbf{d} \leq \mathbf{c}$ with the same domain and halving parameter such that $\|\mathbf{d}\|^{\mathrm{pr}}_n \geq$*

$\|\mathbf{c}\|_n^{\mathrm{pr}} - \frac{|\operatorname{supp}\mathbf{c}|}{n^B_{(n,0)}}$ *and, for any* $i \in \operatorname{supp}\mathbf{c}$ *and* $L \in \mathrm{ht}_n^{\mathrm{pr}}$*, if* $\mathbf{d}(i,L)$ *is non-trivial, then* $\mathbf{d}(i,L) \cap \mathbf{d}(i',L) = \varnothing$ *for all* $i' \in \operatorname{supp}\mathbf{c} \setminus \{i\}$*.*

This motivates the following notion.

**Definition 8.8.** We fix the following terminology.

(1) A pr-compound $\mathbf{c}$ at level $n$ has *separated support* if, for any $i \in \operatorname{supp}\mathbf{c}$ and $L \in \mathrm{ht}_n^{\mathrm{pr}}$, if $\mathbf{c}(i,L)$ is non-trivial, then $\mathbf{c}(i,L) \cap \mathbf{c}(i',L) = \varnothing$ for all $i' \in \operatorname{supp}\mathbf{c} \setminus \{i\}$.
(2) A condition $p \in \mathbb{Q}$ has *separated support* if it satisfies, for all $n < \omega$,
   (i) if $\|p_{\mathrm{pr}}(n)\| = 0$, then $\operatorname{supp} p_{\mathrm{pr}}(n) = \varnothing$;
   (ii) $p_{\mathrm{pr}}(n)$ has separated support; and
   (iii) if $\alpha \in \operatorname{supp} p_{\mathrm{al}}$ and $p(\alpha, L_n^{\mathrm{al}})$ is non-trivial, then $i^*(\alpha) \in \operatorname{supp} p_{\mathrm{pr}}(n)$.
   Note that property (iii) for $\alpha \in \operatorname{supp} p_{\mathrm{lc}}$ already holds by property (8) of Definition 5.4 (on the lc-part of the forcing). This property of the lc-part has not been used so far.

We can now easily obtain conditions in $\mathbb{Q}_{\mathrm{pr}}$ with separated support.

**Lemma 8.9.** *Let* $p \in \mathbb{Q}_{\mathrm{pr}}$*. Then there is some* $q \leq p$ *in* $\mathbb{Q}_{\mathrm{pr}}$ *with the same support and halving parameters satisfying separated support.*

*Proof.* We fix $n_0 \geq \mathrm{trl}(p)$ such that, for all $n \geq n_0$, $\|p(n)\| > 2$ and $|\operatorname{supp} p(n)| \leq n$ (which is possible by Definition 5.1 (8)–(9)). We extend the trunk length to $n_0$; concretely, we choose $\eta \in \mathrm{poss}(p, {<}n_0)$ and let $p^* := p \wedge \eta$.

We now define $q$ by cases: For $n < n_0$, we define $q(n) := p^*(n)$. For $n \geq n_0$, by Corollary 8.7 there is some pr-compound creature $q(n)$ at level $n$ with the same domain, support and halving parameter as $p^*(n)$ such that $q(n) \leq p^*(n)$ satisfies separated support and $\|q(n)\|_n^{\mathrm{pr}} \geq \|p^*(n)\|_n^{\mathrm{pr}} - \frac{|\operatorname{supp} p(n)|}{n^B_{(n,0)}} \geq \|p^*(n)\|_n^{\mathrm{pr}} - \frac{n}{n^B_{(n,0)}} > \|p^*(n)\|_n^{\mathrm{pr}} - 1$.

Finally, we define $q := \langle q(n) \mid n < \omega \rangle$, which is as required. □

**Lemma 8.10.** *If* $p \in \mathbb{Q}$*, then there is some modest* $q \leq p$ *with the same support and halving parameters satisfying separated support.*

*Proof.* By Lemma 8.9, there is some $s \leq p_{\mathrm{pr}}$ in $\mathbb{Q}_{\mathrm{pr}}$ with the same support and halving parameters satisfying separated support. Let $n_0 := \mathrm{trl}(s)$ and let $\eta'$ be the only possibility in $\mathrm{poss}(s, {<}n_0)$, choose some $\eta \in \mathrm{poss}(p, {<}n_0)$ extending $\eta'$ and let $r := (p \wedge \eta) \wedge s$ (see Observation 5.18).

It is clear that conditions (i) and (ii) of Definition 8.8 hold, so we strengthen $r$ to also get (iii). We only need to modify the al-part, so we define $q'$ to be identical to $r$ in the pr-part and the lc-part, with the same trunk lengths, and define $q'$ for each $\alpha \in \operatorname{supp} r_{\mathrm{al}}$ as follows: We can find some $n_\alpha \geq \mathrm{trl}(r, \alpha)$ such that $i^*(\alpha) \in \operatorname{supp} r_{\mathrm{pr}}(n_\alpha)$, so we define $q'(\alpha, L_n^{\mathrm{al}}) := r(\alpha, L_n^{\mathrm{al}})$ for each $n \geq n_\alpha$ and pick some trivial $q'(\alpha, L_n^{\mathrm{al}}) \leq r(\alpha, l_n^{\mathrm{al}})$ for each $n < n_\alpha$.

Finally, by Lemma 5.8, we can find a modest $q \leq q'$ with the same support and trunk lengths as $q'$ and with the same pr-part and lc-part as $q$. This resulting $q$ is as required. □

One direct consequence of separated support is that we can calculate the size of the continuum in Theorem 1.4.

**Lemma 8.11.** *Under* CH*, if* $|S^*| = |S^{\mathrm{pr}}| = \mu = \mu^{\aleph_0}$ *is infinite, then* $\mathbb{Q}$ *forces* $\mathfrak{c} = \mu$.[18]

*Proof.* By Corollary 7.6, $\mathbb{Q}$ forces $\mathfrak{c} \leq \mu$. On the other hand, Lemma 8.10 implies that, in any generic extension, $\langle y_i^* \mid i \in S^{\mathrm{pr}} \rangle$ is a one-to-one sequence of reals, so $\mathbb{Q}$ forces $\mu = |S^{\mathrm{pr}}| \leq \mathfrak{c}$. □

We use separated support to prove that $\mathbb{Q}$ *is* $(\dot{a}_i, \dot{d}_i)$*-bounding over* $\{i\} \cup S_i^{\mathrm{lc}} \cup S_i^{\mathrm{al}}$, which means that $\mathbb{Q}$ forces that any real in $\prod \dot{a}_i$ is localised by some slalom in $\mathcal{S}(\dot{a}_i, \dot{d}_i) \cap \dot{R}(\{i\} \cup S_i^{\mathrm{lc}} \cup S_i^{\mathrm{al}})$ (recall Lemma 7.5 for the notation).

**Lemma 8.12.** *Fix* $i \in S^{\mathrm{pr}}$ *and let* $\dot{r}$ *be a* $\mathbb{Q}$*-name for a function in* $\prod_{L \in \Omega} \dot{a}_i(L)$*. Then, for each* $p \in \mathbb{Q}$*, there are some* $q \leq p$ *in* $\mathbb{Q}$ *and some* $\mathbb{Q}$*-name* $\dot{\varphi}$ *for a slalom in* $\mathcal{S}(\dot{a}_i, \dot{d}_i)$ *such that*

*(i)* $q \Vdash \dot{r} \in^* \dot{\varphi}$*;*
*(ii) for* $n \in \omega$*,* $\dot{\varphi}(L_n^{\mathrm{lc}})$ *is* $L_n^{\mathrm{al}}$*-decided by* $q$ *only using indices in* $\{i\} \cup S_i^{\mathrm{lc}}$*; and*
*(iii)* $\dot{\varphi}(L_n^{\mathrm{al}})$ *is* $n+1$*-decided by* $q$ *only using indices in* $\{i\} \cup S_i^{\mathrm{al}}$*.*

*Proof.* Fix $i \in S^{\mathrm{pr}}$ and let $E_{\mathrm{al}} := \{i\} \cup S_i^{\mathrm{al}}$ and $E_{\mathrm{lc}} := \{i\} \cup S_i^{\mathrm{lc}}$. Assume that $\dot{r}$ is a $\mathbb{Q}$-name for an element of $\prod_{L \in \Omega} a_i(L)$ and let $p$ be a condition in $\mathbb{Q}$. Without loss of generality, we can assume that $i \in \operatorname{supp} p_{\mathrm{pr}}(\operatorname{trl}(p))$, that $p$ rapidly reads $\dot{r}'$ (by Theorem 7.8) and that $p$ is modest with separated support (by Lemma 8.10), where $\dot{r}'$ is a $\mathbb{Q}$-name for a function with domain Ht such that

$$\dot{r}'(L) := \begin{cases} \dot{r}(L) & \text{if } L \in \Omega, \\ 0 & \text{otherwise.} \end{cases}$$

We let $q_{\mathrm{pr}} := p_{\mathrm{pr}}$, $\operatorname{supp} q := \operatorname{supp} p$ and define $q$ to be identical to $p$ below $N_0 := \operatorname{trl}(q) = \operatorname{trl}(p)$. For $L < (N_0, 0)$ in $\Omega$, let $\dot{\varphi}(L)$ be the canonical name for $\varnothing$ (so it is $L$-decided without using any index at all). We now construct $\dot{\varphi}(L)$ and $q$ at level $L$ by induction on $L \in \Omega$, $L \geq N_0$. Fix $L \geq N_0$ in $\Omega$ and assume we have defined $q$ below $L$. We distinguish two cases:

**Case 1:** $L = L_n^{\mathrm{lc}}$. Let $\{\eta_k \mid k < m_n\}$ enumerate $\operatorname{poss}(q, {<}L_n^{\mathrm{lc}})$. For each $k < m_n$, we define

- $S_{n,k} := \{\alpha \in \operatorname{supp} p(L_n^{\mathrm{lc}}) \mid \hat{\eta}_k(\alpha, n) <_{\mathrm{lex}} \hat{\eta}_k(i, n)\}$; and
- $E_{n,k} := \{\alpha \in \operatorname{supp} p(L_n^{\mathrm{lc}}) \mid \hat{\eta}_k(\alpha, n) = \hat{\eta}_k(i, n)\}$.

Due to separated support (and using Definition 5.4 (8)), $E_{n,k} = \operatorname{supp} p(L_n^{\mathrm{lc}}) \cap S_i^{\mathrm{lc}} =: E_n$. On the other hand, by rapid reading, we can define $F_k \colon p(L_n^{\mathrm{lc}}, \eta_k) \to a_{\hat{\eta}_k(i,n)}(L_n^{\mathrm{lc}})$ such that $F_k(\bar{x})$ is the value of $\dot{r}(L_n^{\mathrm{lc}})$ decided by $p \wedge \eta'$, where $\eta'$ is the unique possibility in $\operatorname{poss}(p, {<}L_n^{\mathrm{al}})$ determined by $\eta_k$ and $\bar{x}$.

[18] Technically, we do not need to assume $|S^{\mathrm{pr}}| = |S^*|$. By a method similar to the proof of Lemma 8.6, we can force $\mathfrak{c}$ to be bigger than $|S^{\mathrm{pr}}|$ and $|S_i^{\mathrm{t}}|$ for all $i \in S^{\mathrm{pr}}$ and $\mathrm{t} \in \{\mathrm{lc}, \mathrm{al}\}$, which implies $\mathfrak{c} \geq |S^*|$.

By induction on $k < m_n$, we construct a decreasing sequence of lc-compound creatures $\langle \mathbf{d}_k \mid k < m_n \rangle$ stronger than $p(L_n^{\mathrm{lc}})$ with the same domain, halving parameters and pr-indices, and a $\mathbb{Q}$-name $\dot{s}_{n,k}$ for a subset of $\dot{a}_i(L_n^{\mathrm{lc}})$ as follows:

In step $k$, we define $\mathbf{d} := \mathbf{d}_{k-1}$, where $\mathbf{d}_{-1} := p(L_n^{\mathrm{lc}})$. By (pr10) and Lemma 6.3 applied to $(\hat{\eta}_k(i,n), \ell_n^{\mathrm{mx}})$, there is a compound creature $\mathbf{d}_k \leq \mathbf{d}$ with the same domain, halving parameter and pr-indices such that $F_k{\restriction}_{\mathrm{poss}(\mathbf{d}_k, \hat{\eta}_k(L_n^{\mathrm{lc}}))}$ depends only on $(S_{n,k} \cup E_n) \times I_n^*$. Denote by $F_k' \colon \prod_{(\alpha,\ell) \in (S_{n,k} \cup E_n) \times I_n^*} \mathbf{d}_k(\alpha, \hat{\eta}_k(\alpha, n), \ell) \to a_{\hat{\eta}_k(i,n)}(L_n^{\mathrm{lc}})$ the corresponding function that gives the same values as $F_k$.

We now define a $\mathbb{Q}$-name $\dot{s}_{n,k}$ such that, for each $q' \in \mathbb{Q}$ with $\mathrm{trl}(q') > n$ and $\{i\} \cup E_n \subseteq \operatorname{supp} q'$, if $q'(i, L) = \{\eta_k(i, L)\}$ for all $L \in \mathrm{ht}_n^{\mathrm{pr}}$ and $q'(\alpha, \ell) = \{v(\alpha, \ell)\} \subseteq \mathbf{d}_k(\alpha, \hat{\eta}_k(\alpha, \ell), \ell)$ for all $\alpha \in E_n$ and $\ell \in I_n^*$, then

$$q' \Vdash \dot{s}_{n,k} = \left\{ F_k'(\bar{u} \cup \langle v(\alpha,\ell) \mid (\alpha,\ell) \in E_n \times I_n^* \rangle) \;\middle|\; \bar{u} \in \prod_{(\alpha,\ell) \in S_k \times I_n^*} \mathbf{d}_k(\alpha, \hat{\eta}_k(\alpha, n), \ell) \right\};$$

otherwise, $q'$ forces $\dot{s}_{n,k} = \varnothing$. Clearly, $\mathbb{Q}$ forces $\dot{s}_{n,k} \subseteq \dot{a}_i(L_n^{\mathrm{lc}})$ and $|\dot{s}_{n,k}| \leq m^S_{\hat{\eta}_k(i,n)^-, \ell_n^{\mathrm{mx}}}$ if $\hat{\eta}_k(i,n)$ is not minimal in $T_n^*$; otherwise, $|\dot{s}_{n,k}| \leq 1$. It is also clear that, for any $q_1 \in \mathbb{Q}$ containing $\{i\} \cup E_n$ in its support, $\dot{s}_{n,k}$ is $L_n^{\mathrm{al}}$-decided by $q_1$ (and in particular by $p$) only using indices in $\{i\} \cup E_n \subseteq E_{\mathrm{lc}}$.

Let $q(L_n^{\mathrm{lc}}) := \mathfrak{d}_{m_n - 1}$. Due to the applications of Lemma 6.3, the norms of all subatoms in $q(L_n^{\mathrm{lc}})$ decrease by at most ${}^{n^P_{<L_n^{\mathrm{lc}}}}/_{n^B_{L_n^{\mathrm{lc}}}} < 1$ when compared to the subatoms of $p(L_n^{\mathrm{lc}})$, so we can conclude that $\|q(L_n^{\mathrm{lc}})\|_n^{\mathrm{lc}} \geq \|p(L_n^{\mathrm{lc}})\|_n^{\mathrm{lc}} - 1$. We define $\dot{\varphi}(L_n^{\mathrm{lc}})$ as a $\mathbb{Q}$-name for $\bigcup_{k<m_n} \dot{s}_{n,k}$. Note that $\mathbb{Q}$ forces that $|\dot{\varphi}(L_n^{\mathrm{lc}})| \leq n^P_{<L_n^{\mathrm{lc}}} \cdot m^S_{\dot{y}_i^*(n)^-, \ell_n^{\mathrm{mx}}} < d_{\dot{y}_i^*(n)}(L_n^{\mathrm{lc}}) = \dot{d}_i(L_n^{\mathrm{lc}})$ by (pr10) and Lemma 5.16. Moreover, $\dot{\varphi}(L_n^{\mathrm{lc}})$ is $L_n^{\mathrm{al}}$-decided by $p$ only using indices in $\{i\} \cup E_n$.

We prove that, for all $\eta \in \mathrm{poss}(q, {<}L_n^{\mathrm{al}})$ (which we can use because $q$ has been defined ${<}L_n^{\mathrm{al}}$), $p \wedge \eta \Vdash \dot{r}(L_n^{\mathrm{lc}}) \in \dot{\varphi}(L_n^{\mathrm{lc}})$. Each $\eta \in \mathrm{poss}(q, {<}L_n^{\mathrm{al}})$ depends on some $\eta_k \in \mathrm{poss}(p, {<}L_n^{\mathrm{lc}})$ and $\bar{x} \in \operatorname{poss} q(L_n^{\mathrm{lc}}, \hat{\eta}_k(L_n^{\mathrm{lc}}))$, so $p \wedge \eta \Vdash \dot{r}(L_n^{\mathrm{lc}}) = F_k(\bar{x}) = F_k'(\bar{x}{\restriction}_{((S_{n,k} \cup E_n) \times I_n^*)})$, which implies that $p \wedge \eta \Vdash \dot{r}(L_n^{\mathrm{lc}}) \in \dot{s}_{n,k}$. Hence $p \wedge \eta \Vdash \dot{r}(L_n^{\mathrm{lc}}) \in \dot{\varphi}(L_n^{\mathrm{lc}})$.

**Case 2:** $L = L_n^{\mathrm{al}}$. We first consider the case when all atomic creatures of $p$ at level $L_n^{\mathrm{al}}$ are trivial, so $\dot{r}(L_n^{\mathrm{al}})$ is $L_n^{\mathrm{al}}$-decided by rapid reading. We define $q$ to be identical to $p$ at this level. Let $F \colon \mathrm{poss}(q, {<}L_n^{\mathrm{al}}) \to a_{t_n^{\mathrm{mx}}}(L_n^{\mathrm{al}})$ such that $p \wedge \eta \Vdash \dot{r}(L_n^{\mathrm{al}}) = F(\eta)$; note that $|\operatorname{ran} F| \leq n^P_{<L_n^{\mathrm{al}}} < d_{t_n^{\mathrm{mn}}}(L_n^{\mathrm{al}})$. So we define $\dot{\varphi}(L_n^{\mathrm{al}})$ as a $\mathbb{Q}$-name for $\dot{a}_i(L_n^{\mathrm{al}}) \cap \operatorname{ran} F$, which is forced to have size ${<}\dot{d}_i(L_n^{\mathrm{al}})$. It is clear that $q \wedge \eta \Vdash \dot{r}(L_n^{\mathrm{al}}) = F(\eta) \in \dot{\varphi}(L_n^{\mathrm{al}})$ for all $\eta \in \mathrm{poss}(q, {<}L_n^{\mathrm{al}})$, and that $\dot{\varphi}(L_n^{\mathrm{al}})$ is $L_n^{\mathrm{lc}}$-decided by $q$ only using the index $i$.

Now, consider the case when $p$ contains some non-trivial atomic creature at $L_n^{\mathrm{al}}$, which means, due to modesty, that there is a unique $\alpha_n \in \operatorname{supp} p$ such that $p(\alpha_n, L_n^{\mathrm{al}})$ is non-trivial. For $\eta \in \mathrm{poss}(p, {<}L_n^{\mathrm{al}})$, we can define

$$F_\eta \colon p(\alpha_n, \hat{\eta}(\alpha_n, n), L_n^{\mathrm{al}}) \to a_{\hat{\eta}(i,n)}(L_n^{\mathrm{al}})$$

such that $F_\eta(s)$ is the value of $\dot r(L^{\mathrm{al}}_n)$ forced by $p$ restricted to the unique element of $\mathrm{poss}(p,{<}n+1)$ obtained from $\eta$ and $s$ (which is possible by rapid reading). Set $i_n := i^*(\alpha_n)$, which is in $\operatorname{supp} p_{\mathrm{pr}}(n)$ by Definition 8.8 (2) (iii).

To find $q(L^{\mathrm{al}}_n)$ and $\dot\varphi(L^{\mathrm{al}}_n)$, we further distinguish two subcases:

**Subcase 2.1:** $i = i_n$. In this case, $\alpha_n \in S^{\mathrm{al}}_i$. We define $q$ to be identical to $p$ at level $L^{\mathrm{al}}_n$ and also define $\dot\varphi(L^{\mathrm{al}}_n)$ as a $\mathbb{Q}$-name such that, for any $q' \in \mathbb{Q}$ with $\mathrm{trl}(q') > n$ and $i, \alpha_n \in \operatorname{supp} q'$, if $q'(i, \mathrm{ht}^{\mathrm{pr}}_n) = \{t\} \subseteq q(i, \mathrm{ht}^{\mathrm{pr}}_n)$ and $q'(\alpha_n, L^{\mathrm{al}}_n) = \{v\} \subseteq p(\alpha_n, t, L^{\mathrm{al}}_n)$, then

$$q' \Vdash \dot\varphi(L^{\mathrm{al}}_n) = \{F_\eta(v) \mid \eta \in \mathrm{poss}(q, {<}L^{\mathrm{al}}_n),\ \hat\eta(\alpha_n, n) = t\};$$

otherwise, $q'$ forces $\dot\varphi(L^{\mathrm{al}}_n) = \varnothing$. Note that $\dot\varphi(L^{\mathrm{al}}_n)$ is $n+1$-decided by $p$ only using indices in $\{i, \alpha_n\} \subseteq E_{\mathrm{al}}$.

It is easy to see that $p \wedge \eta' \Vdash$“$\dot r(L^{\mathrm{al}}_n) \in \dot\varphi(L^{\mathrm{al}}_n)$ and $|\dot\varphi(L^{\mathrm{al}})| \le n^P_{<L^{\mathrm{al}}_n} < \dot d_i(L^{\mathrm{al}}_n)$” for all $\eta' \in \mathrm{poss}(q, {<}n+1)$.

**Subcase 2.2:** $i \neq i_n$. Fix $t \in \mathrm{pss}(p, \alpha_n, L^{\mathrm{al}}_n)$ and enumerate $\{\eta \in \mathrm{poss}(q, {<}L^{\mathrm{al}}_n) \mid \hat\eta(i_n, n) = t\} = \{\eta^t_k \mid k < m^t\}$. By induction, we construct subatomic creatures $\mathbf{c}^t_{m^t} \le \cdots \le \mathbf{c}^t_1 \le \mathbf{c}^t_0 := p(\alpha_n, t, L^{\mathrm{al}}_n)$; so given $\mathbf{c}^t_k$, we let $t' := \hat\eta^t_k(i, n)$ and define $\mathbf{c}^t_{k+1}$ as follows: When $t' <_{\mathrm{lex}} t$, since $a_{t'}(L^{\mathrm{al}}_n) < n^S_{t', L^{\mathrm{al}}_n} < d_t(L^{\mathrm{al}}_n)$ by (pr8) and (pr12), we can apply Lemma 6.2 (b) to $\mathbf{c}^t_k$ and $F_{\eta^t_k}$ to get $\mathbf{c}^t_{k+1} \le \mathbf{c}^t_k$ such that $F_{\eta^t_k}\restriction_{\mathbf{c}^t_{k+1}}$ is constant with value $j^t_k \in a_{t'}(L^{\mathrm{al}}_n)$ and $\|\mathbf{c}^t_{k+1}\| \ge \|\mathbf{c}^t_k\| - \frac{1}{d_t(L^{\mathrm{al}}_n)}$; let $s^t_k := \{j^t_k\}$. In the case that $t' >_{\mathrm{lex}} t$, let $\mathbf{c}^t_{k+1} = \mathbf{c}^t_k$ and $s^t_k := F_{\eta^t_k}[\mathbf{c}^t_k]$, so $|s^t_k| \le n^S_{t, L^{\mathrm{al}}_n}$ by (pr8). (Note that the case $t' = t$ cannot happen because $p$ satisfies separated support.)

To define $q$ at level $L^{\mathrm{al}}_n$, let $q(\alpha_n, t, L^{\mathrm{al}}_n) := \mathbf{c}^t_{m^t}$ for all $t \in \mathrm{pss}(q, \alpha_n, L^{\mathrm{al}}_n)$ and $q(\alpha, L^{\mathrm{al}}_n) := p(\alpha, L^{\mathrm{al}}_n)$ for all $\alpha \neq \alpha_n$ in $\operatorname{supp} p_{\mathrm{al}}$. By construction, $\|q(\alpha_n, L^{\mathrm{al}}_n)\| \ge \|p(\alpha_n, L^{\mathrm{al}}_n)\| - \frac{n^P_{<L^{\mathrm{al}}_n}}{d_{t^{\mathrm{mn}}_n}(L^{\mathrm{al}}_n)} > \|p(\alpha_n, L^{\mathrm{al}}_n)\| - 1$. Define $\dot\varphi(L^{\mathrm{al}}_n)$ as a $\mathbb{Q}$-name such that, for all $q' \in \mathbb{Q}_i$ with $\mathrm{trl}(q') > n$ and $i \in \operatorname{supp} q'$, if $q'(i, \mathrm{ht}^{\mathrm{pr}}_n) = \{t_0\} \subseteq q(i, \mathrm{ht}^{\mathrm{pr}}_n)$, then

$$q' \Vdash \dot\varphi(L^{\mathrm{al}}_n) = \bigcup\{s^t_k \mid t \in \mathrm{pss}(p, \alpha_n, L^{\mathrm{al}}_n),\ k < m^t,\ \hat\eta^t_k(i, n) = t_0\};$$

otherwise $q'$ forces $\dot\varphi(L^{\mathrm{al}}_n) = \varnothing$. In the first case, $q'$ forces $|\dot\varphi(L^{\mathrm{al}}_n)| \le n^P_{<L^{\mathrm{al}}_n} \cdot (n^S_{t_0^-, L^{\mathrm{al}}_n} + 1) < d_{t_0}(L^{\mathrm{al}}_n) = \dot d_i(L^{\mathrm{al}}_n)$ by (pr12) (in the case that $t_0 = t^{\mathrm{mn}}_n$, simply replace $n^S_{t_0^-, L^{\mathrm{al}}_n}$ by 0). Therefore, $\mathbb{Q}$ forces that $|\dot\varphi(L^{\mathrm{lc}}_n)| < \dot d_i(L^{\mathrm{al}}_n)$. Moreover, $\dot\varphi(L^{\mathrm{al}}_n)$ is $L^{\mathrm{lc}}_n$-decided by $p$ only using the index $i$.

We show that $p \wedge \eta' \Vdash \dot r(L^{\mathrm{al}}_n) \in \dot\varphi(L^{\mathrm{al}}_n)$ for all $\eta' \in \mathrm{poss}(q, {<}n+1)$. Each $\eta' \in \mathrm{poss}(q, {<}n+1)$ depends on some $\eta^t_k \in \mathrm{poss}(p, {<}L^{\mathrm{al}}_n)$ and $s \in q(\alpha_n, t, L^{\mathrm{al}}_n)$, where $t = \hat\eta'(\alpha_n, n) = \hat\eta^t_k(\alpha_n, n)$. Let $t' := \hat\eta^t_k(i, n)$. If $t' <_{\mathrm{lex}} t$, then $p \wedge \eta \Vdash \dot r(L^{\mathrm{al}}_n) = F_{\eta^t_k}(s) = j^t_k \in s^t_k \subseteq \dot\varphi(L^{\mathrm{al}}_n)$; otherwise, if $t <_{\mathrm{lex}} t'$, then we have that $p \wedge \eta \Vdash \dot r(L^{\mathrm{al}}_n) \in F_{\eta^t_k}[\mathbf{c}^t_k] = s^t_k \subseteq \dot\varphi(L^{\mathrm{al}}_n)$.

This finishes the construction of $q$ and $\dot\varphi$. It is clear that they are as required. □

As an immediate consequence, we conclude:

**Corollary 8.13.** *Under* CH, *for* $\mathrm{t} \in \{\mathrm{lc}, \mathrm{al}\}$ *and* $i \in S^{\mathrm{pr}}$, *if* $\kappa = \max\{|S^{\mathrm{t}}_i|, \aleph_0\}^{\aleph_0}$, *then* $\mathbb{Q}$ *forces* $\mathfrak{c}^{\forall}_{a^{\mathrm{t}}_i, d^{\mathrm{t}}_i} \le \kappa$. *In the case that* $\mathrm{t} = \mathrm{al}$, *we also have* $\mathfrak{c}^{\forall}_{a^{\mathrm{t}}_i, h^{\mathrm{t}}_i} \le \kappa$.

*Proof.* Let $E_{\mathrm{t}} := \{i\} \cup S^{\mathrm{t}}_i$. By Lemma 7.5, $|\dot{R}(E_{\mathrm{t}})| \leq \kappa$. So it is enough to show that, for each $\mathbb{Q}$-name $\dot{r}$ for a real in $\prod a^{\mathrm{t}}_i$ and each $p \in \mathbb{Q}$, there is some (nice) $\mathbb{Q}$-name $\dot{\varphi}$ for an element of $\mathcal{S}(\dot{a}^{\mathrm{t}}_i, \dot{d}^{\mathrm{t}}_i)$ and some $q \leq p$ continuously reading $\dot{\varphi}$ only using indices in $E_{\mathrm{t}}$ such that $q \Vdash \dot{r} \in^* \dot{\varphi}$.

We define a $\mathbb{Q}$-name $\dot{r}'$ for an element of $\prod_{L\in\Omega} a_i(L)$ by

$$\dot{r}'(L) := \begin{cases} \dot{r}(n) & \text{if } L = L^{\mathrm{t}}_n, \\ 0 & \text{otherwise.} \end{cases}$$

By Lemma 8.12, there is a $q \leq p$ and a $\mathbb{Q}$-name $\dot{\varphi}'$ for a slalom in $\mathcal{S}(\dot{a}_i, \dot{d}_i)$ satisfying (i)–(iii) (where $\dot{\varphi}(L^{\mathrm{t}'}_n)$ with $\mathrm{t}' \in \{\mathrm{lc}, \mathrm{al}\} \smallsetminus \{\mathrm{t}\}$ can be taken as the canonical name for $\{0\}$). Let $\dot{\varphi}(n)$ be a (nice) name for $\dot{\varphi}'(L^{\mathrm{lc}}_n)$ for each $n < \omega$. It is clear that $q \Vdash$"$\dot{\varphi} \in \mathcal{S}(\dot{a}^{\mathrm{lc}}_i, \dot{d}^{\mathrm{al}}_i)$ and $\dot{r} \in^* \dot{\varphi}$" and that $q$ continuously reads $\dot{\varphi}$ only using indices in $E_{\mathrm{t}}$.

The additional claim for $\mathrm{t} = \mathrm{al}$ follows by Lemma 2.10 (2), because $\Vdash d^{\mathrm{al}}_i \leq^* h^{\mathrm{al}}_i$. □

Using Assumption 8.1, the proof of Theorem 1.4 is concluded by Corollary 8.3, Corollary 8.5 and Corollary 8.13.

## 9. Discussion

The forcing construction developed in Section 4–7 can be easily generalised as in [FGKS17], in the sense that we can replace the subatomic families $\mathrm{K}_L$ and $\mathrm{K}^t_\ell$ by any arbitrary subatomic family (even allowing different subatomic families at each $(i, L)$ and $(\alpha, t, \ell)$). However, we decided to proceed with the construction presented in this paper to ease the presentation, particularly in relation to the parameters of the forcing. In the generalised framework, to make sense of results like the existence of conditions with arbitrary support (Corollary 5.11), continuous reading of names (Theorem 7.2) and rapid reading (Theorem 7.8), we would need a more complicated setting for the parameters, which would obscure the presentation. Our forcing construction can be adapted to different situations, but we believe it would be simpler to modify the current construction to other contexts than to follow an utterly obscure abstract construction for the most general case possible.

This forcing construction, as well as those in [GS93, Kel08, KS09, KS12, FGKS17, GK21], satisfies enough rapid reading to guarantee that random reals are never added (see Lemma 7.10). It would be interesting to discover how to modify these constructions to allow for a part adding random reals while ensuring rapid reading outside the random part and force the same results with larger values of $\mathrm{cov}(\mathcal{N})$.

Concerning Cichoń's diagram, our construction forces $\mathrm{cov}(\mathcal{N}) = \mathfrak{d} = \aleph_1$ and the other cardinals equal to the continuum because, as indicated in Section 1, the cardinal $\mathfrak{v}^{\forall}_{b,h}$ is below $\mathrm{non}(\mathcal{E})$ (under some conditions fulfilled by our parameters), and $\mathrm{non}(\mathcal{E})$ is below $\mathrm{non}(\mathcal{M})$ and $\mathrm{non}(\mathcal{N})$.

Recall the following:

**Theorem 9.1** ([FGKS17]). *Under* CH, *if* $\lambda_1 \leq \lambda_3 \leq \lambda_4$ *and* $\lambda_2 \leq \lambda_3$ *are infinite cardinals such that* $\lambda_i^{\aleph_0} = \lambda_i$ *for* $i \in \{1,2,3,4\}$, *then there is some proper* $\omega^\omega$*-bounding poset with* $\aleph_2$*-cc forcing* $\mathrm{cov}(\mathcal{N}) = \mathfrak{d} = \aleph_1$, $\mathrm{non}(\mathcal{M}) = \mathrm{cof}(\mathcal{M}) = \lambda_1$, $\mathrm{non}(\mathcal{N}) = \lambda_2$, $\mathrm{cof}(\mathcal{N}) = \lambda_3$ *and* $\mathfrak{c} = \lambda_4$ *(see Figure 6).*

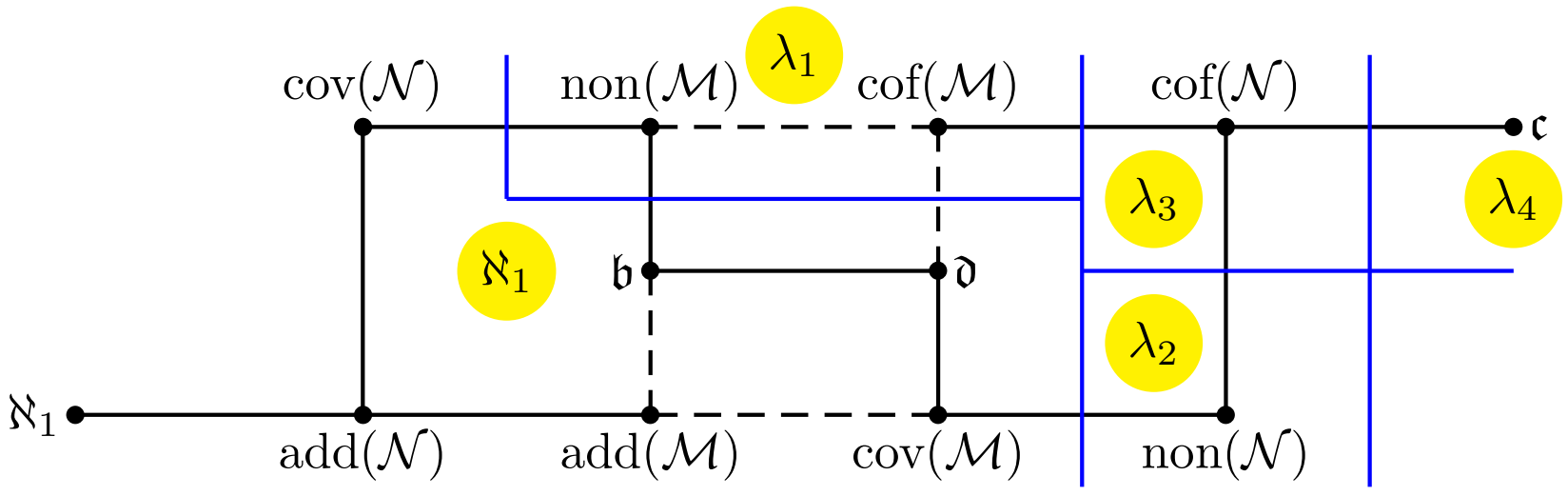

FIGURE 6. The constellation of Cichoń's diagram forced in [FGKS17, GK21].

The construction in [GK21] forces the same and, in addition, forces $\aleph_1$ many different localisation cardinals $\mathfrak{c}^\forall_{b,h}$ (with parameters in the ground model). It seems likely that the methods of [GK21] to force the constellation in Figure 6 and the forcing construction in this paper can be combined to force the same constellation along with $\min\{\lambda_1, \lambda_2\}$ many different cardinals of each one of the six types discussed in this work. A bit more work would be needed to force the constellation in Figure 6 along with $\min\{\lambda_1, \lambda_2\}$ many different cardinals of type $\mathfrak{v}^\forall_{b,h}$, $\lambda_1$ many of both types $\mathfrak{c}^\exists_{b,h}$ and $\mathrm{cov}(\mathcal{I}_f)$, $\lambda_2$ many of both types $\mathfrak{v}^\exists_{b,h}$ and $\mathrm{non}(\mathcal{I}_f)$, and $\lambda_3$ many of type $\mathfrak{c}^\forall_{b,h}$. This result would be quite optimal, because (under certain conditions) $\mathfrak{v}^\forall_{b,h} \leq \mathrm{non}(\mathcal{E})$, $\mathfrak{c}^\exists_{b,h} \leq \mathrm{non}(\mathcal{M})$, $\mathfrak{v}^\exists_{b,h} \leq \mathrm{non}(\mathcal{N})$ and $\mathfrak{c}^\forall_{b,h} \leq \mathrm{cof}(\mathcal{N})$ (as discussed in Section 1).

This work solves many open questions from [KM22], but the questions about additivities and cofinalities of Yorioka ideals remain open:

**Question 9.2.** *Is it consistent with* ZFC *that there are two different cardinals of the form* $\mathrm{add}(\mathcal{I}_f)$, *or two different cardinals of the form* $\mathrm{cof}(\mathcal{I}_f)$*?*

Forcing two different additivities or cofinalities would also indicate a method to separate infinitely (and even continuum) many of them. While the method employed in this paper could be used for the cofinalities, it does not work for the additivities, since $\mathrm{add}(\mathcal{I}_f) \leq \mathfrak{b}$ ([KO08], see Figure 1). For the same reason, continuum many different additivities of Yorioka ideals imply that $\mathfrak{b} = \mathfrak{c}$ is weakly inaccessible: If there are continuum many different cardinals below $\mathfrak{b}$, then $\mathfrak{b} = \mathfrak{c}$. Moreover, since there are continuum many different cardinals below $\mathfrak{c}$, it follows that $\mathfrak{c} = \aleph_\mathfrak{c}$. Consequently, $\mathfrak{c}$ is a limit cardinal and, notably, it is weakly inaccessible (while $\mathfrak{b}$ is regular).

Institute of Discrete Mathematics and Geometry, TU Wien, Wiedner Hauptstrasse 8–10/104, 1040 Wien, Austria

*Email address*: miguel.cardona@mail.huji.ac.il

*URL*: https://sites.google.com/mail.huji.ac.il/miguel-cardona-montoya/home-page

Institute of IT Security Research and Center for Artificial Intelligence, St. Pölten University of Applied Sciences, Campus-Platz 1, 3100 St. Pölten, Austria

*Email address*: mail@l17r.eu

*URL*: https://l17r.eu

Creative Science Course (Mathematics), Faculty of Science, Shizuoka University, 836 Ohya, Suruga-ku, Shizuoka City, Shizuoka Prefecture, 422-8529, Japan

*Email address*: diego.mejia@shizuoka.ac.jp

*URL*: http://www.researchgate.net/profile/Diego-Mejia-14